\documentclass[a4paper,10pt]{article}

\usepackage[english]{babel}
\usepackage{ae}
\usepackage{graphicx}
\usepackage{amsmath, amssymb, amsfonts, mathtools, amsthm}
\usepackage{url}
\usepackage{authblk}

\newtheorem{definition}{Definition}

\newcommand{\argmin}{\operatornamewithlimits{arg\,min}}

\def\vec#1{\mathchoice{\mbox{\boldmath$\displaystyle#1$}}
{\mbox{\boldmath$\textstyle#1$}}
{\mbox{\boldmath$\scriptstyle#1$}}
{\mbox{\boldmath$\scriptscriptstyle#1$}}}

\begin{document}

\title{Decomposition of Optical Flow on the Sphere}

\author[1]{Clemens Kirisits}
\author[1]{Lukas F. Lang}
\author[1,2]{Otmar Scherzer}
\affil[1]{\footnotesize Computational Science Center, University of Vienna, Oskar-Morgenstern-Platz\ 1, 1090 Vienna, Austria}
\affil[2]{Radon Institute of Computational and Applied Mathematics, Austrian Academy of Sciences, Altenberger Str.\ 69, 4040 Linz, Austria}

\maketitle

\begin{abstract}
\noindent
We propose a number of variational regularisation methods for the estimation and decomposition of motion fields on the $2$-sphere. While motion estimation is based on the optical flow equation, the presented decomposition models are motivated by recent trends in image analysis. In particular we treat $u+v$ decomposition as well as hierarchical decomposition. Helmholtz decomposition of motion fields is obtained as a natural by-product of the chosen numerical method based on vector spherical harmonics. All models are tested on time-lapse microscopy data depicting fluorescently labelled endodermal cells of a zebrafish embryo.
\end{abstract}

\section{Introduction}
Motion estimation is a fundamental task for the analysis of spatiotemporal data, the prototypical example of which are sequences of images taken by a camera. The term \emph{optical flow} has been coined to designate the apparent motion in such data. Its accurate and efficient estimation has been a major topic in the fields of computer vision and image processing for more than 30 years. However, the applicability of optical flow algorithms is by no means limited to flat two-dimensional projections of real world scenes. The advance of microscopy techniques has led to a particularly promising application of optical flow: cell motion analysis. Reliable optical flow algorithms supplied with microscopy images of sufficiently high spatial and temporal resolution can obviously help understanding cellular dynamics in transparent organisms, see for example \cite{AmaMyeKel13,MelCamLomRizVer07,QueMenCam10,SchmShaScheWebThi13}.

\begin{figure}
	\centering
	\includegraphics[width=0.49\textwidth]{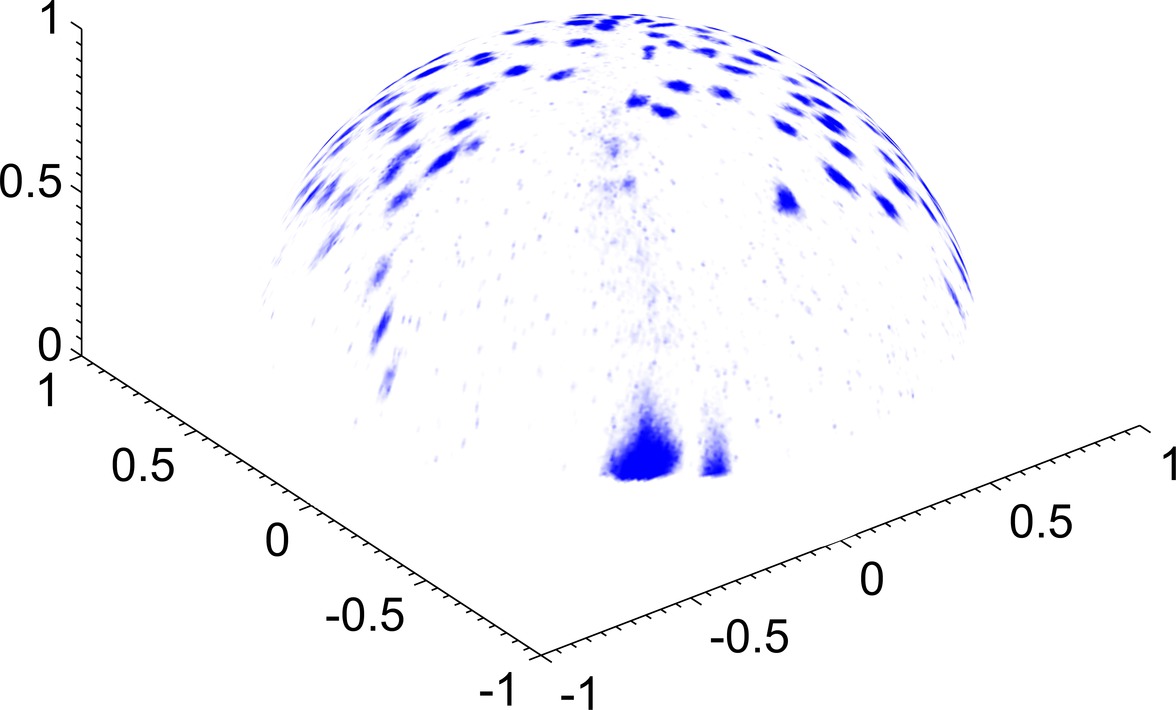} \hfill
	\includegraphics[width=0.49\textwidth]{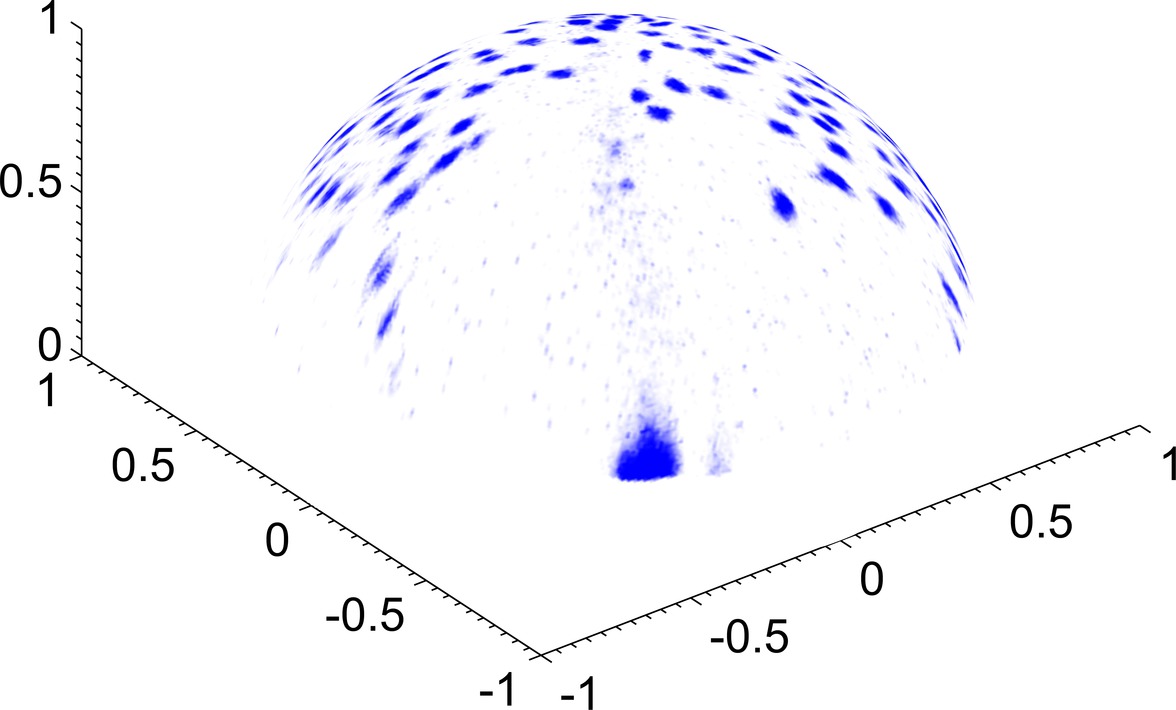} \\ \medskip
	\includegraphics[width=0.49\textwidth]{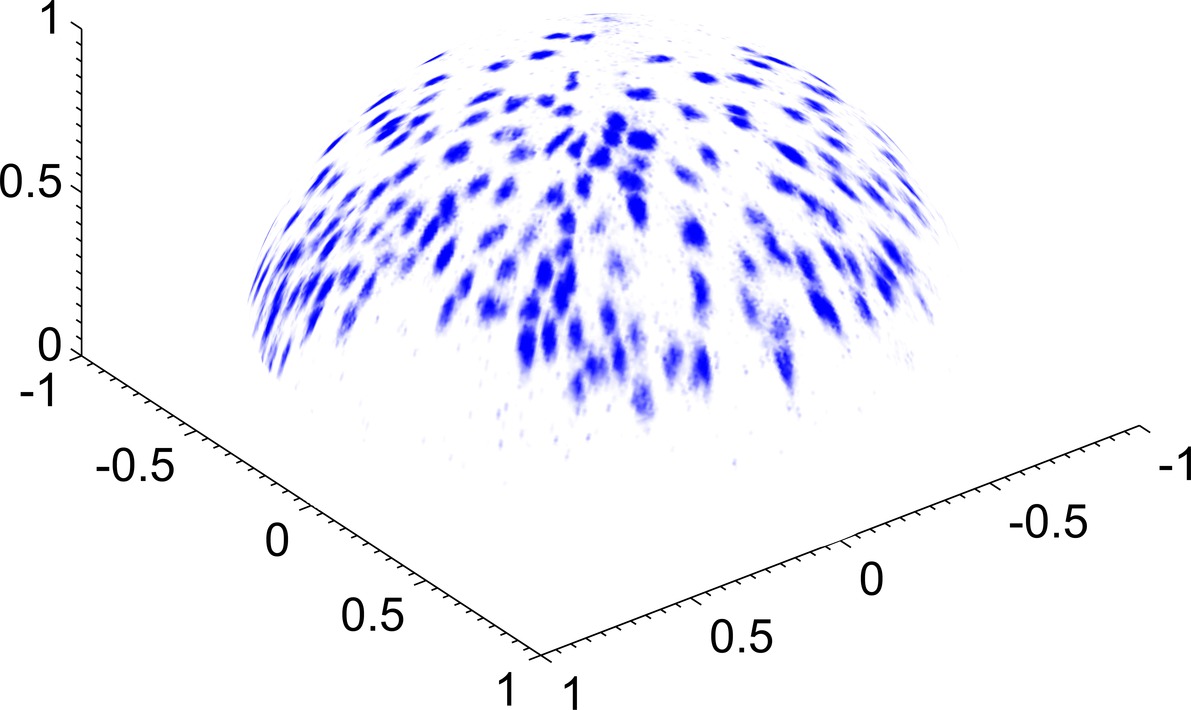} \hfill
	\includegraphics[width=0.49\textwidth]{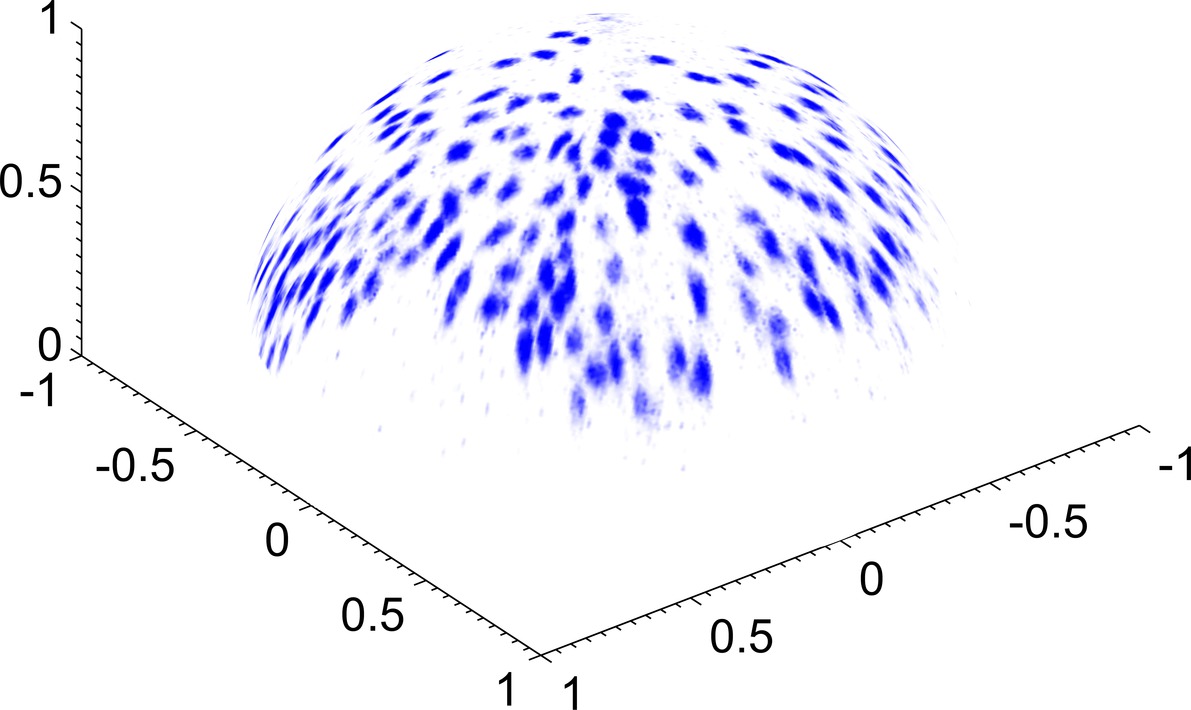} \\
	\caption{Frames no. 57 (left) and 58 (right) of the zebrafish microscopy image sequence. Blue colour indicates fluorescence response. The spherical image is obtained by a radial projection of the unfiltered data onto a fitted sphere. Top and bottom row differ by a rotation of $180$ degrees around the $z$-axis. See Sec.~\ref{sec:experiments} for more details on the data and preprocessing.}
	\label{fig:data3d}
\end{figure}

The particular dataset we are working with in this article depicts a living zebrafish embryo during early embryogenesis. Main feature of this dataset are the embryo's endodermal cells, which have been labelled with a fluorescent protein and are known to develop on the surface of the zebrafish's spherical yolk, see Fig.~\ref{fig:data3d} and Sec.~\ref{sec:data}. The distribution of these cells can be modelled by a nonnegative function $F$ depending on time $t$ and position $x$ on the $2$-sphere, such that the number $F(t,x)$ is directly proportional to the fluorescence response of a point $x$ at time $t$. The models we propose below are motivated, but surely not restricted, to this specific type of data. We argue that the problem of extracting and analysing motion from spherical data is sufficiently general so as to be of potential interest to a wider audience.

Our motion models are based on the optical flow equation
\begin{equation}\label{eq:ofc1}
	\nabla_{\mathcal{S}} F(t,x) \cdot u(t,x) + \partial_t F (t,x) = 0,
\end{equation}
which every vector field $u$ describing the temporal evolution of $F$ should (approximately) satisfy for all $x$ and $t$. Here $\nabla_{\mathcal{S}}$ denotes the surface gradient on the sphere. We derive this equation in more detail in Sec.~\ref{sec:of}. Directly solving the optical flow equation for $u$ is infeasible. We therefore use Tikhonov regularisation to compute an approximate solution to \eqref{eq:ofc1}. Tikhonov regularisation consists in minimising a functional which is a weighted sum of two terms. The first one, usually called data or similarity term, is the squared $L^2$ norm of the left hand side of \eqref{eq:ofc1}. The second term is a regularising functional $\mathcal{R}(u)$, which in this article will always be a Sobolev $H^s$ norm (Sec.~\ref{sec:reg}). These norms are introduced in Sec.~\ref{sec:sobolev}.

Next, we extend the motion estimation method outlined above to two types of decomposition models that are also variational in nature (Sec.~\ref{sec:decomposition}). While the input of those is again $F$, their outputs are now two or more vector fields capturing different structural parts of the total motion. Both models are adaptations of recently proposed image decomposition techniques to the optical flow setting. The first one is a $u+v$ decomposition. Its idea is to replace $u$ in the data term with a sum $u+v$, and then to add two different regularising functionals, one for $u$ and one for $v$. The second one is a hierarchical model. Roughly speaking, one repeatedly minimises an optical flow functional, while the amount of regularisation is constantly decreased. In every iteration, the $u$ in the data term is replaced by a sum $u + \sum_i u_i$, where only $u$ is optimised and the $u_i$ are the results from previous iterations.

Finally, all optimisation problems are solved by projecting them onto finite dimensional spaces spanned by vector spherical harmonics (Sec.~\ref{sec:numerics}). One advantage of this method is that it automatically yields Helmholtz decompositions of all computed vector fields. Mathematical background on (vector) spherical harmonics is presented in Secs.~\ref{sec:sh} and \ref{sec:vsh}. In Sec.~\ref{sec:experiments} we provide details of our implementation, give a more detailed account of the used microscopy data, and show experimental results with these data.

To summarise, the main novelties presented in this article are decomposition models for optical flow on the sphere together with their application to microscopy data.

\subsection*{Related Work}
The first variational optical flow method is usually attributed to \cite{HorSchu81}, where they used an $H^1$ seminorm for regularisation. In \cite{Schn91a} it was shown that this particular choice leads to a well-posed problem. We refer to \cite{AubKor06} for a gentle introduction to optical flow, to \cite{WeiBruBroPap06} for an overview of different optical flow functionals and to \cite{BakSchaLewRotBla11} for a recent survey and benchmark.

Optical flow algorithms have only recently been extended to data defined on non-Euclidean domains. In \cite{ImiSugTorMoc05,TorImiSugMoc05} images defined on the sphere were treated, whereas in \cite{LefBai08} the original functional by Horn and Schunck was generalised to $2$-Riemannian manifolds and well-posedness was verified. Most recently, optical flow on evolving manifolds has been considered in \cite{KirLanSch13,KirLanSch13a_report}.

Horn and Schunck \cite{HorSchu81} numerically solved the variational problem by applying a finite difference scheme to the Euler-Lagrange equations. A similar approach was adopted in \cite{KirLanSch13,KirLanSch13a_report} after parametrisation of the surface. In \cite{LefBai08}, however, the problem was solved by finite elements on a surface triangulation. Finally, we mention the work by Schuster and Weickert \cite{SchuWei07}, where they used projection methods to solve the optical flow equation in the plane. Instead of Tikhonov-regularising their solution, they solely relied on regularisation by discretisation. The main reason for choosing a numerical method based on vector spherical harmonics in this article is that $H^s$-type regularisers, for arbitrary real $s$, are handled very easily in contrast to most other methods.

The aim of $u+v$ image decomposition models, as pioneered in \cite{Mey01}, is to separate the cartoon and texture parts of images. While the cartoon component should capture large-scale structural components and should therefore be piecewise smooth, the texture component is supposed to consist of high-frequency oscillating patterns. Since the original model was promising but hard to implement, a large number of modifications and approximations have been proposed. In some of them the problematic $g$-norm was approximated by an $H^{-1}$ norm \cite{OshSolVes03,VesOsh03}. Recently, $u+v$ models have been extended to the $\mathbb{R}^2$ optical flow setting \cite{AbhBelSch09}. Hierarchical models, originally introduced in \cite{TadNezVes04} for image analysis, have not yet been tried in combination with optical flow. They have, however, the preferable property of producing arbitrarily fine multiscale descriptions of input data. As a concluding remark about vector field decompositions, let us remark that Helmholtz-Hodge decompositions of motion fields have enjoyed a certain degree of attention in recent years, not only in the plane \cite{KohMemSchn03,YuaSchSte07,YuaSchSte09}, but also on surfaces \cite{KhaLefAmmBai11}.

Applying optical flow algorithms to cell microscopy data has become increasingly popular lately. See for example \cite{AmaMyeKel13,KirLanSch13,KirLanSch13a_report,MelCamLomRizVer07,SchmShaScheWebThi13} and the references therein. We highlight the article \cite{SchmShaScheWebThi13}, where also endodermal cells of zebrafish embryos have been analysed. There the authors point out that, although of immense importance for developmental biology, only little is known about the motion behaviour of this type of cells.

\section{Notation and Background}
\subsection{Scalar Spherical Harmonics}\label{sec:sh}
Let
\begin{equation*}
	\mathcal{S} = \{ x \in \mathbb{R}^3 : |x| = 1 \}
\end{equation*}
be the two-sphere embedded in $\mathbb{R}^3$. For functions $F : \mathcal{S} \to \mathbb{R}$ we define the Laplace-Beltrami operator by 
\begin{equation*}
	\Delta_\mathcal{S} F = - \Delta \bar F,
\end{equation*}
where $\Delta$ is the usual Laplacian of $\mathbb{R}^3$ and $\bar F (x) = F(x/|x|)$ is the radially constant extension of $F$ to $\mathbb{R}^3 \setminus \{0\}$. The eigenvalues of $\Delta_\mathcal{S}$ are
\begin{equation}\label{eq:ev}
	\lambda_n = n(n+1), \quad n \in \mathbb{N}_0.
\end{equation}
The corresponding eigenspaces $\mathrm{Harm}_n$ have dimension $2n+1$ and are mutually orthogonal in $L^2(\mathcal{S})$. Their direct sum equals $L^2(\mathcal{S})$. Every eigenfunction $Y_n \in \mathrm{Harm}_n$ lies in $C^\infty(\mathcal{S})$ and is called \emph{(scalar) spherical harmonic of degree} $n$. Their name derives from the equivalent characterisation of $\mathrm{Harm}_n$ as the restriction to $\mathcal{S}$ of the space of harmonic polynomials $P:\mathbb{R}^3 \to \mathbb{R}$ that are homogeneous of degree $n$. From now on
\begin{equation} \label{eq:scalarbasis}
	\{Y_{nj} : n \in \mathbb{N}_0,\, 1 \le j \le 2n+1\}
\end{equation}
always refers to a particular orthonormal basis of $L^2(\mathcal{S})$ consisting of real-valued spherical harmonics. For numerical experiments we use the so-called fully normalised spherical harmonics, see \cite[Sec.~5.2]{Mic13} for a detailed construction.

\subsection{Vector Spherical Harmonics}\label{sec:vsh}
Denote by $\nu$ the outward unit normal to $\mathcal{S}$ and let
\begin{equation*}
	\nabla_{\mathcal{S}} F = \nabla \bar F
\end{equation*}
be the surface gradient of $F : \mathcal{S} \to \mathbb{R}$.
\begin{definition} \label{def:vsh}
	Let $n \in \mathbb{N}_0$ and $Y_n \in \mathrm{Harm}_n$. Whenever a function $y : \mathcal{S} \to \mathbb{R}^3$ that does not vanish identically admits one of the following three representations
	\begin{equation}\label{eq:vspharm}
	y =
	\begin{cases}
		y_{n}^{(1)} \coloneqq	Y_{n} \nu, \\
		y_{n}^{(2)}  \coloneqq  \nabla_{\mathcal{S}} Y_{n}, \\
		y_{n}^{(3)}  \coloneqq  \nabla_{\mathcal{S}} Y_{n} \times \nu,
	\end{cases}
\end{equation}
then $y = y_{n}^{(i)}$ is called a \emph{vector spherical harmonic of degree $n$ and type $i$}. For obvious reasons we refer to types $2$ and $3$ as \emph{tangential vector spherical harmonics}.
\end{definition}
Note that there is no tangential spherical harmonic of degree $0$.

\begin{figure}
\centering
	\includegraphics[width=0.49\textwidth]{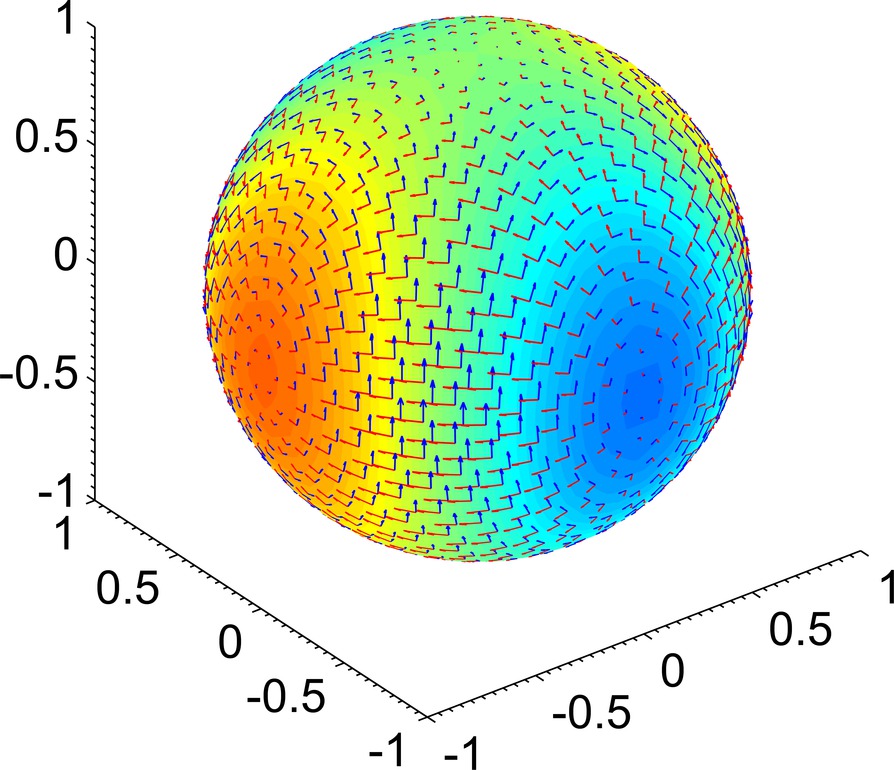} \hfill
	\includegraphics[width=0.49\textwidth]{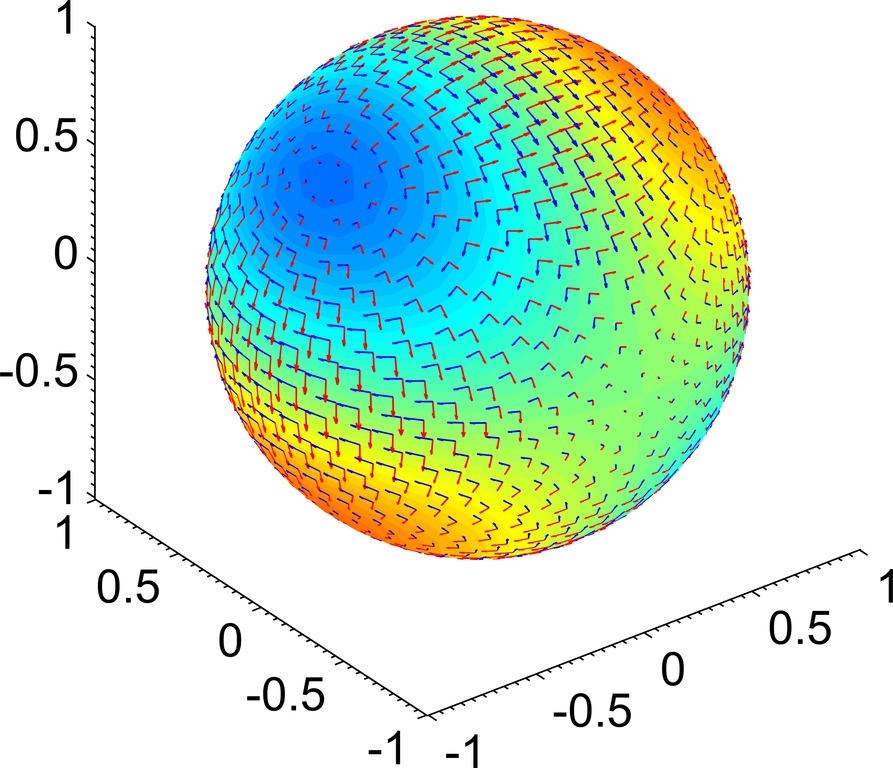} \\ \medskip
	\includegraphics[width=0.49\textwidth]{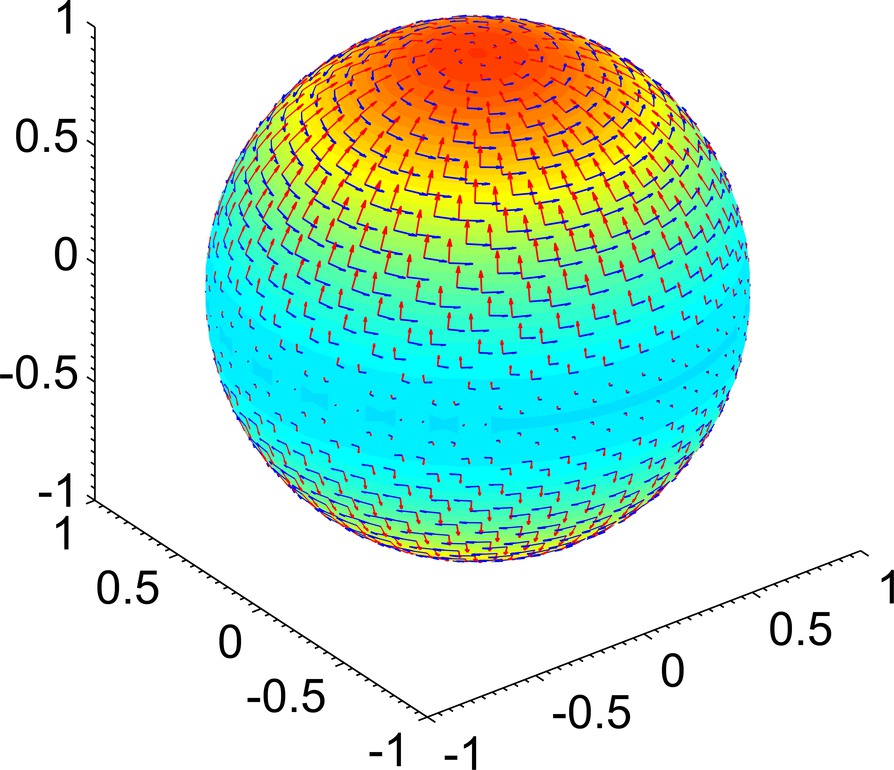} \hfill
	\includegraphics[width=0.49\textwidth]{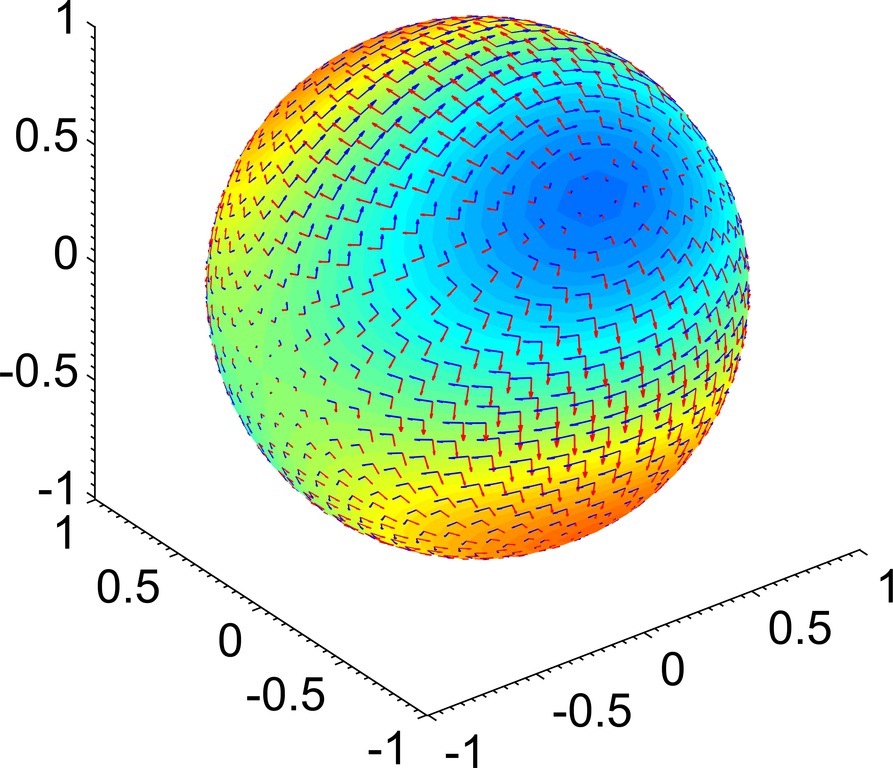} \\ \medskip
	\includegraphics[width=0.65\textwidth]{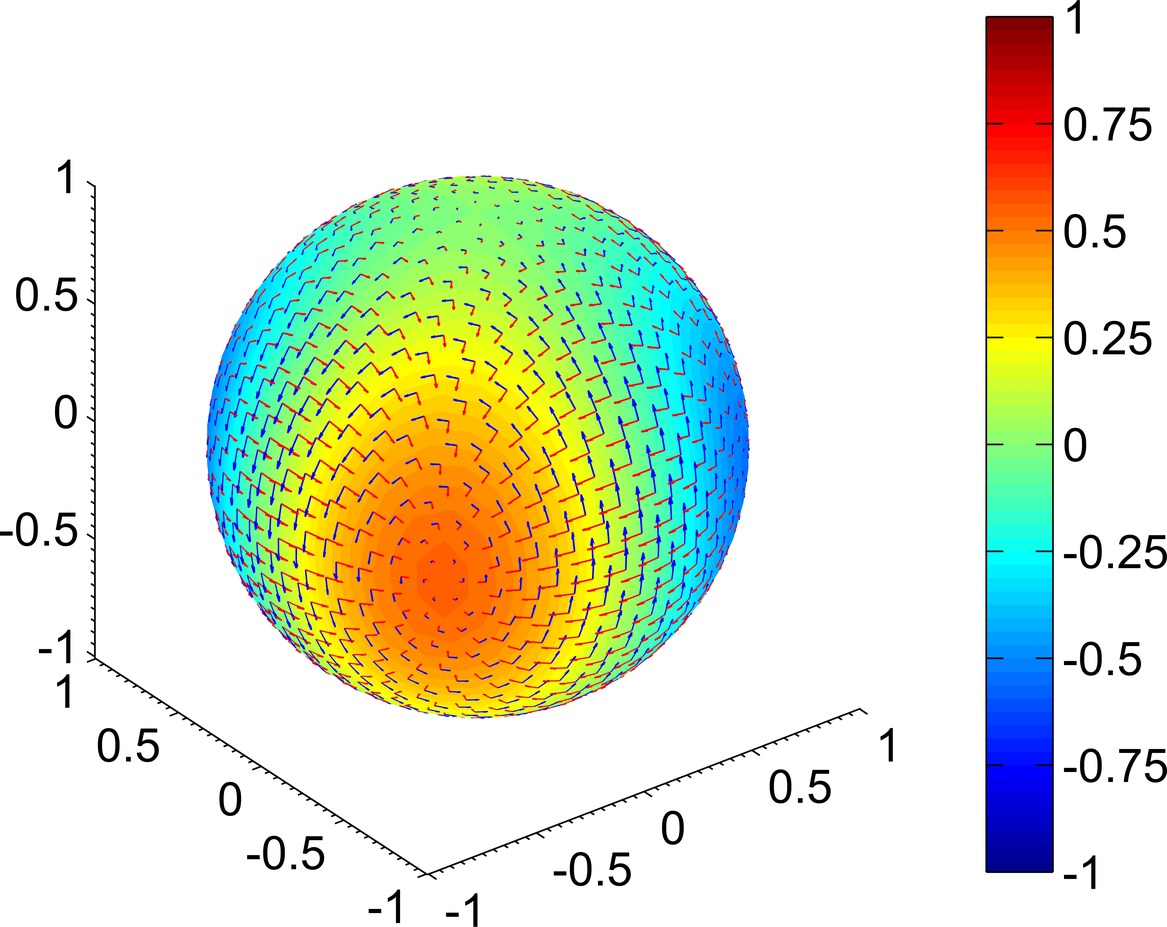}
	\caption{Fully normalised scalar and corresponding vector spherical harmonics of degree $n = 2$. Scalar spherical harmonics are depicted using a blue to red colour bar. Type $2$ vector spherical harmonics are visualised with red arrows, type $3$ ones with blue arrows. Note that the length of the vectors has been scaled for better illustration.}
	\label{fig:vspharm}
\end{figure}

We are mainly interested in the space $L^2(\mathcal{S}, T\mathcal{S})$ of square integrable tangent vector fields on $\mathcal{S}$ endowed with the inner product
\begin{equation*}
	\langle u,v \rangle = \int_\mathcal{S} u \cdot v \, \mathrm{d}\mathcal{S},
\end{equation*}
where $\mathrm{d}\mathcal{S}$ is the usual surface measure on the sphere. An orthonormal basis of this space is obtained from \eqref{eq:scalarbasis} by setting
\begin{equation}\label{eq:basis}
	\begin{aligned}
		y_{nj}^{(2)}	&=	\lambda_n^{-1/2} \nabla_{\mathcal{S}} Y_{nj}, \\
		y_{nj}^{(3)}	&=	\lambda_n^{-1/2} \nabla_{\mathcal{S}} Y_{nj} \times \nu
	\end{aligned}
\end{equation}
for all $n\in \mathbb{N}$ and $1 \le j \le 2n+1$. In Fig.~\ref{fig:vspharm} a handful of the elements of both bases \eqref{eq:scalarbasis} and \eqref{eq:basis} are depicted. Every $v \in L^2(\mathcal{S}, T\mathcal{S})$ has the following Fourier series representation
\begin{gather*}
	 \sum^3_{i=2} \sum^\infty_{n=1} \sum^{2n+1}_{j=1} \langle v,y_{nj}^{(i)} \rangle y_{nj}^{(i)}.
\end{gather*}
In particular, we have Parseval's identity
\begin{gather*}
	\| v \|_{L^2(\mathcal{S}, T \mathcal{S})}^2	= \sum_{i,n,j} \langle v, y_{nj}^{(i)} \rangle^2.
\end{gather*}
For a comprehensive and unified treatment of both scalar and vector spherical harmonics we refer to \cite{FreSchr09}.

\subsection{Sobolev Spaces on the Sphere} \label{sec:sobolev}
For an arbitrary real number $s$, the space $H^s(\mathcal{S})$ is commonly defined as the domain of $\Delta_{\mathcal{S}}^{s/2}$. See \cite[p.~37]{LioMag72a} or \cite[Sec.~6.2]{Mic13} for example. In this section we introduce the spaces $H^s(\mathcal{S}, T\mathcal{S})$ by means of the vectorial counterpart of $\Delta_{\mathcal{S}}$.

For tangent vector fields $v$ we define the Laplace-Beltrami operator by
\begin{equation*}
	\vec{\Delta}_{\mathcal{S}} v = \mathrm{P} \Delta_{\mathcal{S}} v,
\end{equation*}
where application of $\Delta_{\mathcal{S}}$ to $v$ is understood componentwise and $\mathrm{P} = \mathrm{P}(x)$ is the orthogonal projector onto the tangent plane $T_x \mathcal{S}$, compare \cite[Def.~5.26]{FreSchr09}. The tangential vector spherical harmonics introduced in Def.~\ref{def:vsh} are eigenfunctions of this operator to the same eigenvalues as their scalar counterparts: If we let
\begin{equation*}
	\mathrm{harm}_n = \mathrm{span} \Bigl\{ y_{nj}^{(i)} : 1 \le j \le 2n+1, \, i=2,3 \Bigr\},
\end{equation*}
then
\begin{equation*}
	\vec{\Delta}_{\mathcal{S}} y_n = \lambda_n y_n
\end{equation*}
for every $y_n \in \mathrm{harm}_n$. The $\lambda_n$ are as defined in \eqref{eq:ev}, the only difference being that now $n > 0$ and therefore the spectrum of $\vec{\Delta}_{\mathcal{S}}$ is strictly positive. Applying functional calculus, we formally define the $s$-th power of $\vec{\Delta}_{\mathcal{S}}$ by
\begin{equation*}
	\vec{\Delta}_{\mathcal{S}}^s v  = \sum_{i,n,j} \lambda_n ^s \langle v, y_{nj}^{(i)} \rangle y_{nj}^{(i)}.
\end{equation*}
Finally, for every $s\in\mathbb{R}$, set
\begin{equation} \label{eq:norm}
	\| v \|_{H^s(\mathcal{S}, T\mathcal{S})}^2 \coloneqq \| \vec{\Delta}_{\mathcal{S}}^{s/2} v \|_{L^2(\mathcal{S}, T\mathcal{S})}^2 = \sum_{i,n,j} \lambda_n^s \langle v, y_{nj}^{(i)} \rangle^2.
\end{equation}
Note that, in contrast to the scalar setting, this functional is an actual norm, rather than only a seminorm. We therefore define, for every real $s$, $H^s(\mathcal{S}, T\mathcal{S})$ as the space of all distributions $v \in C^\infty (\mathcal{S}, T\mathcal{S})'$ for which the series in \eqref{eq:norm} is finite. Clearly, if $(\mu_n)_n$ is any sequence satisfying
\begin{equation}\label{eq:comparable}
	c \mu_n \le \lambda_n^s \le C \mu_n
\end{equation}
for two positive constants $c$, $C$ and for all $n$, then replacing $\lambda_n^s$ with $\mu_n$ in \eqref{eq:norm} leads to an equivalent norm and thus to the same space. For every sequence of positive numbers $\mu_n$ we denote the resulting norm simply by $\|\cdot\|_{\mu_n}$.

\section{Decomposition Models for Optical Flow}
\subsection{Optical Flow on the Sphere} \label{sec:of}
Let $I = [0, T] \subset \mathbb{R}$ be a time interval. We assume to be given a scalar time-varying (brightness) function
\begin{equation*}
	F: I \times \mathcal{S} \rightarrow \mathbb{R}.
\end{equation*}
The problem of estimating optical flow consists in tracking the temporal evolution of the data $F$ by means of a time-dependent vector field. Our optical flow model is based on the so-called \emph{brightness constancy assumption}: We assume existence of a function $\phi : I \times \mathcal{S} \rightarrow \mathcal{S}$ satisfying
\begin{align} 
	F(t, \phi(t, x))	&= F(0, x), \label{eq:bca} \\
	\phi(0,x)			&= x, \nonumber
\end{align}
for all $x$ and $t$. Intuitively this means that for every starting point $x$ on the sphere, the function $F$ remains constant along the trajectory $t \mapsto \phi(t, x)$. In addition we require that $\phi(t,\cdot)$ is a diffeomorphism of $\mathcal{S}$ for every $t$. The first equation in~\eqref{eq:bca} implies that
\begin{equation*}
	\frac{\mathrm{d}}{\mathrm{d}t} F(t, \phi(t, x)) = \nabla_{\mathcal{S}} F (t, \phi(t, x)) \cdot \partial_t \phi(t,x) + \partial_{t} F(t, \phi(t, x)) = 0.
\end{equation*}
This equation is typically written in terms of the vector field $u: I \times \mathcal{S} \rightarrow T\mathcal{S}$ whose integral curves are the trajectories $t \mapsto \phi(t, x)$, which is defined by the equation $u(t, \phi(t, x)) = \partial_t \phi(t,x)$. The resulting \emph{optical flow equation} reads
\begin{equation}\label{eq:ofc}
	\nabla_{\mathcal{S}} F \cdot u + \partial_{t} F = 0.
\end{equation}

\subsection{Regularisation}\label{sec:reg}
Solving the optical flow equation directly is problematic. In general, a solution to \eqref{eq:ofc} need not exist, and if it exists, it cannot be unique. The typical remedy is Tikhonov regularisation, where one minimises a functional of the form
\begin{equation}\label{eq:reg}
	\| \nabla_{\mathcal{S}} F \cdot u + \partial_t F \|^2_{L^2(I\times S)} + \alpha \mathcal{R}(u).
\end{equation}
with $\mathcal{R}$ being a regularising functional that incorporates a-priori knowledge about desirable solutions. The parameter $\alpha > 0$ controls the amount of regularisation. In the context of optical flow one usually tries to enforce spatial (and temporal) smoothness on the solution. A natural candidate for $\mathcal{R}$ would then be the squared Sobolev $H^1(I\times\mathcal{S}, T\mathcal{S})$ (semi-)norm, which penalises first derivatives in space and time equally, compare \cite{KirLanSch13a_report,WeiSchn01b}

Another, and in fact more popular, possibility is to drop time regularisation, in which case minimisation of \eqref{eq:reg} is equivalent to minimising
\begin{equation}\label{eq:regh1}
	\| \nabla_{\mathcal{S}} F \cdot u + \partial_t F \|^2_{L^2(S)} + \alpha \|u\|^2_{H^1(\mathcal{S}, T\mathcal{S})}
\end{equation}
for each instant $t$ separately. This corresponds to the original approach of Horn and Schunck \cite{HorSchu81}. From now on we denote the above data term by $\mathcal{D}(u,F)$. Instead of \eqref{eq:regh1} we consider the more general class of optical flow functionals
\begin{equation}\label{eq:rega_n}
	\mathcal{E}_{\mu_n} (u) \coloneqq \mathcal{D}(u,F) + \| u \|^2_{\mu_n}.
\end{equation}
The regularisation parameter is omitted, as it can be considered a constant factor in the sequence $(\mu_n)_n$. Functional \eqref{eq:rega_n} forms the basic optical flow setting of this article. All variational models considered here are extensions of \eqref{eq:rega_n}.

\subsection{Optical Flow Decomposition}\label{sec:decomposition}
The following two decomposition models are inspired by techniques that have recently been developed in the context of image analysis. The fact that motion estimation based on \eqref{eq:ofc} can be viewed as denoising of vector-valued images suggests the translation of said image decomposition models to the optical flow setting  \cite{AbhBelSch09}.

\subsubsection*{$u+v$ Models}
The aim is now not to extract one, but two vector fields $u$ and $v$ in such a way that they capture different structural parts of the total motion $u + v$ of $F$. The idea is to solve the following variational problem
\begin{gather}
	\min_{u,v} \mathcal{E}_{\mu_n, \nu_n} (u,v) \nonumber
	\intertext{where the functional $\mathcal{E}_{\mu_n, \nu_n}$ is defined as}	
	\mathcal{E}_{\mu_n, \nu_n}(u,v) \coloneqq \mathcal{D}(u+v,F) + \| u \|^2_{\mu_n} + \| v \|^2_{\nu_n}. \label{eq:u+v}
\end{gather}
Choosing, for instance, the two regularisers to be $H^1(\mathcal{S},T\mathcal{S})$ and $H^{-1}(\mathcal{S},T\mathcal{S})$ norms, respectively, would lead to a model which, in spirit, comes closest to the image decomposition model considered in \cite{VesOsh03}. Generalising \eqref{eq:u+v} to a decomposition into $k\in\mathbb{N}$, instead of two, constituents is possible, but will not be considered here \cite{FruPonSch11}.

\subsubsection*{Hierarchical Models}
Hierarchical image decomposition models have been introduced in \cite{TadNezVes04}. There, an original image is decomposed by repeatedly applying denoising steps. The input of one such step is the residual of the previous one. In every step the degree of regularisation is decreased. In contrast to $u+v$ models, hierarchical decomposition models provide multiscale descriptions of the input data.

This iterative procedure can be transferred to the optical flow setting as follows. Let $\| \cdot \|^2_{\mu^{(k)}_n}$, $k \in \mathbb{N}$, be a sequence of norms as defined in Section \ref{sec:sobolev}, such that
\begin{equation} \label{eq:hiersequ}
	\mu^{(k+1)}_n \le \mu^{(k)}_n
\end{equation}
for all $n$ and $k$. For every such sequence of sequences we propose the following iterative scheme,
\begin{equation}\label{eq:hierachical}
	u_k =
	\begin{dcases}
		\argmin_u \mathcal{E}_{\mu^{(1)}_n} (u), & \text{if } k = 1, \\
		\argmin_u \mathcal{D} \big( u + \textstyle{\sum^{k-1}_{i=1} u_i},F \big) + \| u \|^2_{\mu^{(k)}_n}, & \text{if } k > 1. \\
	\end{dcases}
\end{equation}
The resulting sequence of accumulated solutions
\begin{equation*}
	\Big\{ u^{(k)} \coloneqq \sum^{k}_{i=1} u_i : k \in \mathbb{N} \Big\}
\end{equation*}
provides a multiscale representation which, with an appropriate choice of sequences $\mu^{(k)}_n$, can be made arbitrarily fine. 

The hierarchical model as formulated above is a slight generalisation of the originally proposed one, in the sense that the sequences of regularising functionals considered in \cite{TadNezVes04} are always of the form $\alpha^{-k} \mathcal{R}(\cdot)$. That is, the regularising functional of step $k$ is the same as the one of previous steps save for a smaller regularisation parameter. Requirement \eqref{eq:hiersequ} allows a more general setup.

\subsubsection*{Helmholtz Decomposition}
The Helmholtz decomposition theorem states that every continuously differentiable tangent vector field $w$ on the sphere can be uniquely represented as the sum of a consoidal (curl-free) and a toroidal (divergence-free) vector field. More precisely, there exist two uniquely determined tangent vector fields $w^{(2)}$ and $w^{(3)}$ satisfying
\begin{align*}
	\nabla_{\mathcal{S}} \cdot \big( w^{(2)} \times \nu \big) &= 0, \\
	\nabla_{\mathcal{S}} \cdot w^{(3)} &= 0, \\
	w^{(2)} + w^{(3)} &= w.
\end{align*}
See \cite[Sec.~5.3]{FreSchr09}, for example.

The projection method we use to numerically solve the variational problems presented above leads to solutions $w$ that are finite Fourier sums
\begin{equation*}
	w = \sum_{i,n,j} w_{nj}^{i} y_{nj}^{(i)},
\end{equation*}
where $w_{nj}^{i} \in \mathbb{R}$. Now, from the definition of basis \eqref{eq:basis}, and the fact that $\nabla_{\mathcal{S}} \cdot (\nabla_{\mathcal{S}} Y \times \nu ) = 0$ for all sufficiently smooth functions $Y$, the Helmholtz decomposition of $w$ is obtained immediately
\begin{align*}
	w &= \sum_{n,j} w_{nj}^{2} y_{nj}^{(2)} + \sum_{n,j} w_{nj}^{3} y_{nj}^{(3)} \\
	  &= 	\underbracket{ \nabla_{\mathcal{S}} \Big( \sum_{n,j} w_{nj}^{2}\lambda_n^{-1/2} Y_{nj} \Big)}_{w^{(2)}} +
	  		\underbracket{ \nabla_{\mathcal{S}} \Big( \sum_{n,j} w_{nj}^{3}\lambda_n^{-1/2} Y_{nj} \Big) \times \nu }_{w^{(3)}}.
\end{align*}

\section{Numerical Solution}\label{sec:numerics}
In the first subsection below we describe the numerical optimisation of the optical flow functional \eqref{eq:rega_n} based on the optical flow equation \eqref{eq:ofc} and explain the modifications necessary for the hierarchical decomposition. Subsection \ref{sec:u+v_num} is dedicated to the $u+v$ decomposition model. Finally, we explain how the resulting spherical integrals are approximated (Sec.~\ref{sec:discretisation}).

For convenience we relabel the orthonormal basis \eqref{eq:basis} using a single index $p\in\mathbb{N}$ and write, for example, $u = \sum_p u_p y_p$ from now on.

\subsection{Optical Flow}
Let $s \in \mathbb{R}$ and let $(\mu_n)_n$ be a sequence comparable to $(\lambda^s_n)_n$ in the sense of \eqref{eq:comparable}. Then, a minimiser of $\mathcal{E}_{\mu_n}$, if it exists, has to be in $H^s(\mathcal{S}, T\mathcal{S})$. We solve the problem of finding
\begin{equation*}
	\min_{u \in H^s(\mathcal{S}, T\mathcal{S})} \mathcal{E}_{\mu_n}(u)
\end{equation*}
by a projection method. That is, we let $u$ range only over a finite-dimensional subspace $\mathcal{U}$ of $H^s(\mathcal{S}, T\mathcal{S})$, where
\begin{equation*}
	\mathcal{U} = \mathrm{span}\{y_p : p \in I_\mathcal{U}\}
\end{equation*}
and $I_\mathcal{U} \subset \mathbb{N}$ is a finite index set. The unknown vector field can now be written as
\begin{equation}\label{eq:ansatz}
	u = \sum_{p \in I_\mathcal{U}} u_p y_p
\end{equation}
and the problem of finding an optimal $u \in H^s(\mathcal{S}, T\mathcal{S})$ simplifies to a minimisation problem over $\mathbb{R}^{|I_\mathcal{U}|}$. Plugging \eqref{eq:ansatz} into the optical flow functional gives
\begin{equation}\label{eq:discretefunctional}
	\mathcal{E}_{\mu_n}(u) = \int_\mathcal{S} \Big( \sum_{p \in I_\mathcal{U}} u_p (\nabla_\mathcal{S} F \cdot y_p) + \partial_t F \Big)^2 \, \mathrm{d}\mathcal{S} + \sum_{p \in I_\mathcal{U}} \mu_p u_p^2,
\end{equation}
which is minimal, if the optimality conditions $\partial \mathcal{E}_{\mu_n} / \partial u_p = 0$, for all $p\in I_{\mathcal{U}}$, are satisfied. They read
\begin{equation*}
	\sum_{q \in I_\mathcal{U}} u_q \int_\mathcal{S} (\nabla_\mathcal{S} F \cdot y_p) (\nabla_\mathcal{S} F \cdot y_q) \, \mathrm{d}\mathcal{S} + \mu_p u_p = - \int_\mathcal{S} \partial_t F \nabla_\mathcal{S} F \cdot y_p \, \mathrm{d}\mathcal{S}, \quad p\in I_{\mathcal{U}},
\end{equation*}
or in matrix-vector form
\begin{equation}\label{eq:linsys}
	(A + D) w = b,
\end{equation}
where  $w = (u_1,\ldots, u_{|I_\mathcal{U}|})^\top$ is the vector of unknown coefficients, the elements of matrix $A = (a_{pq})_{pq}$ read
\begin{equation}\label{eq:a_pq}
	a_{pq} = \int_\mathcal{S} (\nabla_{\mathcal{S}} F \cdot y_p) (\nabla_{\mathcal{S}} F \cdot y_q) \, \mathrm{d}\mathcal{S},
\end{equation}
$D = \mathrm{diag}(\mu_1,\ldots, \mu_{|I_\mathcal{U}|})$ is a diagonal matrix that corresponds to the regularisation term and the right hand side $b = (b_1,\ldots, b_{|I_\mathcal{U}|})^\top$ is given by
\begin{equation}\label{eq:b_p}
	b_p = -\int_\mathcal{S} \partial_t F \nabla_{\mathcal{S}} F \cdot y_p \, \mathrm{d}\mathcal{S}.
\end{equation}
With a slight abuse of notation we identified the set $I_{\mathcal{U}}$ with $\{1, \ldots, |I_{\mathcal{U}}| \}$ in the definitions of $A,D,w,b$. We continue to do so below.

\subsubsection*{Hierarchical Decomposition}
The hierarchical model only needs a minor modification for the case $k>1$. We can rewrite the data term from \eqref{eq:hierachical} as
\begin{equation*}
	\mathcal{D} \big( u + u^{(k-1)},F \big) = \| \nabla_{\mathcal{S}} F \cdot u + \partial_t \tilde{F} \|^2_{L^2(\mathcal{S})},
\end{equation*}
where $\partial_t \tilde{F} = \partial_t F + \nabla_\mathcal{S} F \cdot u^{(k-1)}$. Therefore, only the right hand side of the optimality system \eqref{eq:linsys} has to be updated in every step. For simplicity we can assume that the approximation space $\mathcal{U}$ is the same in every step, so that $u^{(k-1)}$ has the representation $\sum_{p\in I_{\mathcal{U}}} c^{k-1}_p y_p$, where the coefficients $c^{k-1}_p$ are already known from previous steps. Letting $b^k$ denote the right hand side of the optimality system for step $k$, we calculate
\begin{align*}
	b^k_p	&= -\int_\mathcal{S} \partial_t \tilde{F} \nabla_{\mathcal{S}} F \cdot y_p \, \mathrm{d}\mathcal{S} \\
			&= b_p - \sum_{q \in I_{\mathcal{U}}} c^{k-1}_q a_{pq},
\end{align*}
or simply
\begin{equation*}
	b^k = b - Ac^{k-1}.
\end{equation*}

\subsection{$u+v$ Decomposition}  \label{sec:u+v_num}
The projection approach explained above is easily adapted to the $u+v$ decomposition problem. Again, let $r \neq s$ be real numbers and choose two sequences $(\mu_n)_n$, $(\nu_n)_n$ so that $\| \cdot \|_{\mu_n}$ is a norm for $H^r(\mathcal{S}, T\mathcal{S})$ and $\| \cdot \|_{\nu_n}$ is one for $H^s(\mathcal{S}, T\mathcal{S})$. Now, we solve
\begin{equation*}
	\min_{(u,v) \in \mathcal{U} \times \mathcal{V}} \mathcal{E}_{\mu_n, \nu_n} (u,v),
\end{equation*}
where
\begin{align*}
	\mathcal{U} &= \mathrm{span}\{y_p : p \in I_\mathcal{U}\}, \\
	\mathcal{V} &= \mathrm{span}\{y_p : p \in I_\mathcal{V}\} 
\end{align*}
are finite dimensional spaces. Proceeding as in the previous section, we obtain the following optimality conditions
\begin{align*}
	\sum_{p \in I_\mathcal{U}} u_p a_{kp} + \sum_{q \in I_\mathcal{V}} v_q a_{kq} +	\mu_k u_k &= - b_k, \quad k \in I_\mathcal{U}, \\
	\sum_{p \in I_\mathcal{U}} u_p a_{\ell p} + \sum_{q \in I_\mathcal{V}} v_q a_{\ell q} +	\nu_\ell v_\ell &= - b_\ell, \quad \ell \in I_\mathcal{V},
\end{align*}
where the coefficients $a_{pq}$ and $b_p$ are as defined in \eqref{eq:a_pq} and \eqref{eq:b_p}, respectively. Concatenating the two coefficient vectors $(u_p)_p$ and $(v_q)_q$ into a single vector $w\in\mathbb{R}^{|I_\mathcal{U}| + |I_\mathcal{V}|}$, so that the $u_p$ occupy the first $|I_\mathcal{U}|$ entries while the $v_q$ occupy the last $|I_\mathcal{V}|$ entries, the linear system reads
\begin{equation*}
	\tilde{A} w = \tilde{b}.
\end{equation*}
The matrix $\tilde{A}$ is given by
\begin{align*}
	\tilde{A} = \begin{pmatrix}
		U + D_1 & W \\
		W^{\top} & V + D_2
	\end{pmatrix},
\end{align*}
where
\begin{align*}
	U &= (a_{pq})_{p,q\in I_\mathcal{U}}, \\
	V &= (a_{pq})_{p,q\in I_\mathcal{V}}, \\
	W &= (a_{pq})_{p\in I_\mathcal{U},q\in I_\mathcal{V}}, \\
	D_1 &= \mathrm{diag}(\mu_1,\ldots,\mu_{|I_\mathcal{U}|}), \\
	D_2 &= \mathrm{diag}(\nu_1,\ldots,\nu_{|I_\mathcal{V}|}),
\end{align*}
and $\tilde{b}$ is concatenated from two versions of $b$ in the same way as $w$.

\subsection{Evaluation of Integrals} \label{sec:discretisation}
It remains to discuss the numerical evaluation of the integrals \eqref{eq:a_pq}, \eqref{eq:b_p}. First, we approximate the $2$-sphere $\mathcal{S}$ with a polyhedron $\hat{\mathcal{S}} = (\mathcal{V}, \mathcal{T})$ defined by a set $\mathcal{V} = \{v_{1}, \dots v_{m}\} \subset \mathcal{S}$ of vertices and a set $\mathcal{T} = \{T_{1}, \dots, T_{n}\} \subset \mathcal{V} \times \mathcal{V} \times \mathcal{V}$ of triangular faces. Each triangle $T_{i} \in \mathcal{T}$ is associated with a tuple $(i_{1}, i_{2}, i_{3})$ identifying the corresponding vertices $(v_{i_{1}}, v_{i_{2}}, v_{i_{3}})$. How the triangulated sphere is obtained in practice, is explained in Sec.~\ref{sec:acquisition}.

Second, in every experiment data $F$ are given only at the vertices and for two time steps $t=0$ and $t=1$. We set $F_0(\cdot) \coloneqq F(0,\cdot)$ and $F_{1}(\cdot) \coloneqq F(1,\cdot)$ and extend both functions to all of $\hat{\mathcal{S}}$ by linear interpolation on every triangle. Thus we obtain two continuous piecewise linear functions $\hat F_0$, $\hat F_{1}$. The time derivative of $F$ is approximated by a simple forward difference
\begin{equation*}
	\partial_t \hat F = \hat F_1 - \hat F_0,
\end{equation*}
which is again piecewise linear on $\hat{\mathcal{S}}$. The surface gradient $\nabla_{\mathcal{S}} F$ is replaced by a vector field $\nabla_{\hat{\mathcal{S}}} \hat F$ that is constant on every triangle. It is given by

\begin{equation*}
	\nabla_{\hat{\mathcal{S}}} \hat F |_{T_i} = \left(\hat F(v_{i_{1}}) - \hat F(v_{i_{2}})\right) \frac{h_{i_{2}}}{|h_{i_{2}}|^{2}} + \left(\hat F(v_{i_{1}}) - \hat F(v_{i_{3}})\right) \frac{h_{i_{3}}}{|h_{i_{3}}|^{2}},
\end{equation*}
where $h_{i_{j}} \in \mathbb{R}^{3}$ is the height vector of the triangle $T_{i}$ pointing from vertex $v_{i_{j}}$ to the opposite side, compare \cite[Sec.~3.3.3]{BotKobPauAllLev10}.

Finally, we also replace the fully normalised scalar spherical harmonics $Y_{nj}$ by their piecewise linear approximations $\hat Y_{nj}$ defined on $\hat{\mathcal{S}}$. As before we obtain piecewise constant approximations $\hat y_p$ of $y_p$. The resulting approximated integrals read
\begin{align*}
	\int_\mathcal{S} (\nabla_{\mathcal{S}} F \cdot y_p) (\nabla_{\mathcal{S}} F \cdot y_q) \, \mathrm{d}\mathcal{S}
		&\approx	\int_{\hat{\mathcal{S}}} (\nabla_{\hat{\mathcal{S}}} \hat F \cdot \hat y_p) (\nabla_{\hat{\mathcal{S}}} \hat F \cdot \hat y_q) \, \mathrm{d}\hat{\mathcal{S}} \\
		&=			\sum_{i=1}^n (\nabla_{\hat{\mathcal{S}}} \hat F|_{T_i} \cdot \hat y_p|_{T_i}) (\nabla_{\hat{\mathcal{S}}} \hat F|_{T_i} \cdot \hat y_q|_{T_i}) A_i,
\end{align*}
where $A_i$ denotes the area of $T_i$, and
\begin{align*}
	\int_\mathcal{S} \partial_t F \nabla_{\mathcal{S}} F \cdot y_p \, \mathrm{d}\mathcal{S}
		&\approx	\int_{\hat{\mathcal{S}}} \partial_t \hat F \nabla_{\hat{\mathcal{S}}} \hat F \cdot \hat y_p \, \mathrm{d}\hat{\mathcal{S}}  \\
		&=			\sum_{i=1}^n  (\nabla_{\hat{\mathcal{S}}} \hat F|_{T_i} \cdot \hat y_p|_{T_i}) \int_{T_i} \partial_t \hat F \, \mathrm{d}T_i \\
		&=			\sum_{i=1}^n  (\nabla_{\hat{\mathcal{S}}} \hat F|_{T_i} \cdot \hat y_p|_{T_i}) \frac{A_i}{3}\sum_{j=1}^3 \partial_t \hat F (v_{i_j}).
\end{align*}

\section{Experiments}\label{sec:experiments}

\begin{figure}
	\centering
	\includegraphics[width=0.49\textwidth]{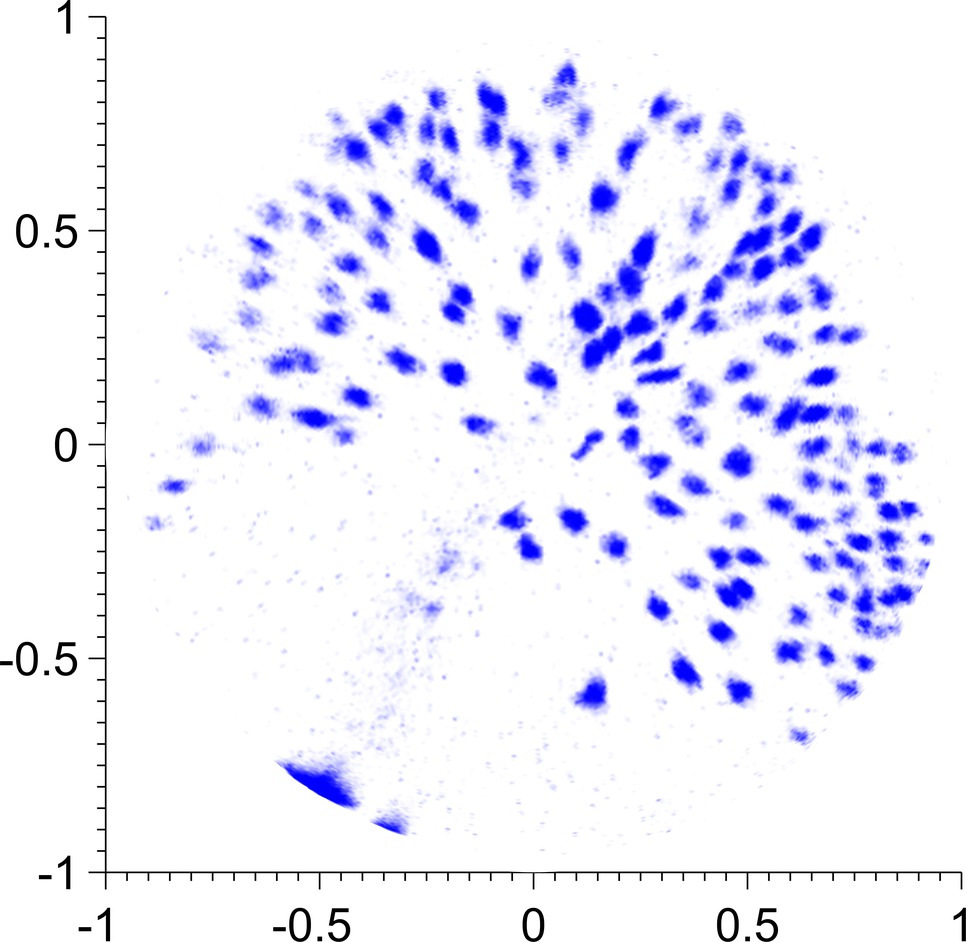} \hfill
	\includegraphics[width=0.49\textwidth]{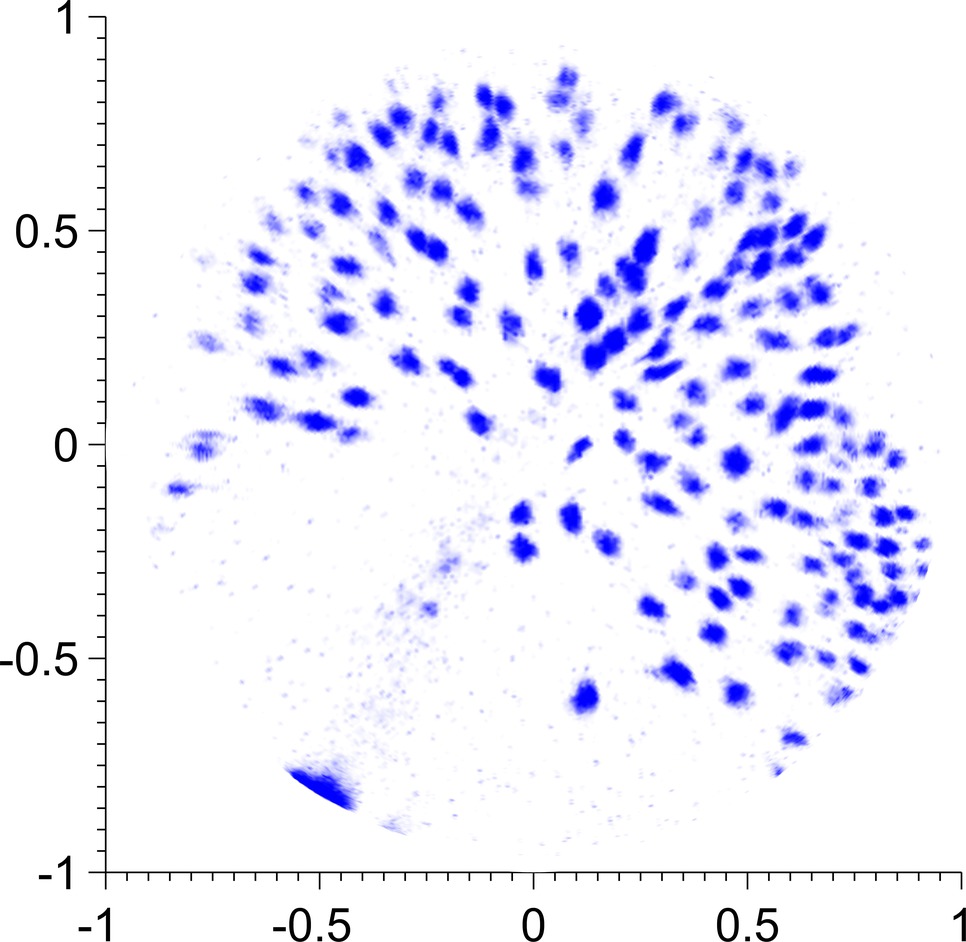}
	\caption{Top view of frames no. 57 (left) and 58 (right) of the processed zebrafish images. The embryo's body axis runs from bottom left to top right.}
	\label{fig:data2d}
\end{figure}

\subsection{Description of Microscopy Data} \label{sec:data}
The data which motivated the study of the proposed decomposition models are time-lapse volumetric (4D) images. The obtained sequence depicts a live zebrafish embryo during the gastrula period, approximately five to ten hours after fertilisation. With the help of confocal laser-scanning microscopy, endoderm cells expressing a green fluorescence protein were recorded separately from the background. For details on the imaging techniques and the fluorescence marker we refer to~\cite{MegFra03} and~\cite{MizVerHeaKurKik08}, respectively.

The sequence obtained by the microscope captures a cuboid region of approximately $860 \times 860 \times 340 \, \mu \mathrm{m}^{3}$. The spatial resolution is $512 \times 512 \times 44$~voxels and the intensity is in the range $\{0, \dots, 255\}$. A total number of 77 images were taken, one every $240 \, \mathrm{s}$. In the following, the microscopic data will be denoted by
\begin{equation*}
	F^{\delta} \in \{0, \dots, 255\}^{77 \times 512 \times 512 \times 44}.
\end{equation*}

During this early stage, endodermal cell proliferation is known to take place on a so-called monolayer~\cite{WarNus99}. In other words, cells move and divide without stacking and, as the yolk is ball-shaped, admit for the extraction of a spherical image sequence. For further explanations and numerous illustrations of the developmental process of zebrafish embryos we refer to~\cite{KimBalKimUllSchi95}.

\subsection{Acquisition of Spherical Data} \label{sec:acquisition}
We extracted spherical images from the dataset by first fitting a sphere to the approximate cell centres in each pair of consecutive frames. For simplicity we restrict our attention to one such pair of frames which we denote by $F_0^\delta, F_1^\delta$. Cell centres are typically characterised by local maxima in intensity and can be found by applying a Gaussian filter and simple thresholding. Without loss of generality the radius of the fitted sphere is assumed to be $1$. In a second step, we created a point grid $\mathcal{V} \subset \mathcal{S}$ starting from an icosahedron inscribed in the sphere. In each iteration every triangular face is split into four sub-triangles by connecting the edge midpoints with each other and projecting them onto the sphere. Thus, the total number of faces is $20 \cdot 4^{k}$, where $k$ is the number of refinements. In our experiments we found that $k = 7$ iterations suffice.

In order to project the volumetric time-lapse data $F_j^{\delta}$, onto the grid $\mathcal{V}$, we define
\begin{equation*}
	\hat{F}_j(v_{i}) := \max_{c \in [1-\epsilon, 1+\epsilon]} \bar{F}_j^{\delta}(c v_{i}),
\end{equation*}
for the said pair of consecutive frames $j=0,1$. Here $\bar{F}^{\delta}_j(x)$ is a piecewise linear extension to $\mathbb{R}^3$ of $F_j^{\delta}$ and $\epsilon >0$ is sufficiently large. Deviations of the monolayer from a perfect sphere are thereby corrected. The obtained data are subsequently scaled to the range $[0, 1]$. Note that, in contrast to our previous work~\cite{KirLanSch13}, here we consider the unfiltered microscopy data for optical flow estimation.

The support of the obtained data is contained in the northern hemisphere. Thus, it suffices to consider only half of the triangulation leading to a total number of around 164000 faces. Figure~\ref{fig:data3d} depicts a sample of two frames $\hat{F}_j$. In Fig.~\ref{fig:data2d}, a top view of the two frames is shown.

\subsection{Visualisation of Tangent Vector Fields}
\begin{figure}
	\centering
	\includegraphics[width=0.32\textwidth]{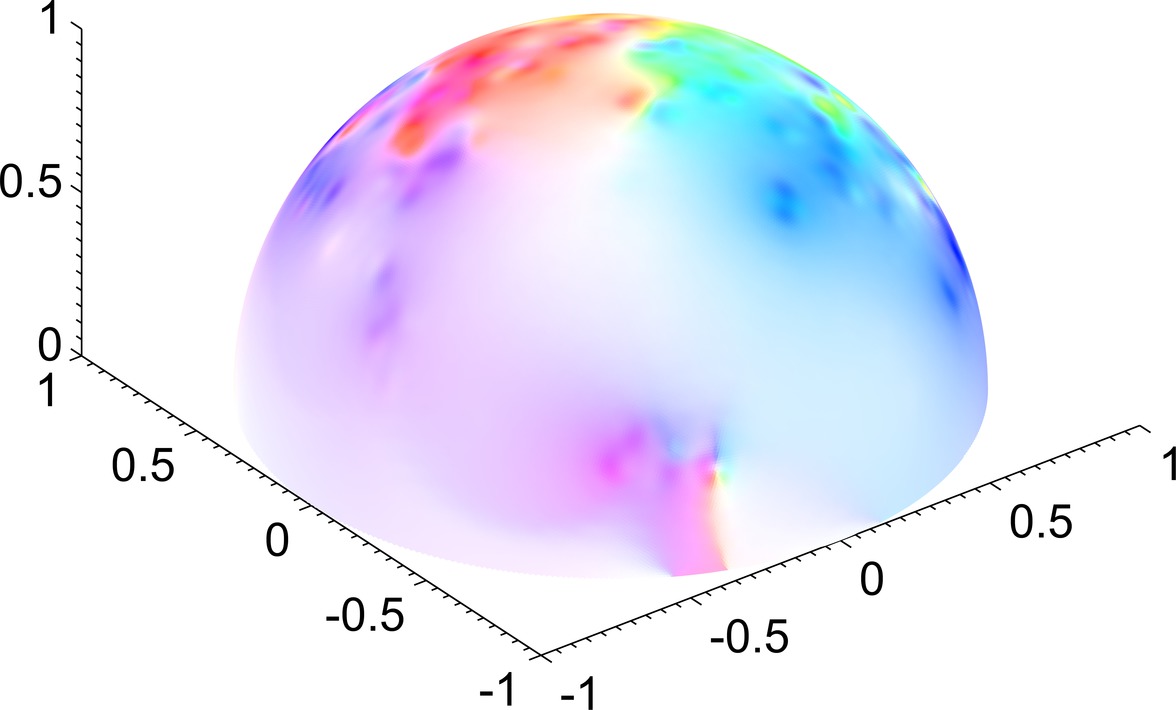} \hfill
	\includegraphics[width=0.15\textwidth]{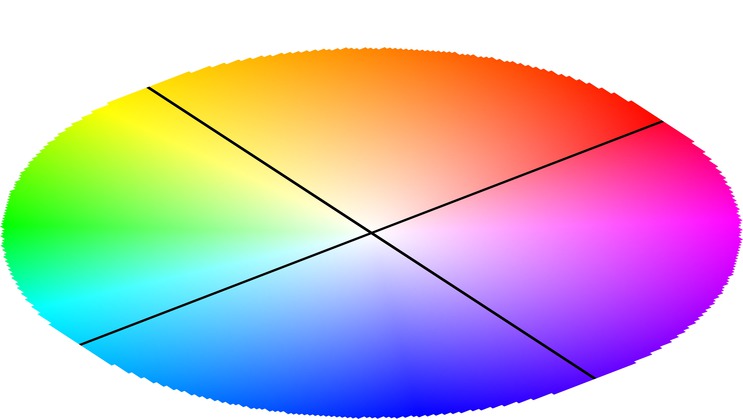} \hfill
	\includegraphics[width=0.32\textwidth]{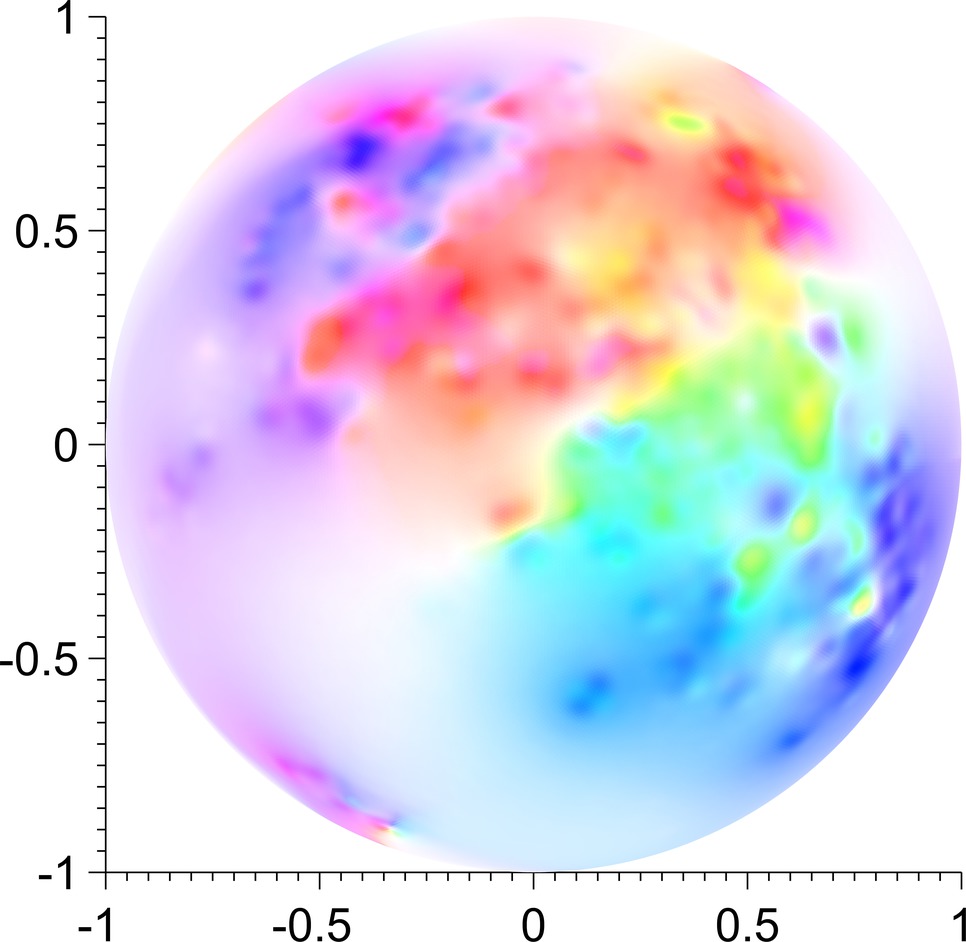} \hfill
	\includegraphics[width=0.15\textwidth]{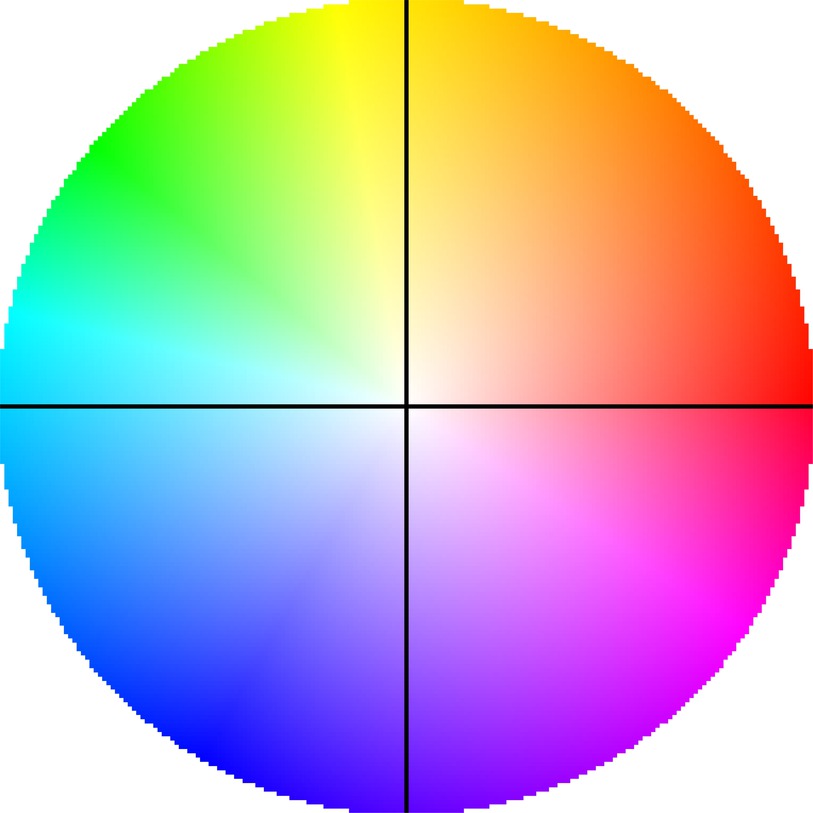}
	\caption{The first image shows a tangential velocity field with the adjusted colour-coding. The second image depicts the colour space in the unit circle. The third and fourth images portray the vector field and the colour space, respectively, but in a top view.}
	\label{fig:colourcoding}
\end{figure}

In order to visualise our results we will apply the standard flow colour-coding~\cite{BakSchaLewRotBla11} using a colour disk. Figure~\ref{fig:colourcoding} (rightmost image) depicts this colour space. Each vector is assigned a colour determined by its angle and length. However, this colour-coding is defined for planar vector fields only. As a possible remedy we suggest to first project tangential velocities to the plane and then correct the length. To this end, let $\mathrm{P}_{x_3}: ( x_{1}, x_{2},x_{3})^\top \mapsto ( x_{1}, x_{2},0)^\top$ be the orthogonal projector of $\mathbb{R}^3$ onto the $x_1$-$x_2$-plane. Accordingly, given a tangent vector field $v$, the planar vector field which we visualise is
\begin{equation*}
	\frac{|v|}{|\mathrm{P}_{x_3}v|} \mathrm{P}_{x_3} v.
\end{equation*}
This construction is chosen so that it preserves the length of $v$. This additional rescaling is different to~\cite{KirLanSch13}. The resulting colour image is finally mapped back onto the hemisphere. Figure.~\ref{fig:colourcoding} shows a tangent vector field visualised with the proposed approach.\footnote{Some figures may appear in colour only in the online version of this article.} From now on we will visualise velocity fields only in top view, as in the right hand side of Fig.~\ref{fig:colourcoding}. In addition, for every figure the colour disk's radius $R$ was chosen to be equal to the length of the longest vector under consideration. Specific values of $R$ are given in Table~\ref{tab:radii}.

\begin{table}
	\scriptsize
	\centering
	\begin{tabular}{@{}l@{\hspace{1em}}@{}l@{\hspace{1em}}@{}l@{\hspace{1em}}@{}l@{\hspace{1em}}@{}l@{\hspace{1em}}@{}l@{\hspace{1em}}@{}l@{\hspace{1em}}@{}l@{\hspace{1em}}@{}l@{\hspace{1em}}@{}l@{\hspace{1em}}@{}l@{}}
		Figure & \ref{fig:of1} left & \ref{fig:of1} mid & \ref{fig:of1} right & \ref{fig:of2} left & \ref{fig:of2} mid & \ref{fig:of2} right & \ref{fig:u+v} & \ref{fig:u+v2} & \ref{fig:hierarchical} & \ref{fig:hierarchical2} \\
		\hline
		$R$ & $0.0081$ & $0.0046$ & $0.0014$ & $0.0110$ & $0.0045$ & $0.0009$ & $0.0204$ & $0.0102$ & $0.0280$ & $0.0267$
	\end{tabular}
	\caption{Radii $R$ of the colour disks used in the different experiments below.}
	\label{tab:radii}
\end{table}

As a second way of illustrating steady velocity fields, we employ streamlines, see e.g.~\cite{WeiErl05}. In all our experiments we consider time as fixed and compute the optical flow $v$ for one pair of frames, cf.~Sec.~\ref{sec:reg}. Given a tangential velocity field $v$ and a starting point $x_{0} \in \mathcal{S}$, a streamline $\gamma(\cdot, x_{0})$ on $\mathcal{S}$ is the solution to the ordinary differential equation
\begin{equation}
	\begin{aligned}
		\partial_{\tau} \gamma(\tau, x_{0}) & = v(\gamma(\tau, x_{0})), \\
		\gamma(0, x_{0}) & = x_{0}.
	\end{aligned}
\label{eq:streamline}
\end{equation}
Numerically, we approximated~\eqref{eq:streamline} by solving
\begin{equation*}
	\begin{aligned}
		\hat{\gamma}(\tau + 1, x_{0}) & = \hat{\gamma}(\tau, x_{0}) + h v(\hat{\gamma}(\tau, x_{0})), \\
		\hat{\gamma}(0, x_{0}) & = x_{0},
	\end{aligned}
\end{equation*}
where $h$ is a step size, for a number of approximately 1300 initial points $x_{0} \in \hat{\mathcal{S}}$ and $\tau =50$ iterations. The step size was chosen as $h = (10{\lVert v \rVert}_{L^{\infty}(\mathcal{S}, T\mathcal{S})})^{-1}$. This use of integral curves is different from~\cite{KirLanSch13a_report}, where we computed approximate cell trajectories in a nonsteady velocity field. The visualisation by means of the colour coding is rich in detail and is even capable of indicating individual cell motion. Nevertheless, it fails to deliver intuition about the Helmholtz decomposition. Streamlines provide the anticipated effect.

\subsection{Experimental Results} \label{sec:results}
We performed numerous experiments and minimised functionals~\eqref{eq:rega_n},~\eqref{eq:u+v}, and~\eqref{eq:hierachical} as outlined in Sec.~\ref{sec:numerics} for the two frames shown in Fig.~\ref{fig:data3d} and Fig.~\ref{fig:data2d}. In all experiments the finite-dimensional spaces introduced in Sec.~\ref{sec:numerics} were chosen as
\begin{equation*}
	\mathcal{U} = \mathcal{V} = \mathrm{span} \Bigl\{y_{nj}^{(i)}: 1 \le n \le 100, 1 \le j \le 2n + 1, i = 2, 3\Bigr\}.
\end{equation*}
All resulting linear systems were solved using the Generalized Minimal Residual Method (GMRES) on an Intel Xeon E5-1620 $3.6 \, \mathrm{GHz}$ workstation with $128 \, \mathrm{GB}$~RAM. Solutions converged to a relative residual of 0.02 within 100 iterations. The runtime was governed by the evaluation of the integrals, cf.~Sec.~\ref{sec:discretisation}, and amounts to approximately five hours for the chosen bases and the chosen triangulation. Nevertheless, once the integrals are computed they can be used in all of the proposed models and the linear systems can be solved in a few seconds for different parameters and different norms. Our Matlab implementation and the data are available on our website.\footnote{\url{http://www.csc.univie.ac.at}}

\subsubsection{Optical Flow}

\begin{figure}
	\centering
	\includegraphics[width=0.32\textwidth]{of-flow2-frames-114-116-unfiltered-1-100-7-89} \hfill
	\includegraphics[width=0.32\textwidth]{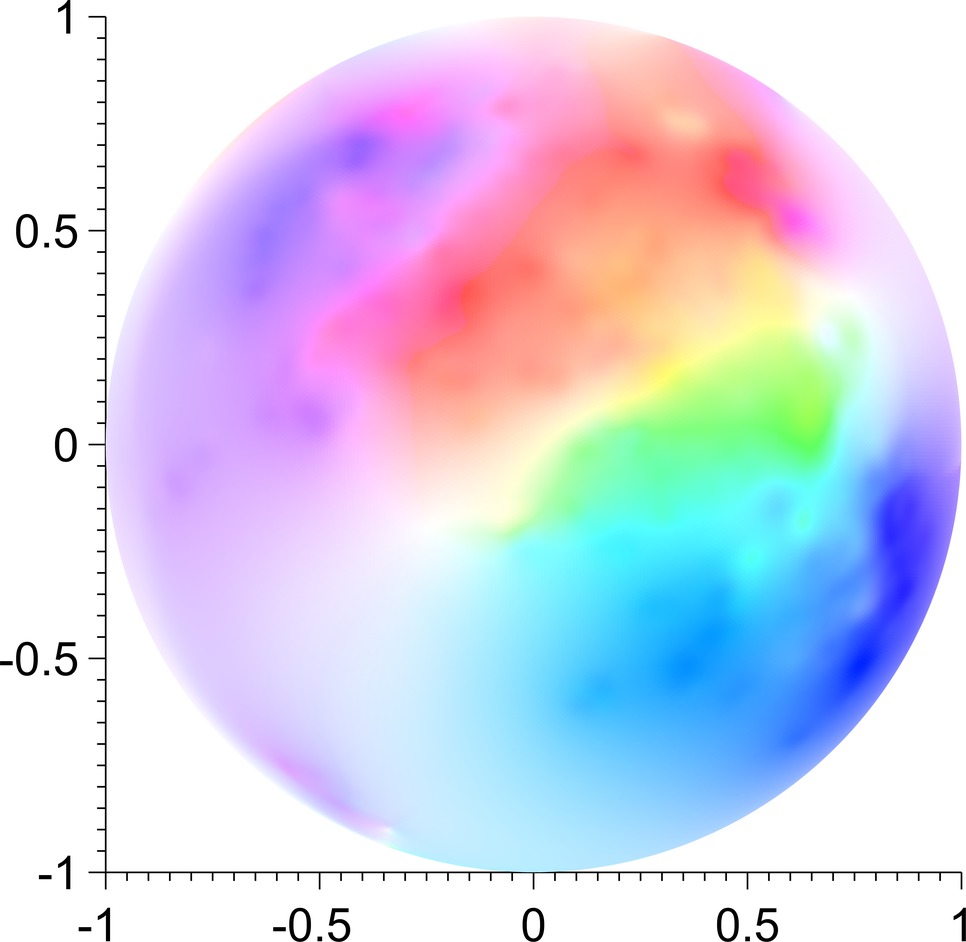} \hfill
	\includegraphics[width=0.32\textwidth]{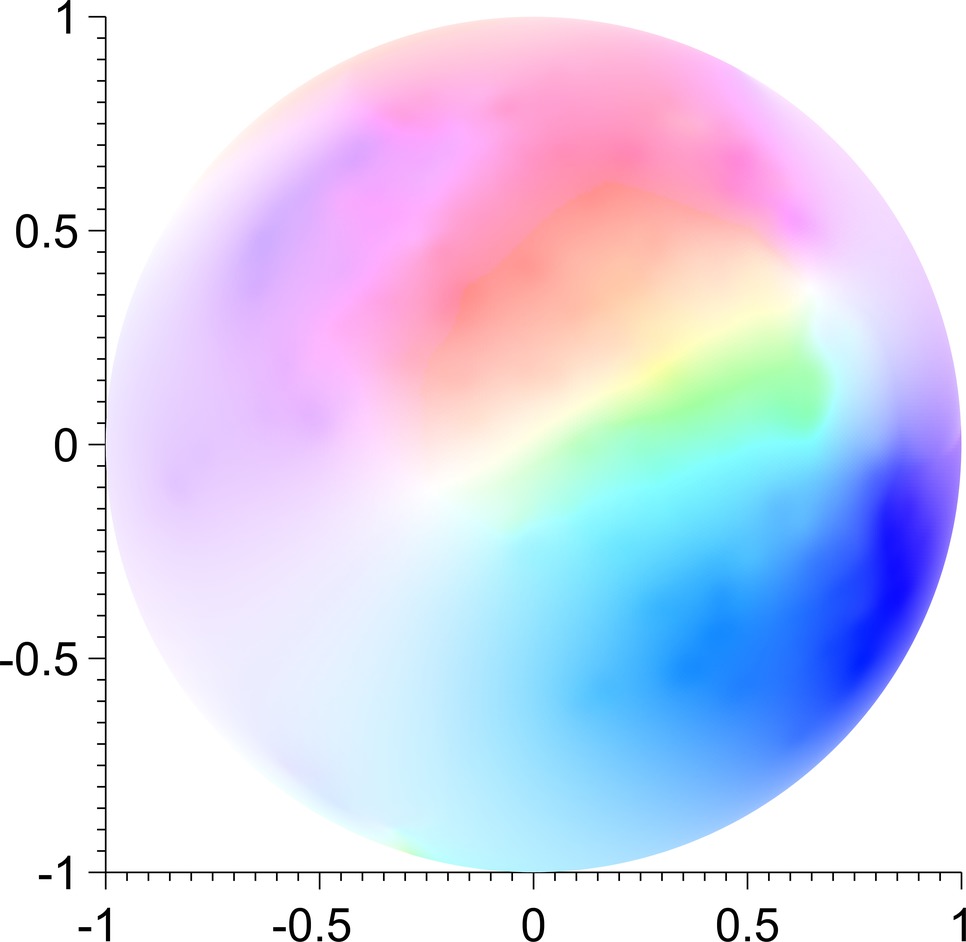} \\
	\smallskip
	\includegraphics[width=0.32\textwidth]{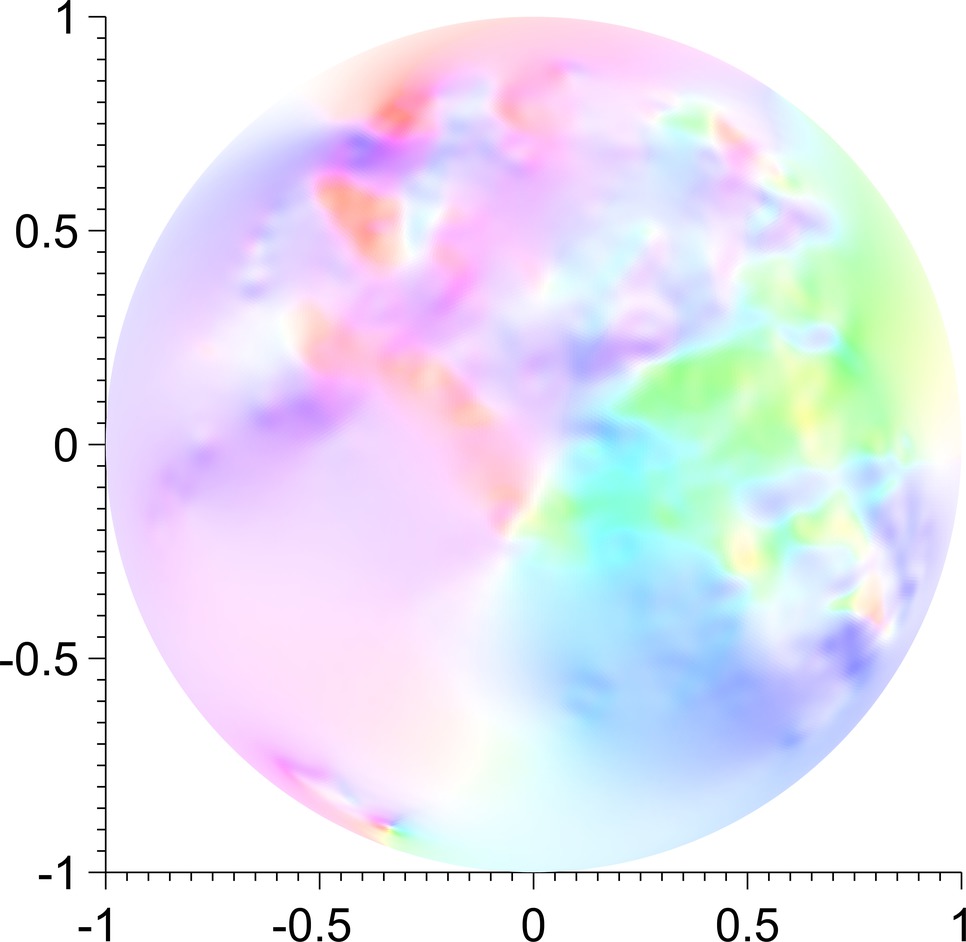} \hfill
	\includegraphics[width=0.32\textwidth]{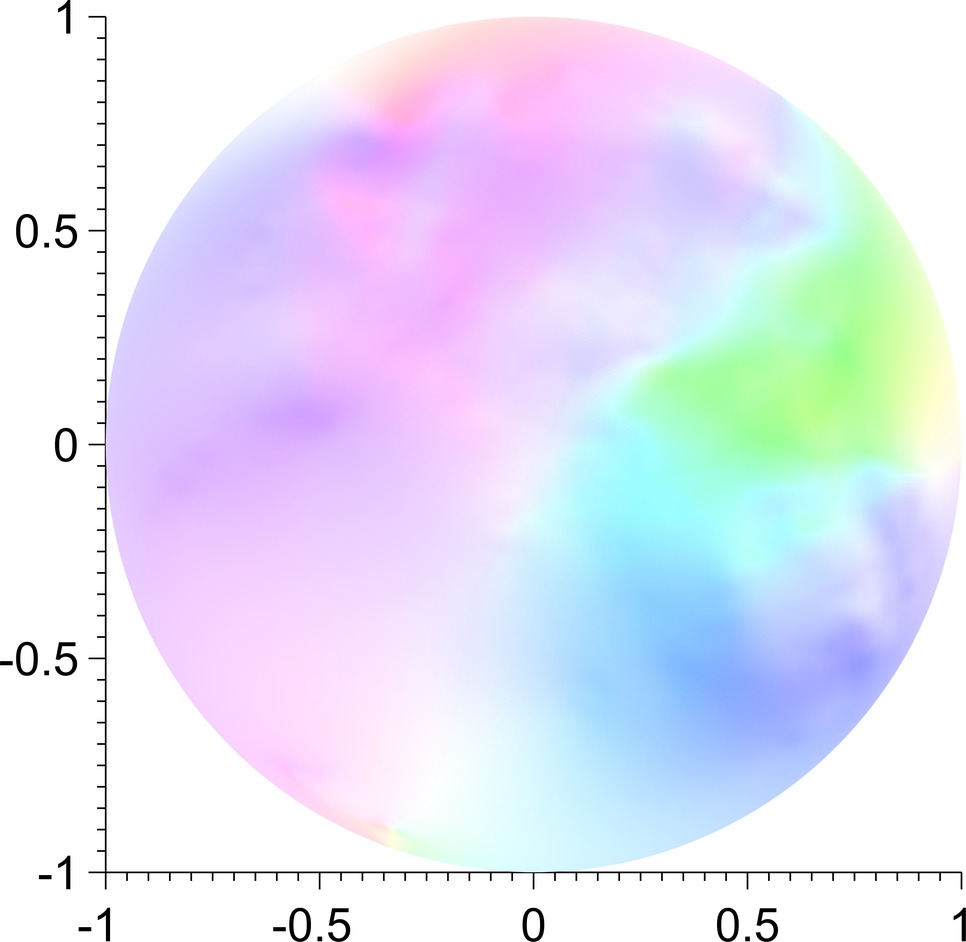} \hfill
	\includegraphics[width=0.32\textwidth]{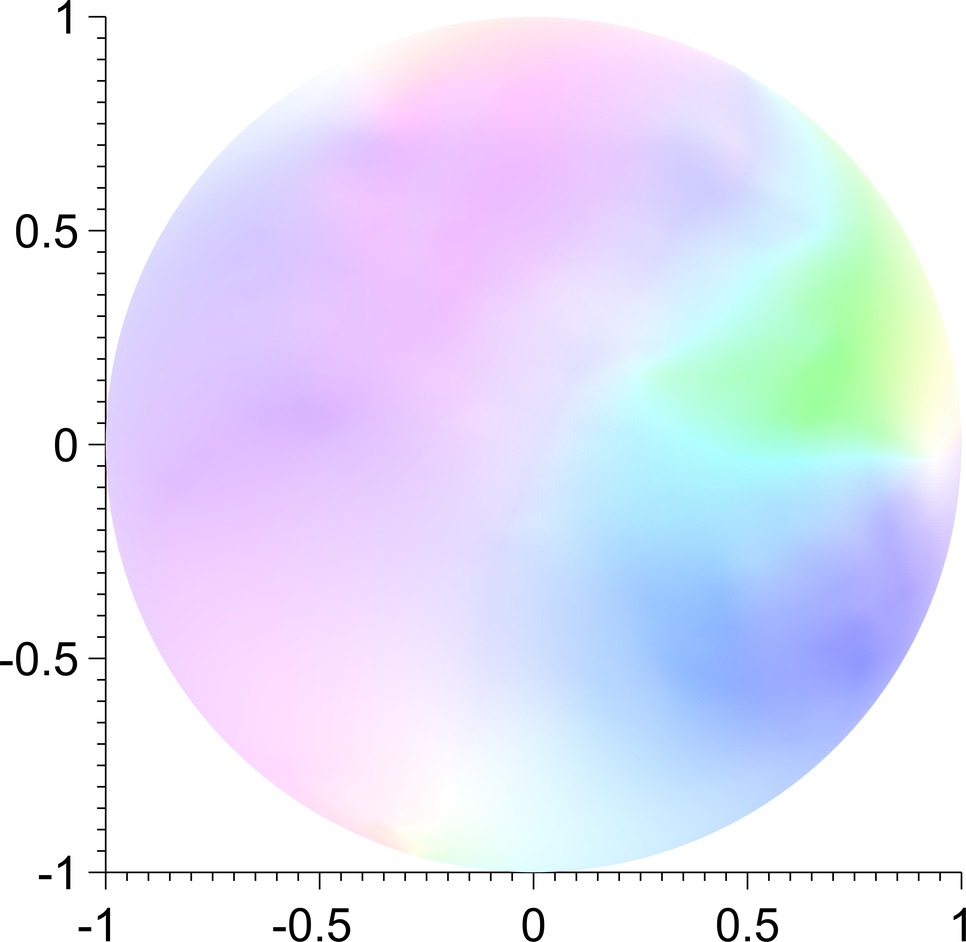} \\
	\smallskip
	\includegraphics[width=0.32\textwidth]{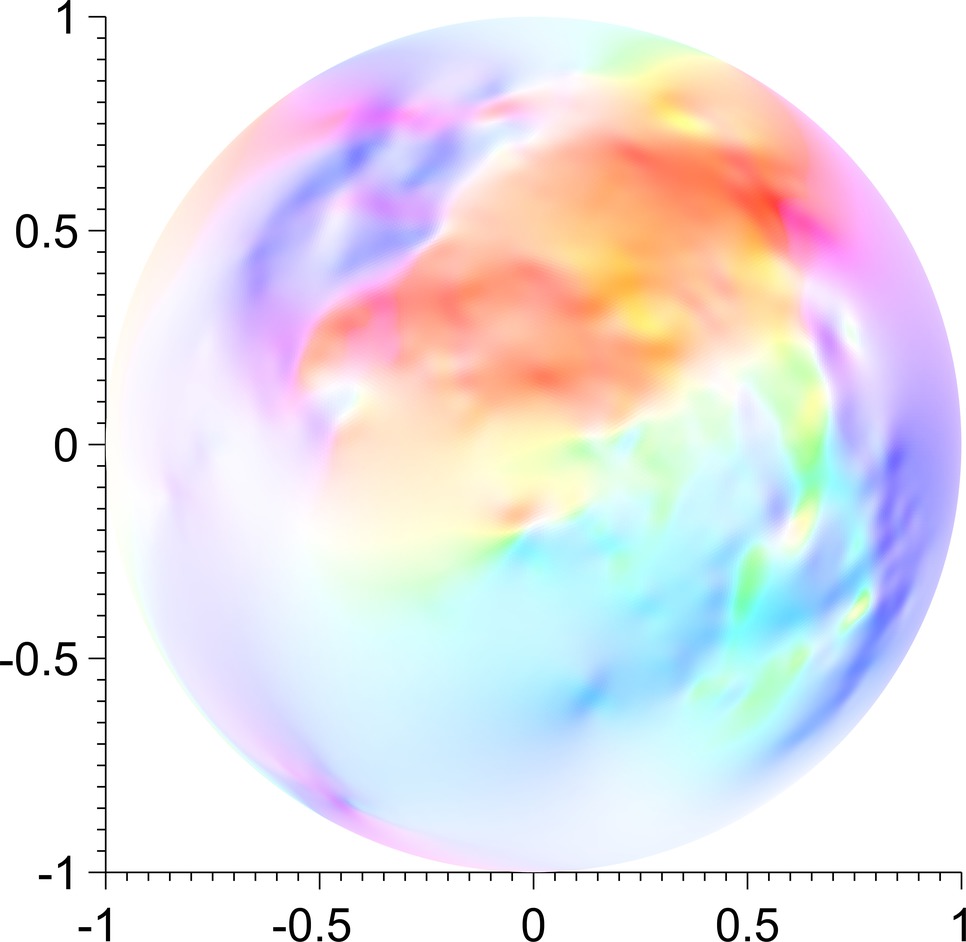} \hfill
	\includegraphics[width=0.32\textwidth]{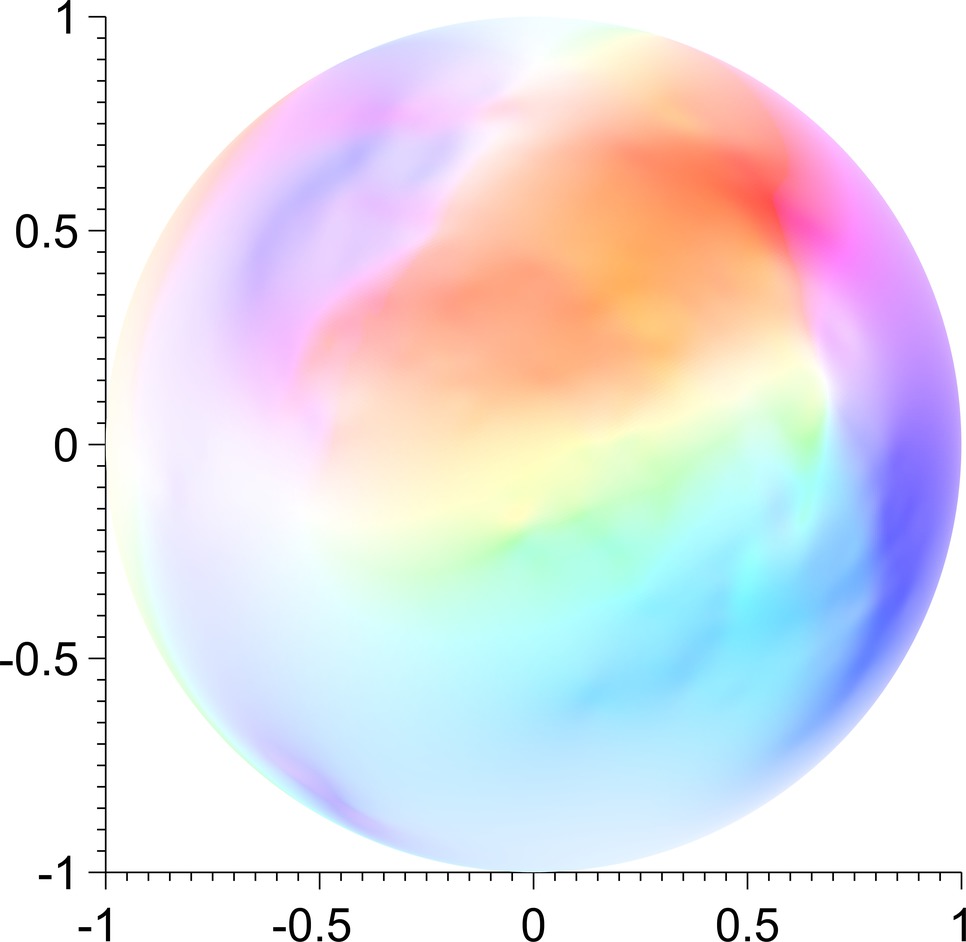} \hfill
	\includegraphics[width=0.32\textwidth]{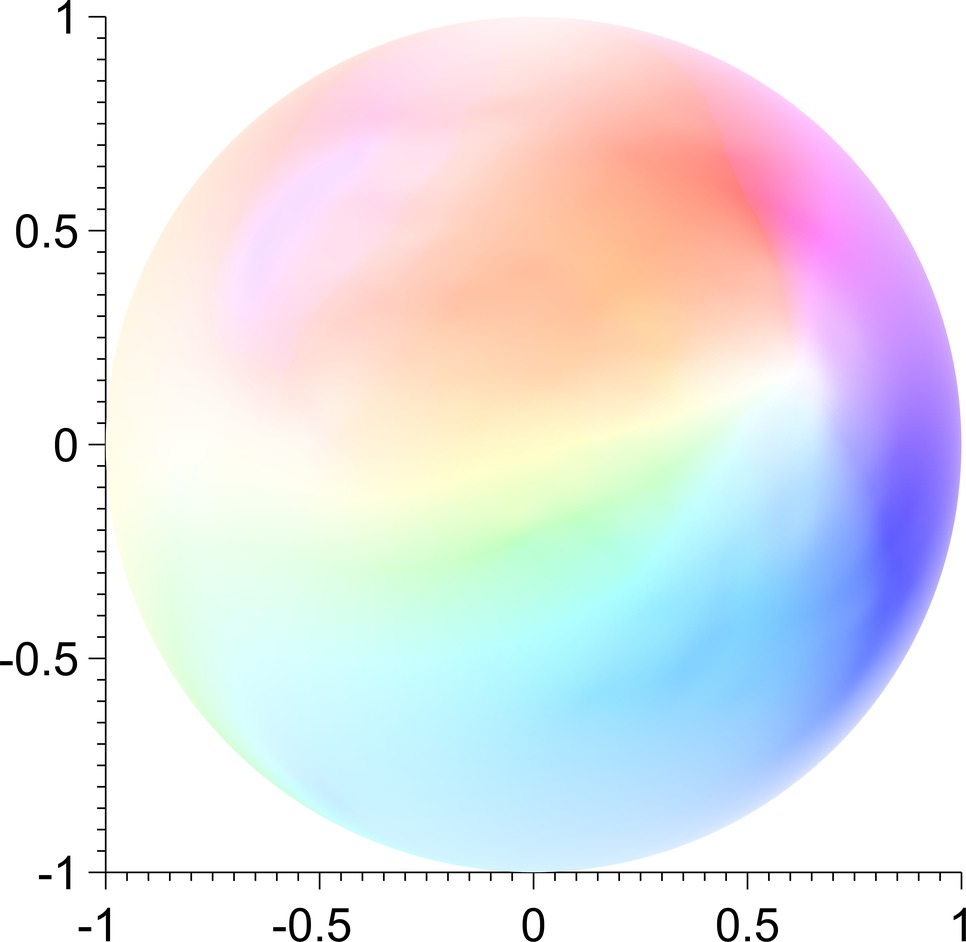}
	\caption{Minimiser of $\mathcal{E}_{\mu_n}$ (top), with $\mu_{n} = \alpha \lambda_{n}^{s}$, $s = 1$, and increasing values $\alpha = 1$, $\alpha = 10$, and $\alpha = 100$ from left to right. The middle row depicts the curl-free component whereas the bottom row depicts the divergence-free component. The embryo's body axis roughly runs from bottom left to top right in all images.}
	\label{fig:of1}
\end{figure}

\begin{figure}
	\centering
	\includegraphics[width=0.32\textwidth]{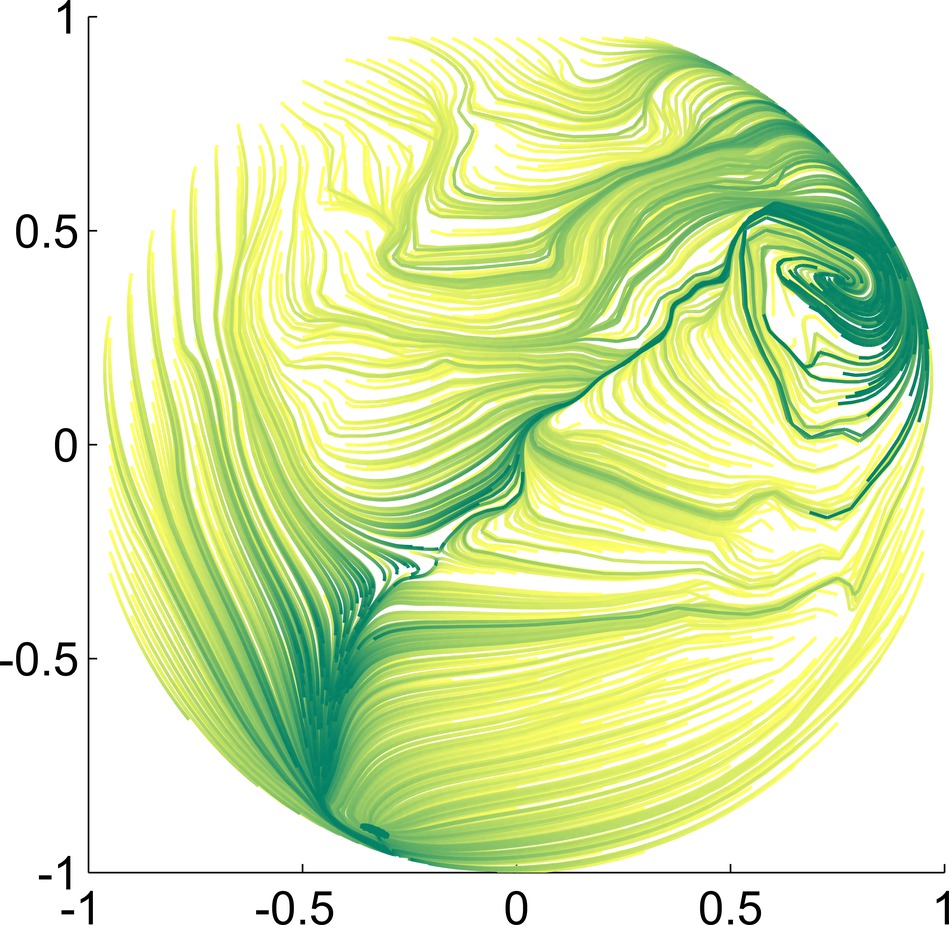} \hfill
	\includegraphics[width=0.32\textwidth]{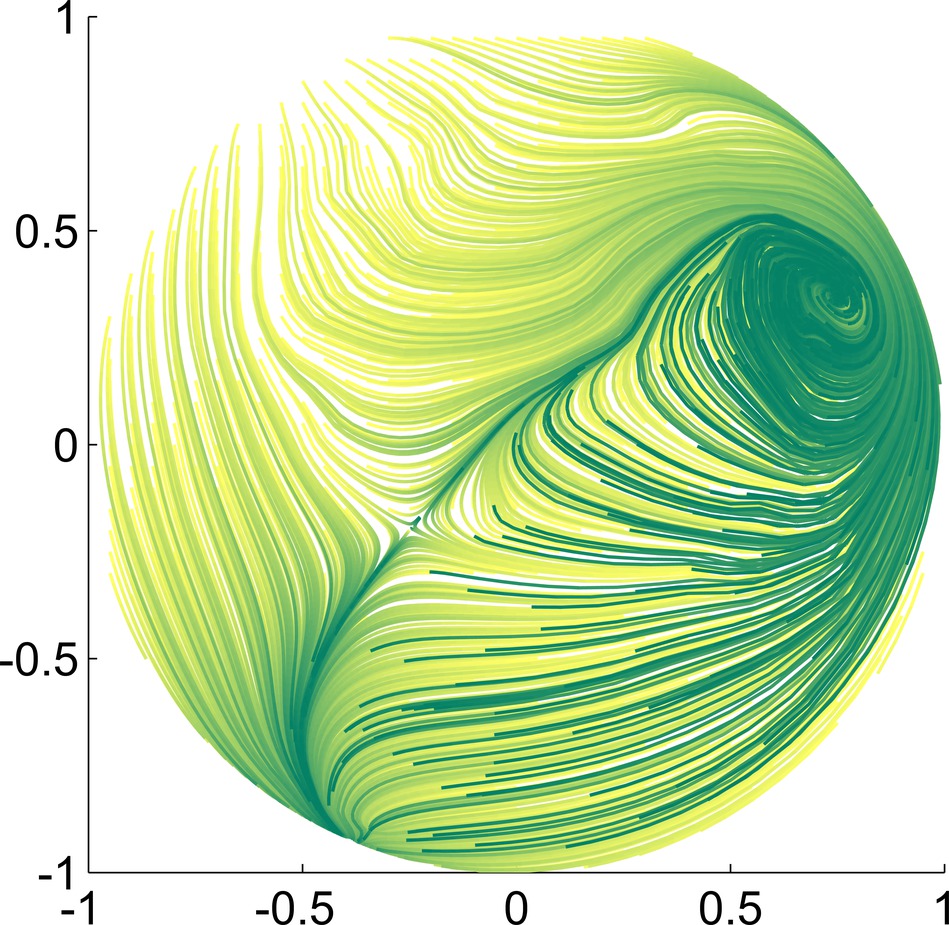} \hfill
	\includegraphics[width=0.32\textwidth]{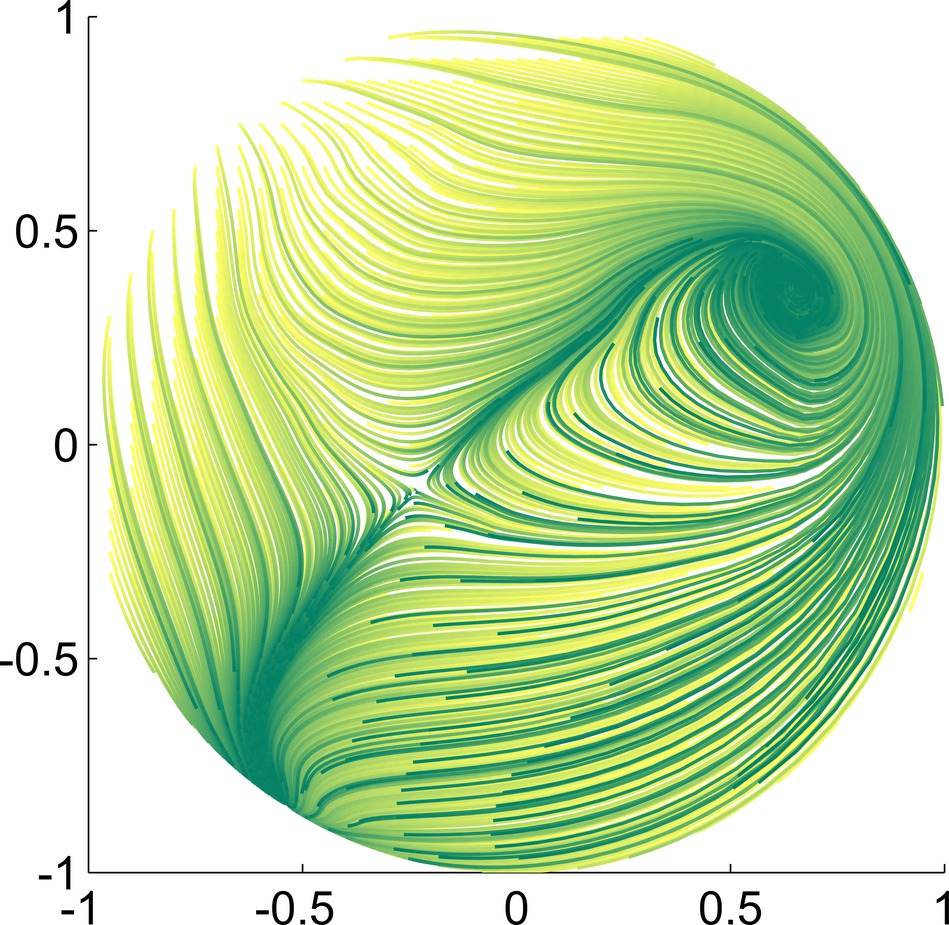} \\
	\smallskip
	\includegraphics[width=0.32\textwidth]{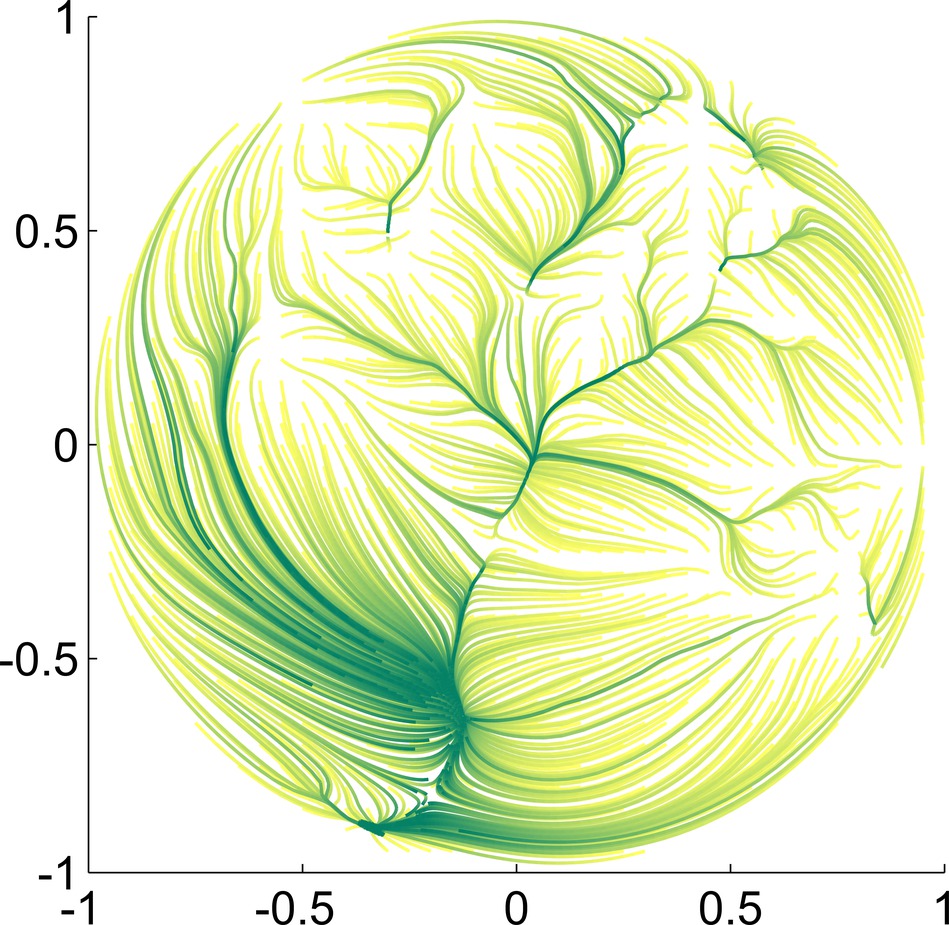} \hfill
	\includegraphics[width=0.32\textwidth]{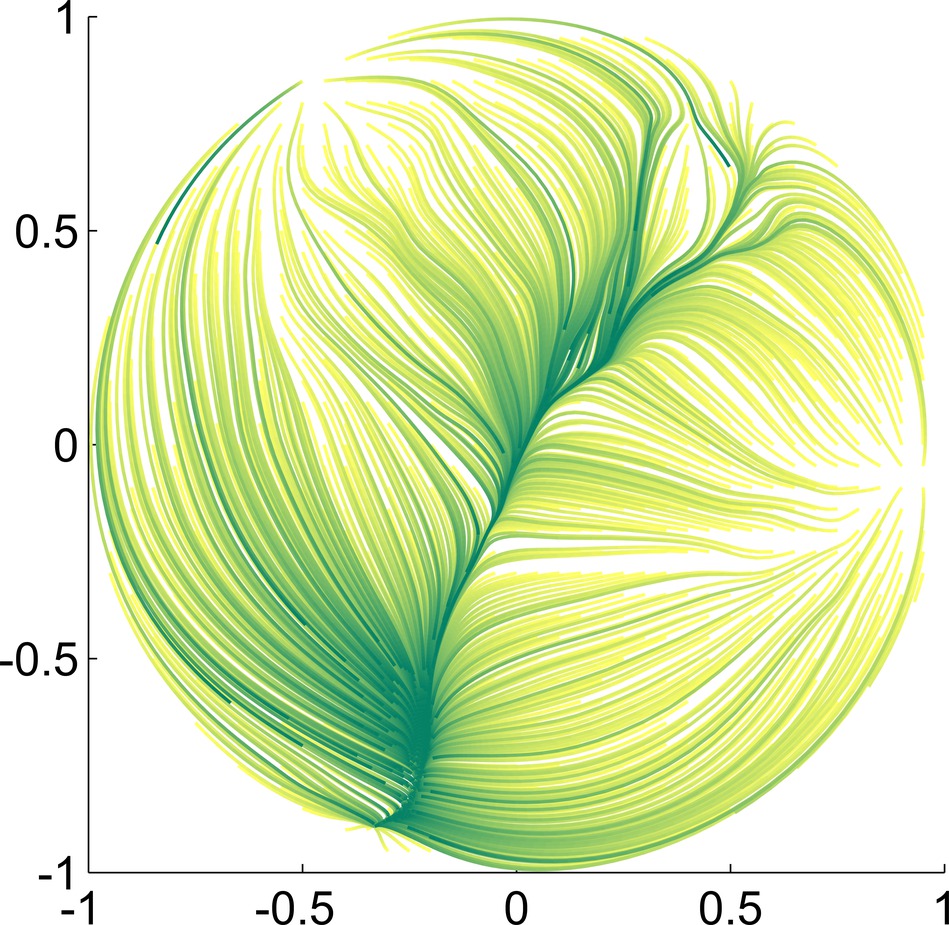} \hfill
	\includegraphics[width=0.32\textwidth]{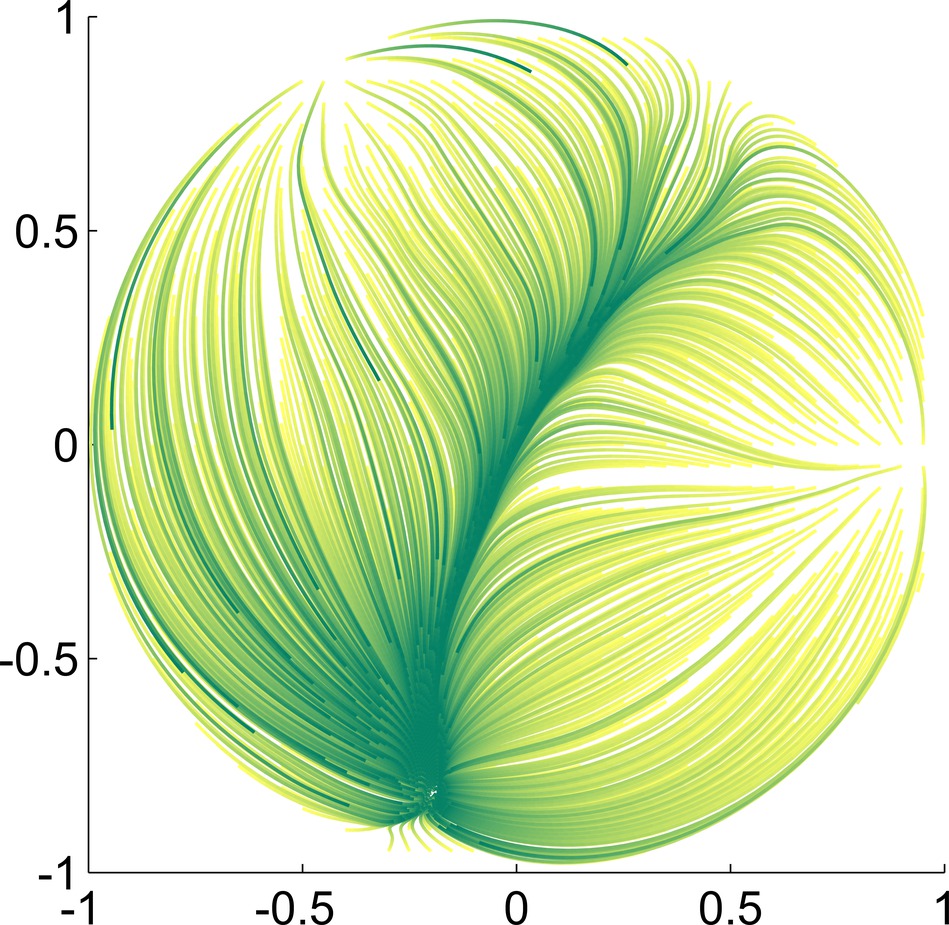} \\
	\smallskip
	\includegraphics[width=0.32\textwidth]{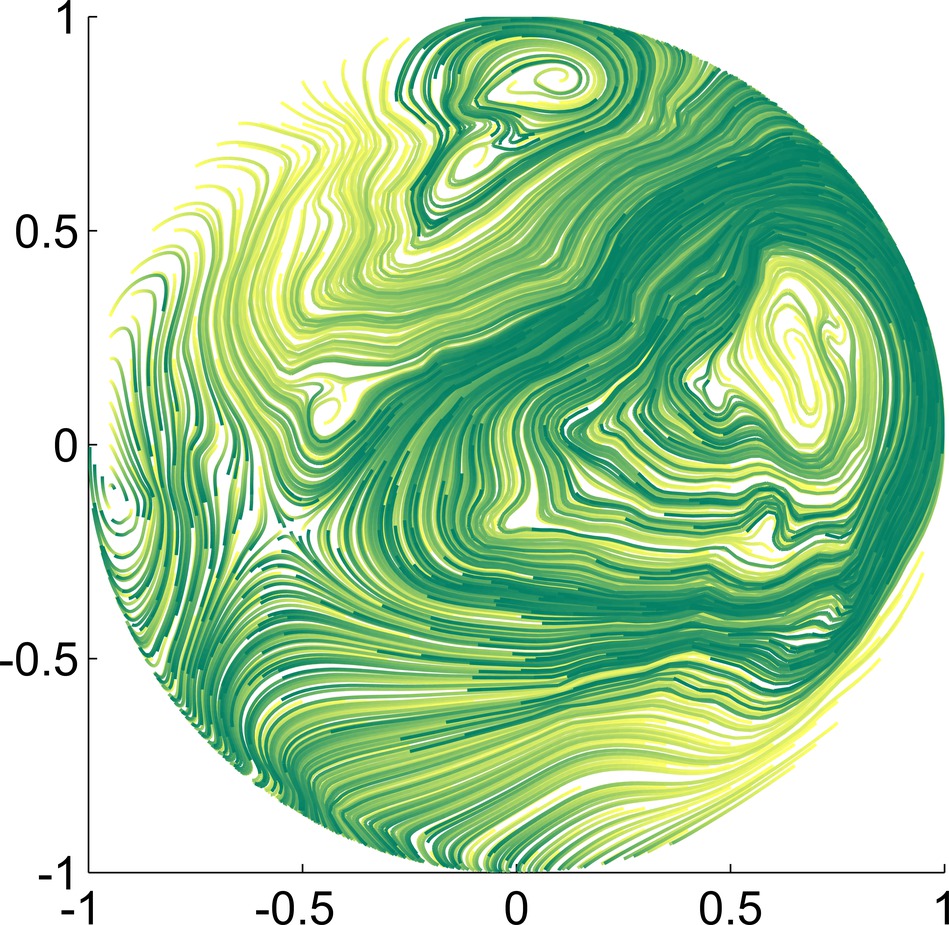} \hfill
	\includegraphics[width=0.32\textwidth]{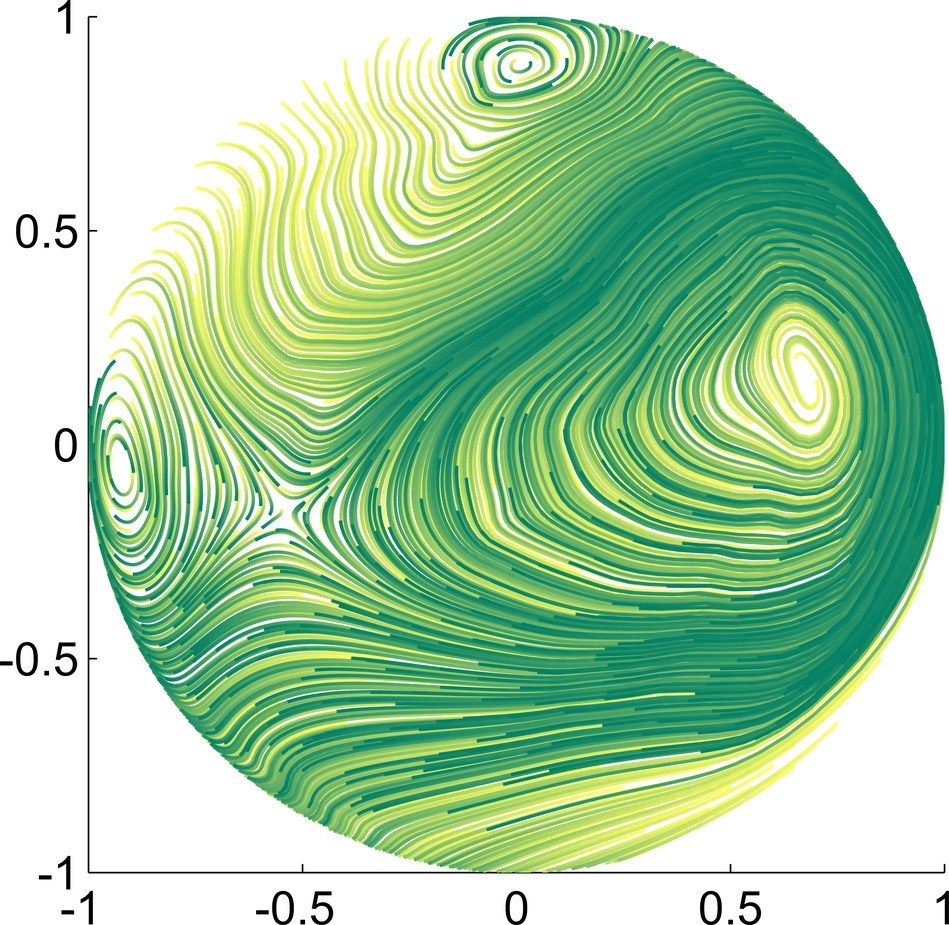} \hfill
	\includegraphics[width=0.32\textwidth]{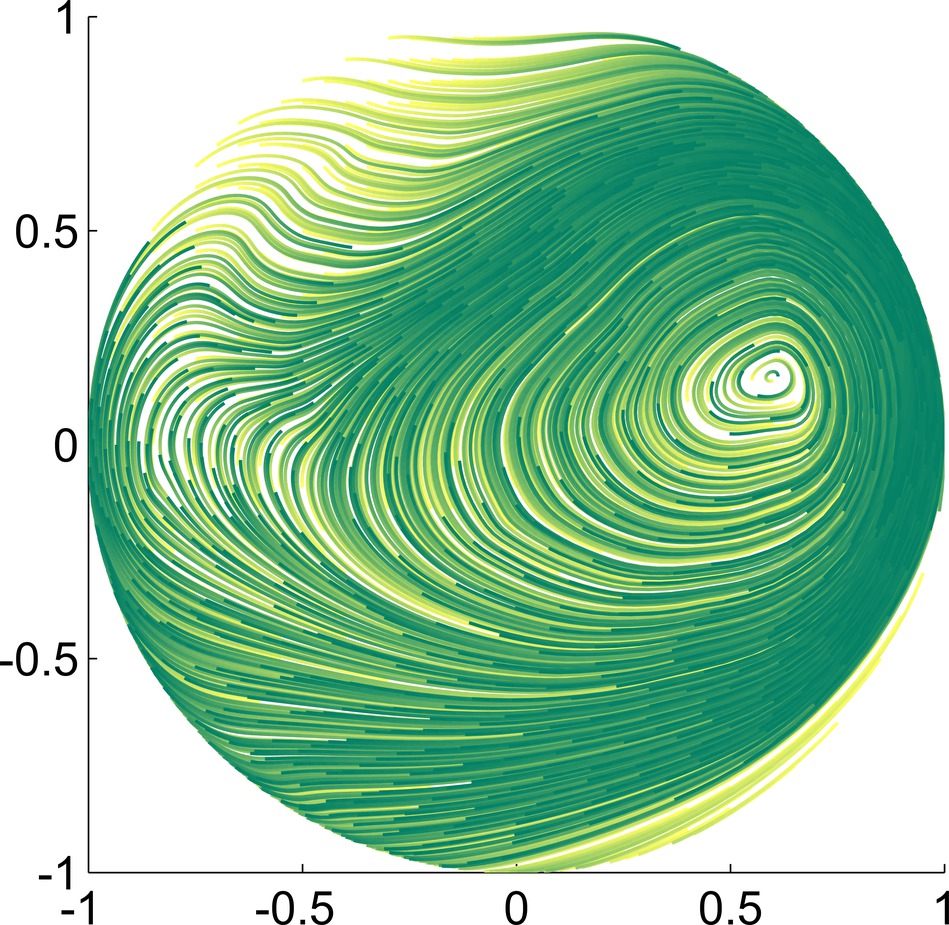}
	\caption{Streamlines illustrating the velocity fields from Fig.~\ref{fig:of1}. With increasing parameter $\tau$ the streamlines change colours from yellow (bright) to green (dark). The top row shows the total motion whereas middle and bottom rows depict the curl-free and divergence-free parts of the Helmholtz decomposition, respectively. Images are arranged in accordance with Fig.~\ref{fig:of1}.}
	\label{fig:of1stream}
\end{figure}

\begin{figure}
	\centering
	\includegraphics[width=0.32\textwidth]{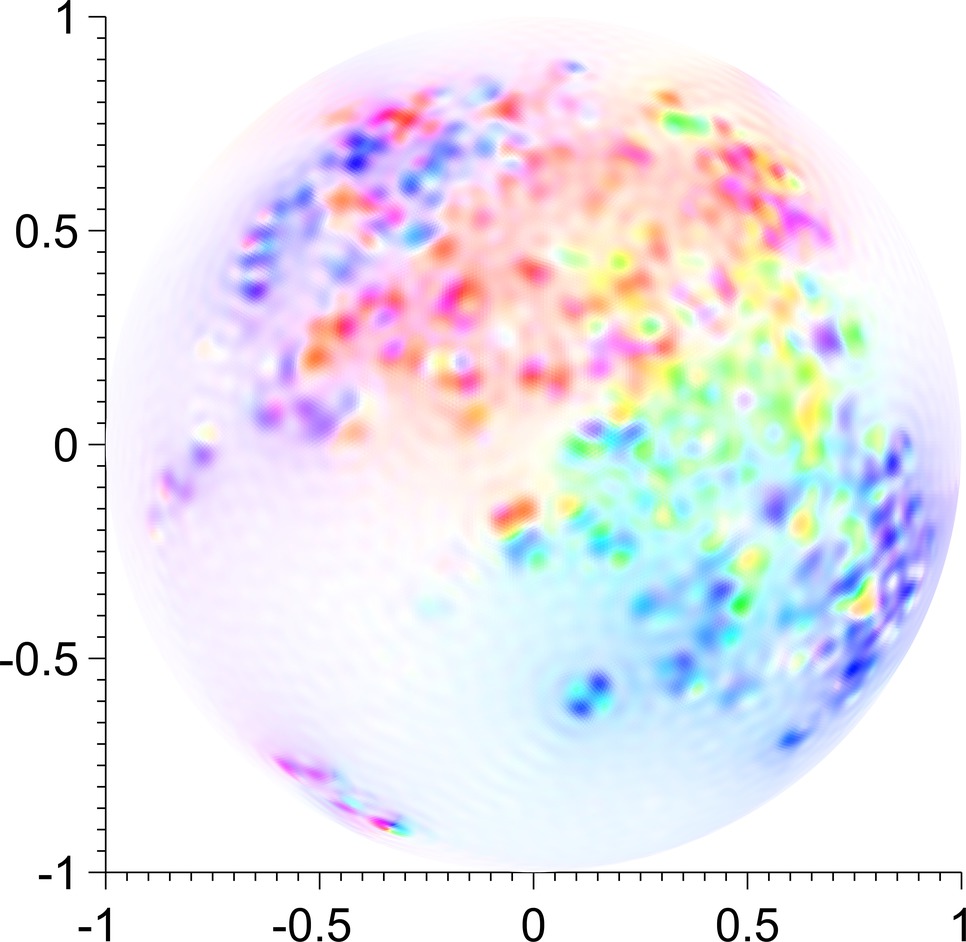} \hfill
	\includegraphics[width=0.32\textwidth]{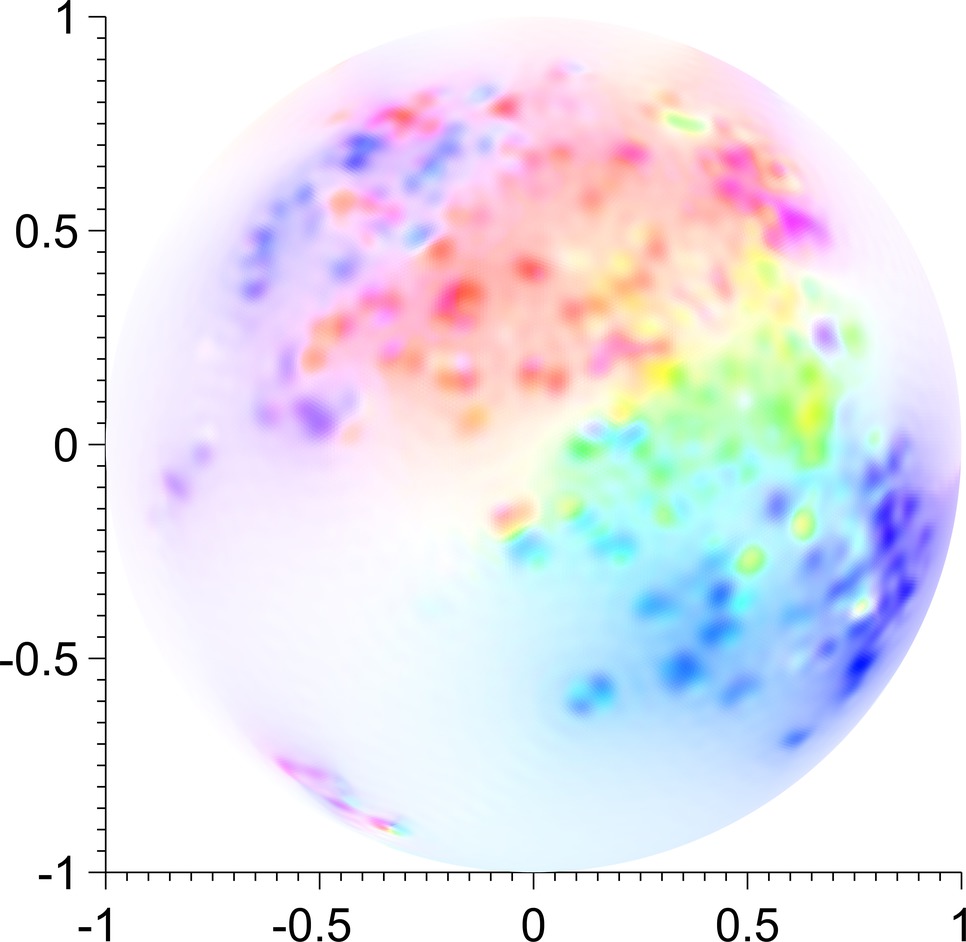} \hfill
	\includegraphics[width=0.32\textwidth]{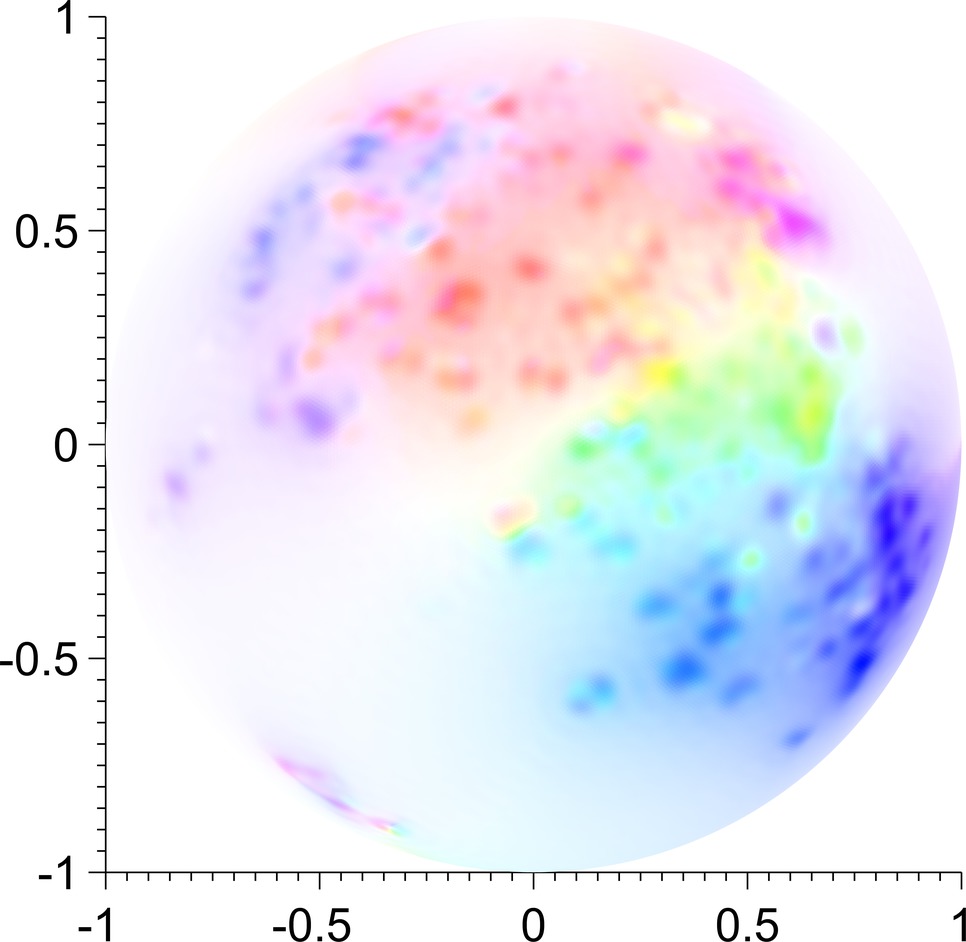} \\
	\smallskip
	\includegraphics[width=0.32\textwidth]{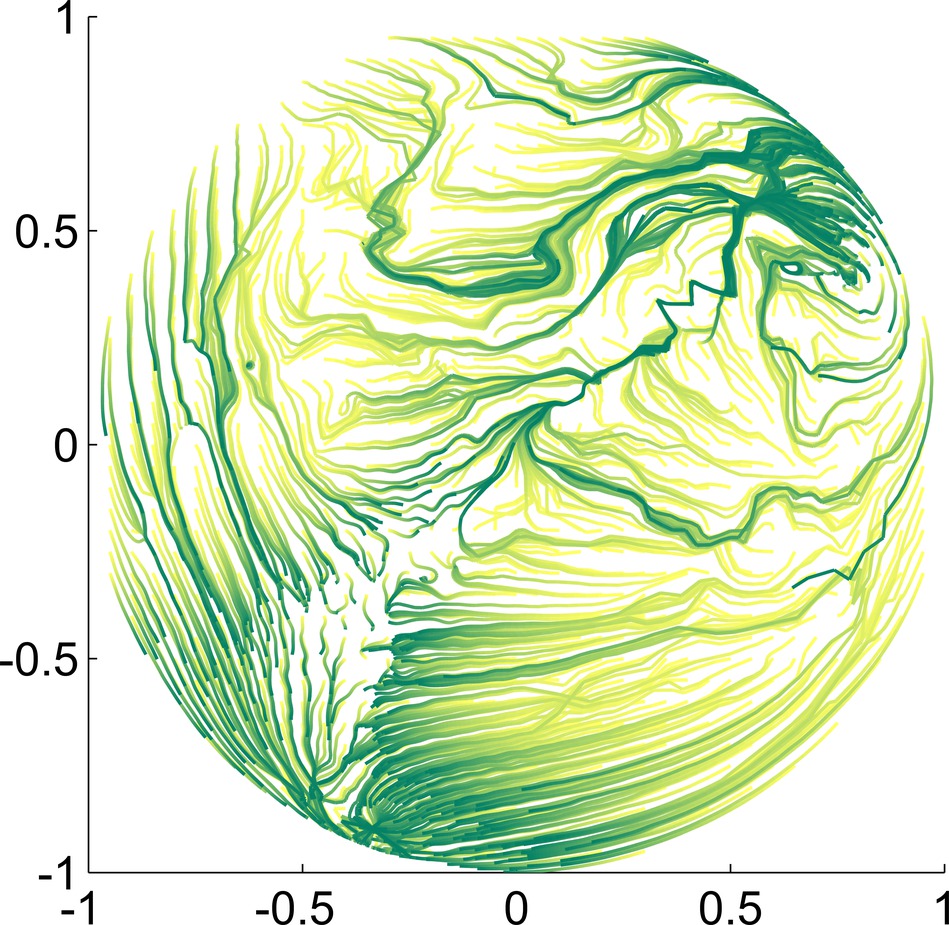} \hfill
	\includegraphics[width=0.32\textwidth]{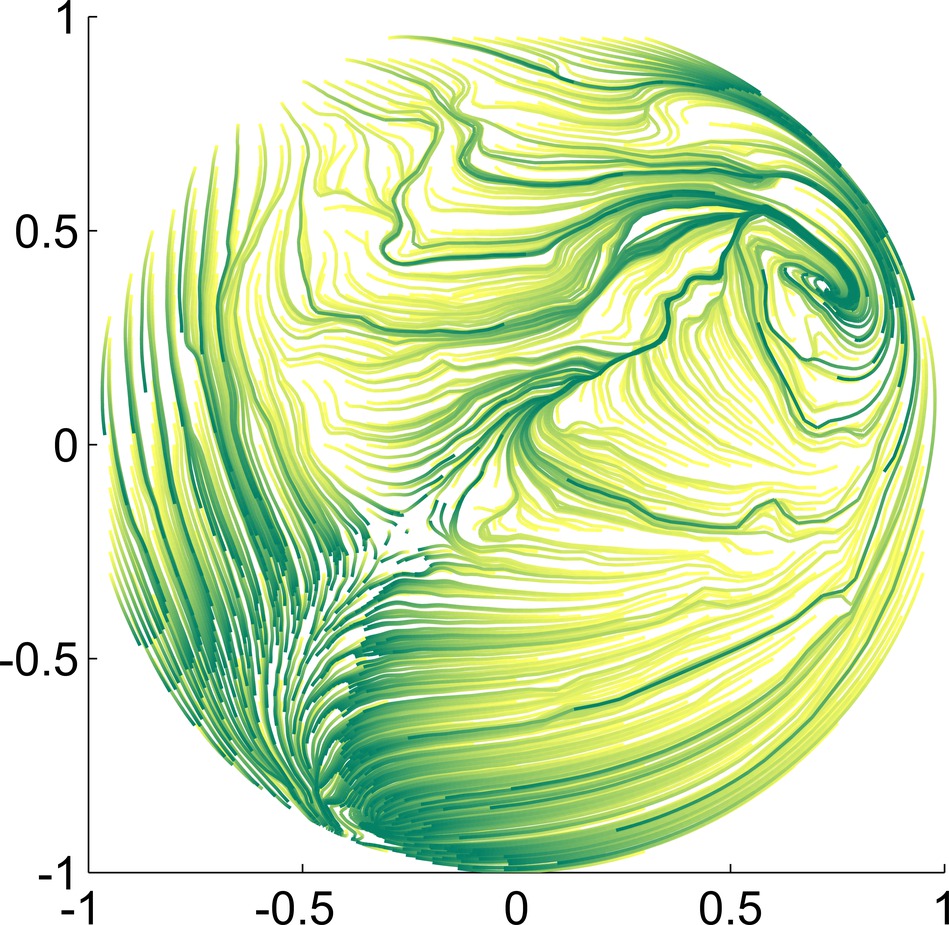} \hfill
	\includegraphics[width=0.32\textwidth]{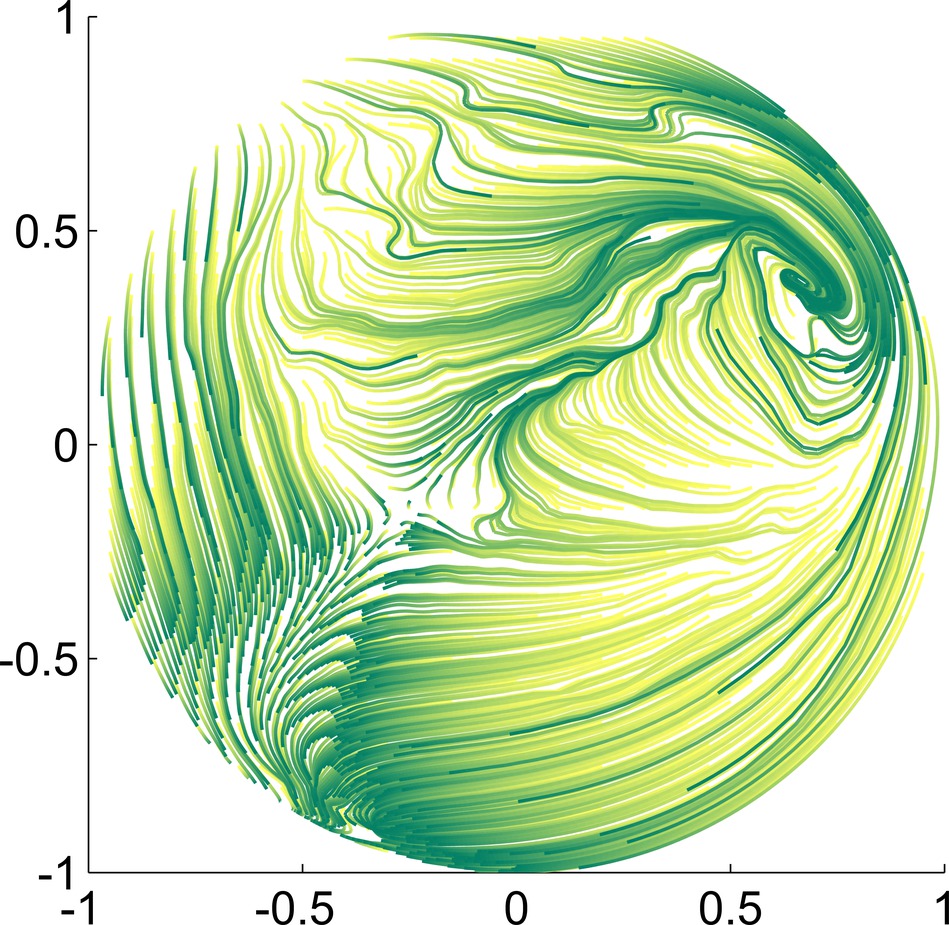}
	\caption{Colour-coded velocity fields for sequence $\mu_{n} = \alpha \lambda_{n}^{s}$ are shown in the top row. Parameters were chosen as $s = 0.5$ and $\alpha = 10$, $\alpha = 10^{2}$, and $\alpha = 10^{3}$ and are arranged from left to right. The bottom row depicts the corresponding streamlines.}
	\label{fig:of2}
\end{figure}

In the first experiment, we minimised functional $\mathcal{E}_{\mu_n}$ as defined in~\eqref{eq:rega_n} for $\mu_{n} = \alpha \lambda_{n}^{s}$ and different values of $s$ and $\alpha$. Figure~\ref{fig:of1}, top row, depicts the optical flow field for $s = 1$ and values $\alpha = 1$, $\alpha = 10$, and $\alpha = 100$. The presented results are in accordance with our findings in~\cite{KirLanSch13,KirLanSch13a_report}. As explained in Sec.~\ref{sec:decomposition}, a Helmholtz decomposition $w = w^{(2)} + w^{(3)}$ is obtained immediately. Figure~\ref{fig:of1}, middle row, shows $w^{(2)}$ whereas Fig.~\ref{fig:of1}, bottom row, shows $w^{(3)}$. Furthermore, in Fig.~\ref{fig:of1stream}, streamlines for the same velocity fields are portrayed and the individual plots are arranged accordingly. In addition, Fig.~\ref{fig:of2} shows the velocity fields for parameters $s = 0.5$, $\alpha = 10$, $\alpha = 10^{2}$, and $\alpha = 10^{3}$.

\subsubsection{$u + v$ Decomposition}

\begin{figure}
	\centering
	\includegraphics[width=0.32\textwidth]{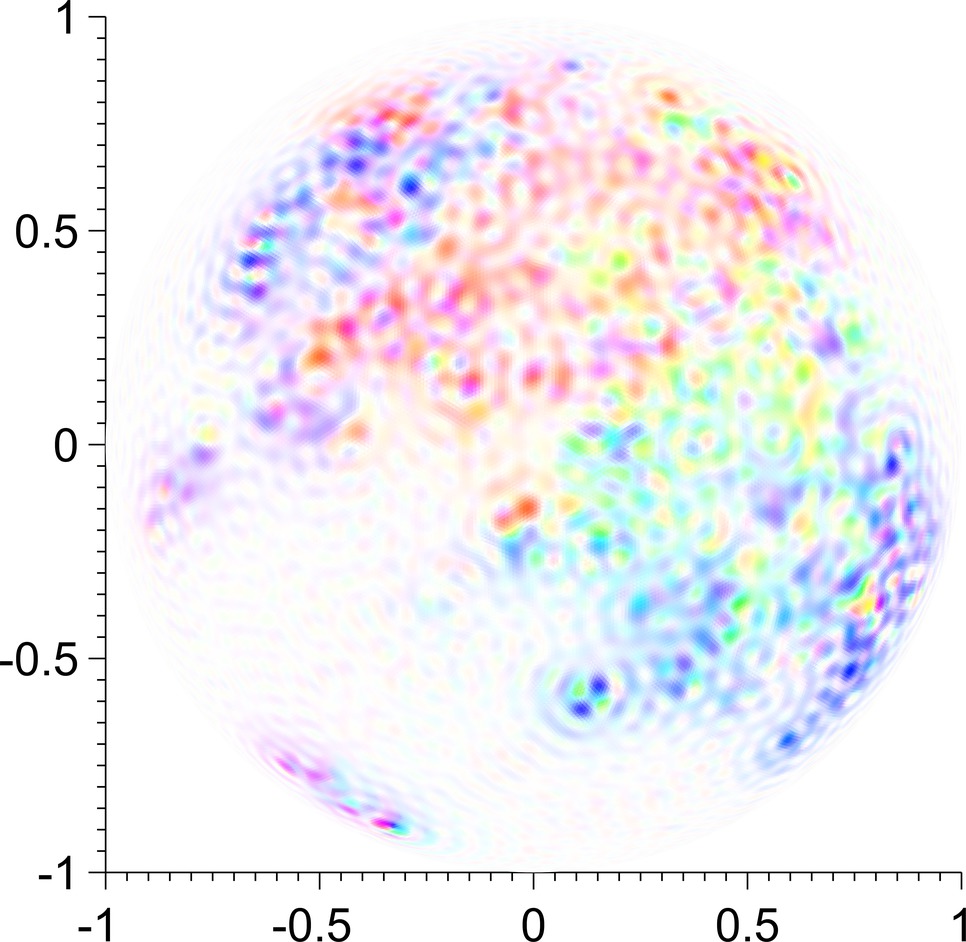} \hfill
	\includegraphics[width=0.32\textwidth]{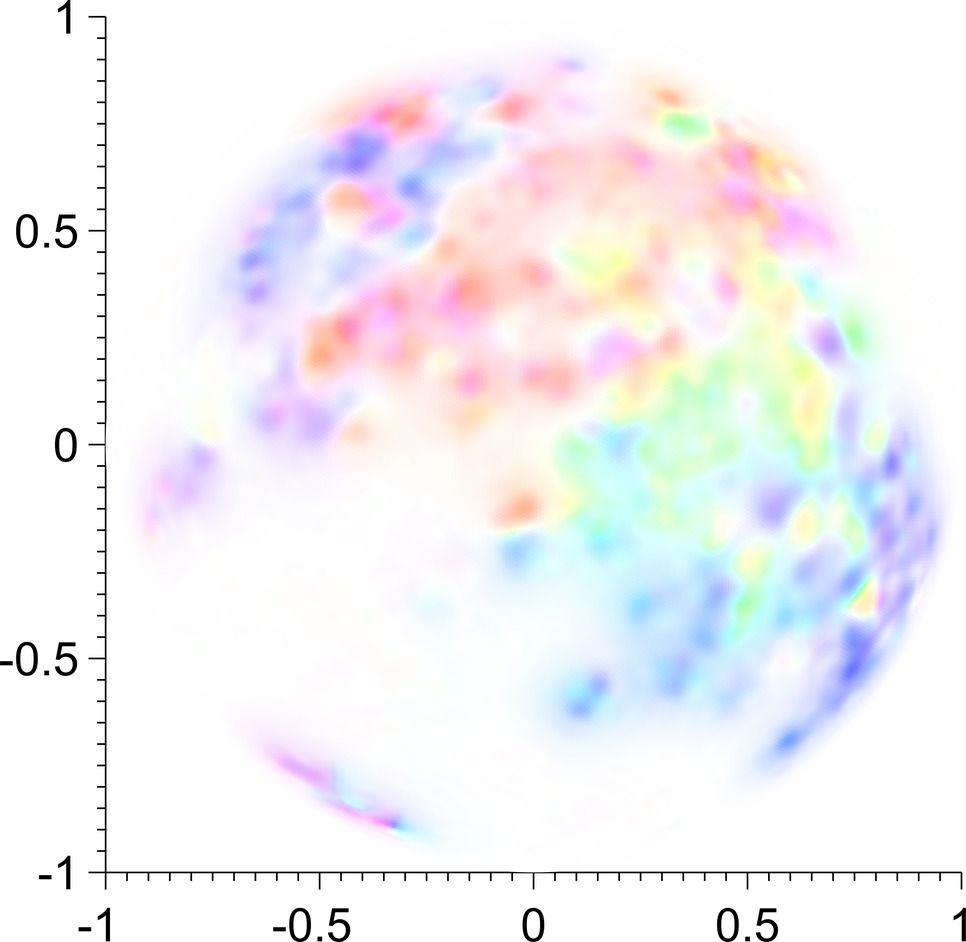} \hfill
	\includegraphics[width=0.32\textwidth]{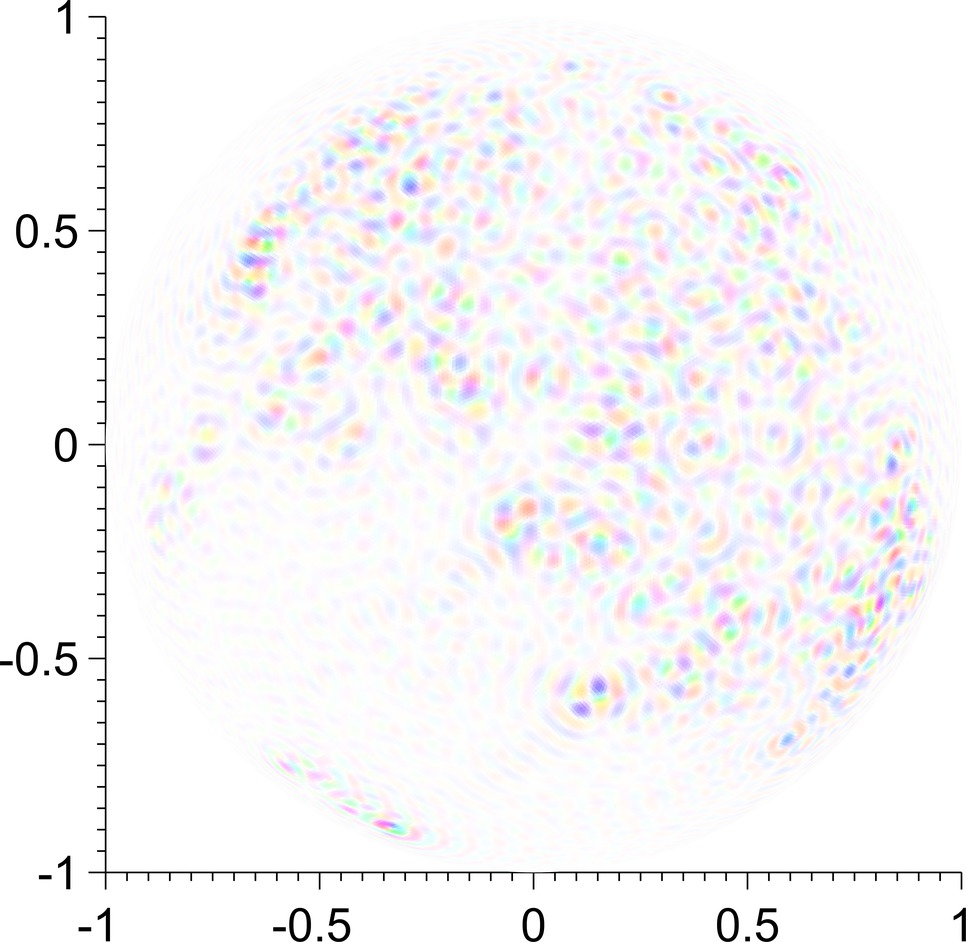}
	\caption{Decomposition of the total motion $u + v$ (left) into structural parts $u$ (middle) and $v$ (right). Sequences were chosen as $\mu_{n} = \alpha \lambda_{n}^{r}$ and $\nu_{n} = \beta \lambda_{n}^{s}$ and parameters were set to $r = 1$, $s = -1$, $\alpha = 10^{-1}$, and $\beta = 10^{6}$.}
	\label{fig:u+v}
\end{figure}

\begin{figure}
	\centering
	\includegraphics[width=0.32\textwidth]{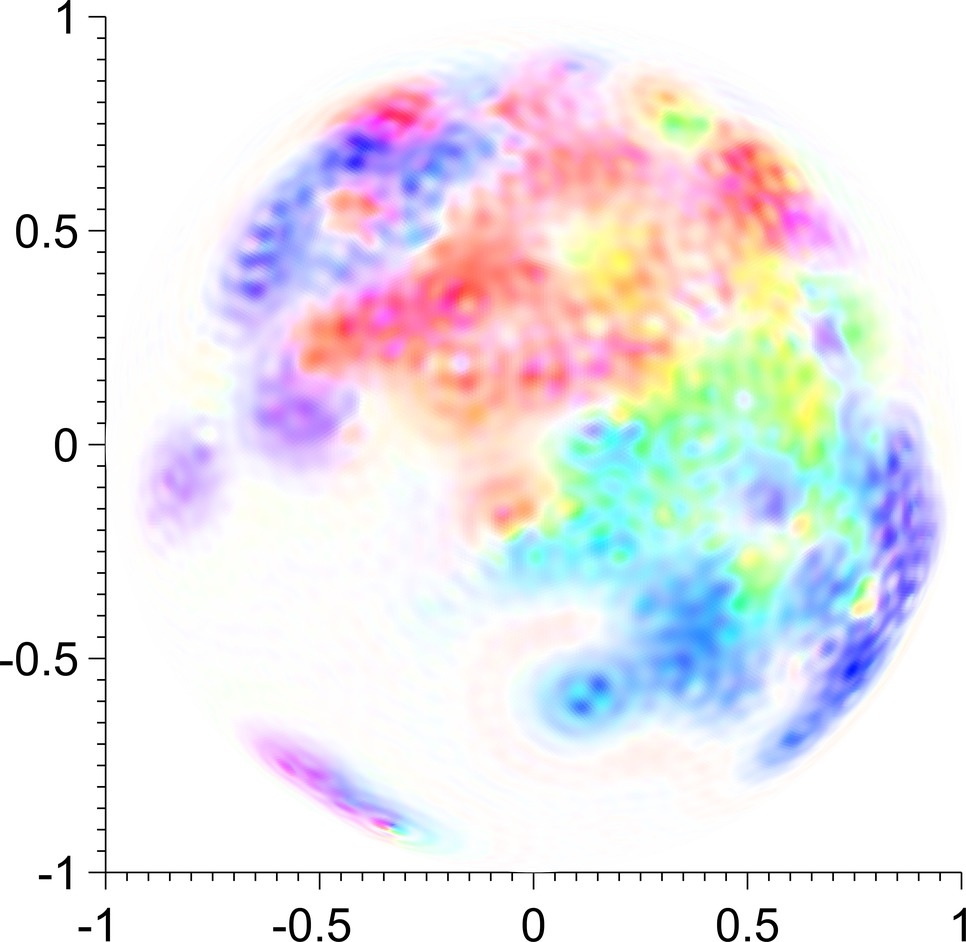} \hfill
	\includegraphics[width=0.32\textwidth]{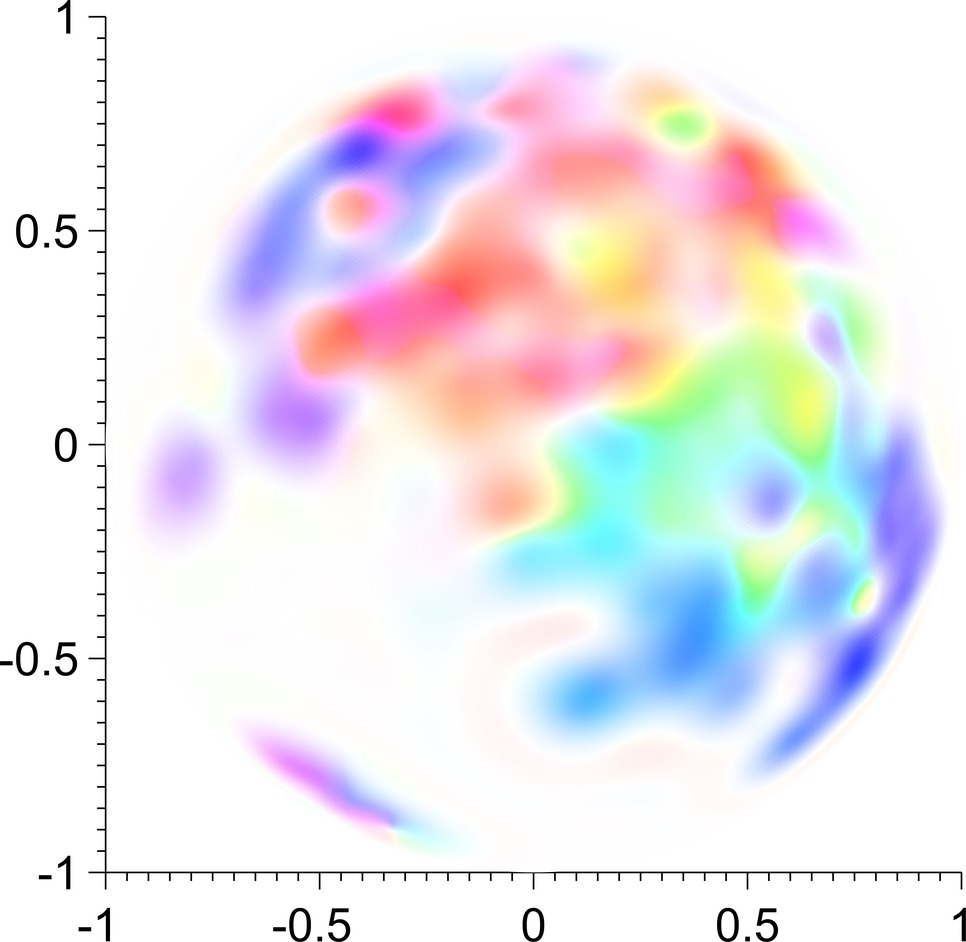} \hfill
	\includegraphics[width=0.32\textwidth]{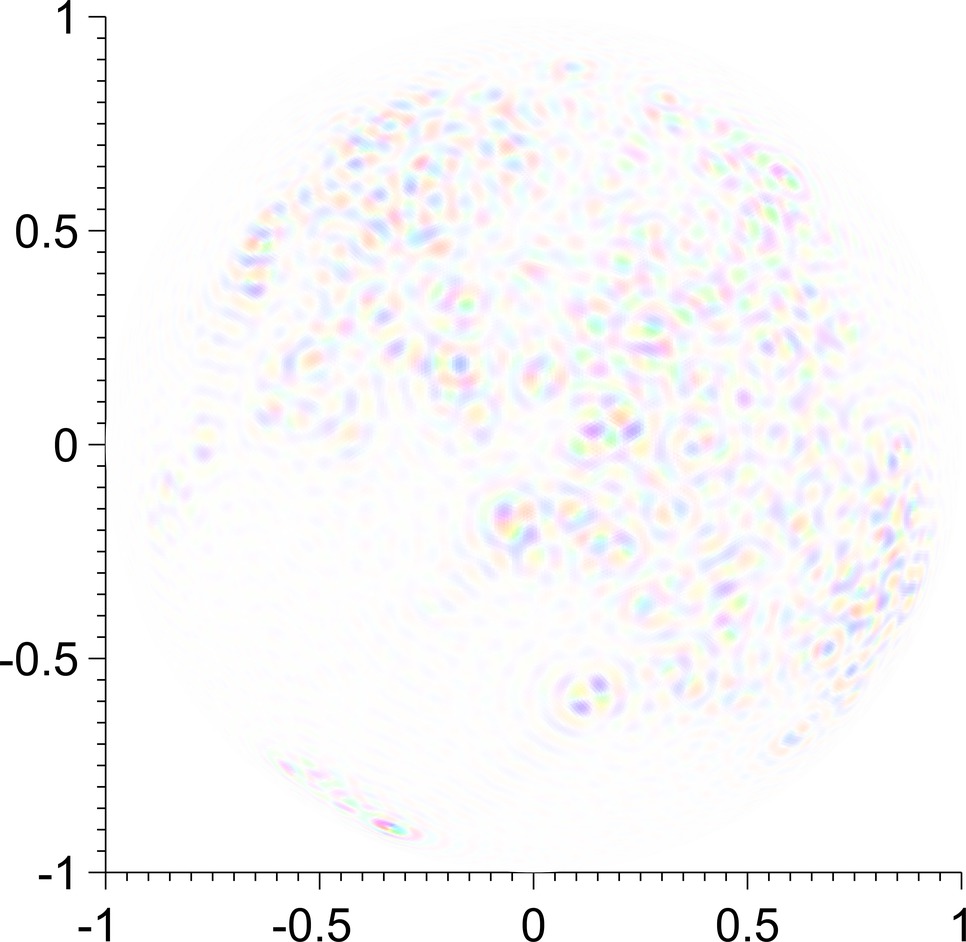}
	\caption{Decomposition of the total motion $u + v$ (left) into structural parts $u$ (middle) and $v$ (right). Sequences were chosen as $\mu_{n} = \alpha \lambda_{n}^{r}$ and $\nu_{n} = \beta \lambda_{n}^{s}$ and parameters were set to $r = 2$, $s = -1$, $\alpha = 10^{-3}$, and $\beta = 10^{7}$.}
	\label{fig:u+v2}
\end{figure}

In a next experiment, we computed a minimiser for functional~\eqref{eq:u+v} in order to obtain a $u + v$ decomposition of the optical flow. The sequences $(\mu_{n})_{n}$ and $(\nu_{n})_{n}$ were chosen as $\mu_{n} = \alpha \lambda_{n}^{r}$ and $\nu_{n} = \beta \lambda_{n}^{s}$, respectively. In Figs.~\ref{fig:u+v} and~\ref{fig:u+v2}, the resulting decomposition is shown for two different parameter settings. The motion field in Fig.~\ref{fig:u+v} was obtained with parameters $r = 1$, $s = -1$, $\alpha = 10^{-1}$, and $\beta = 10^{6}$. As anticipated, $u$ and $v$ capture different structural parts of the motion. $u$ is sufficiently smooth whereas $v$ contains spatial oscillations.

\begin{figure}
	\centering
	\frame{\includegraphics[width=0.49\textwidth]{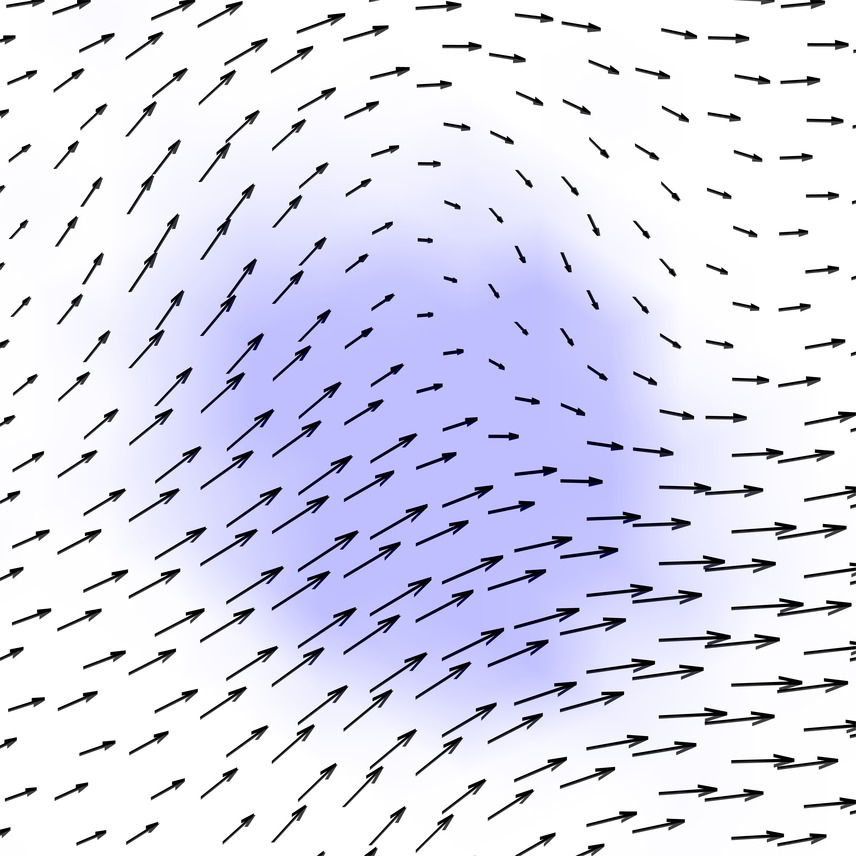}} \hfill
	\frame{\includegraphics[width=0.49\textwidth]{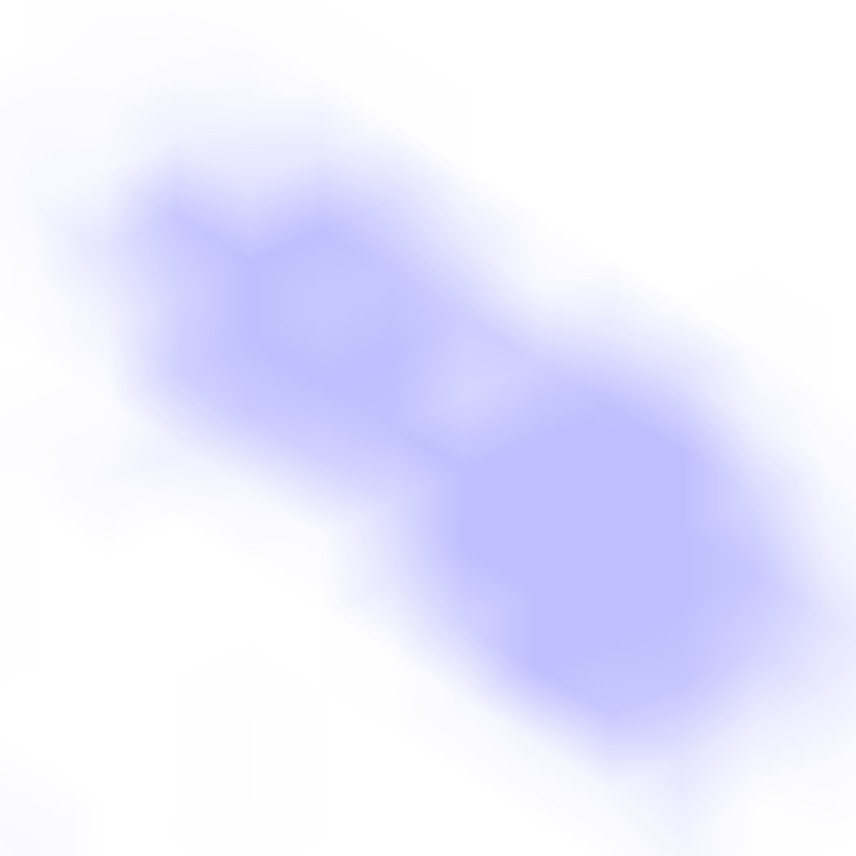}} \\ \medskip
	\frame{\includegraphics[width=0.49\textwidth]{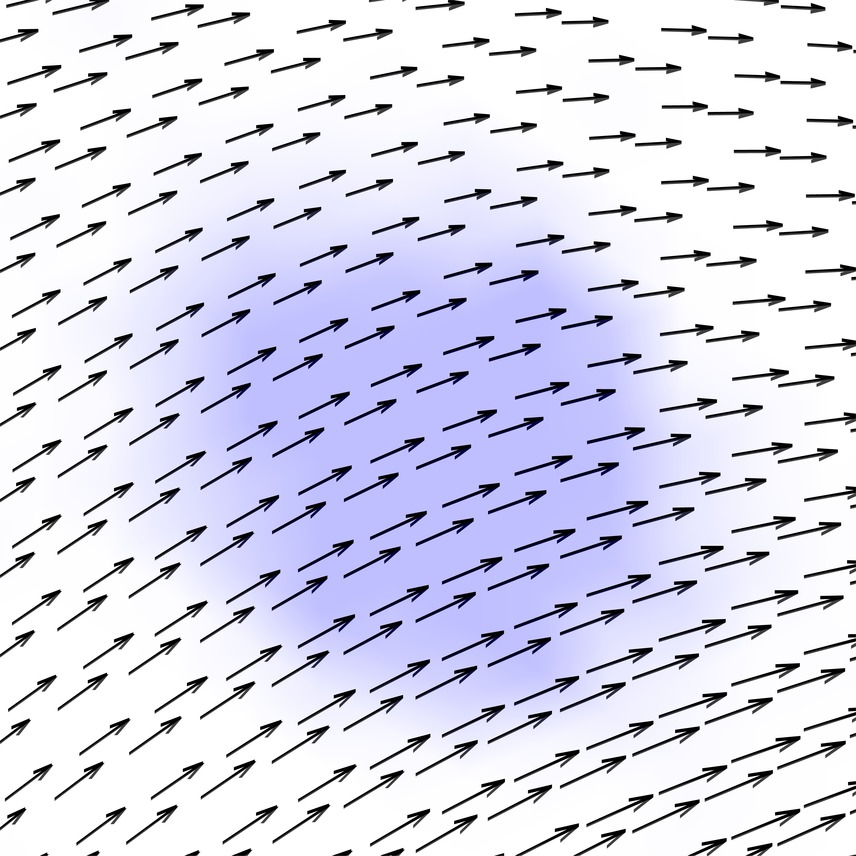}} \hfill
	\frame{\includegraphics[width=0.49\textwidth]{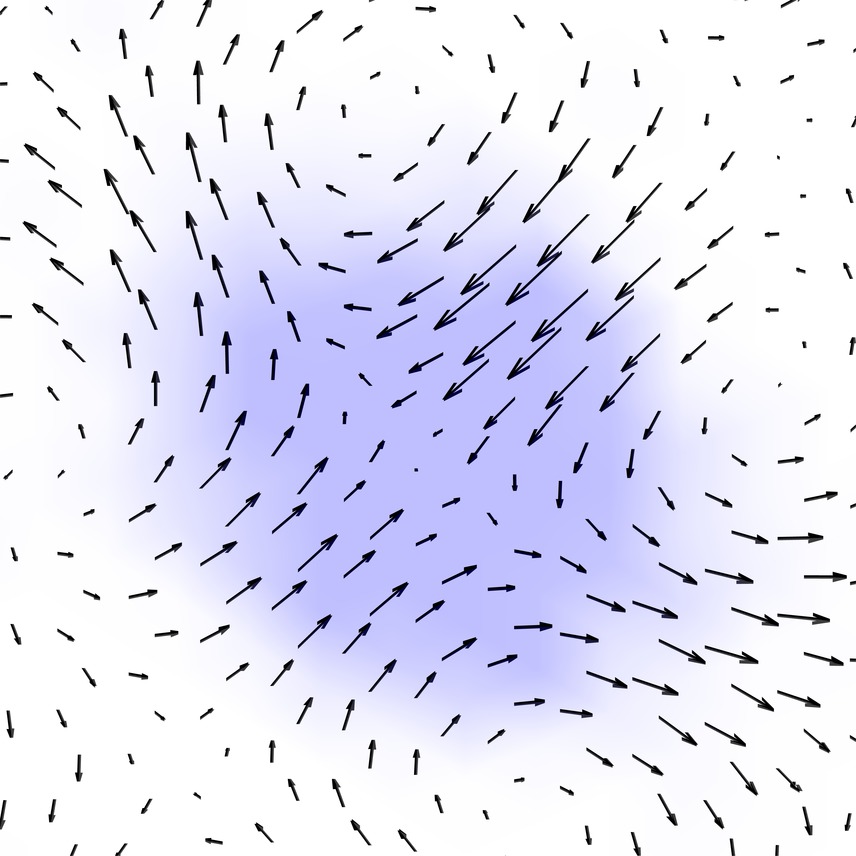}}
	\caption{Detailed view of a cell division. The same parameters as in Fig.~\ref{fig:u+v2} are used. The top left image depicts $\hat{F}_{0}$ with $u + v$ superimposed. The top right image shows $\hat{F}_{1}$. The bottom left image illustrates $u$ whereas the bottom right image shows $v$. For better illustration, $\hat{F}_{0}$ and $\hat{F}_{1}$ have been lightened up and vectors have been scaled.}
	\label{fig:detail}
\end{figure}

The result in Fig.~\ref{fig:u+v2} was computed by setting $r = 2$, $s = -1$, $\alpha = 10^{-3}$, and $\beta = 10^{7}$. Expectedly, the velocity field $u$ is, by choice of $r$, smoother than in the previous setting. In addition, Fig.~\ref{fig:detail} illustrates the characteristics of the velocity fields during a cell division in more detail. While $u$ is smooth, $v$ clearly indicates the cell division.

\subsubsection{Hierarchical Decomposition}

\begin{figure}
	\centering
	\includegraphics[width=0.32\textwidth]{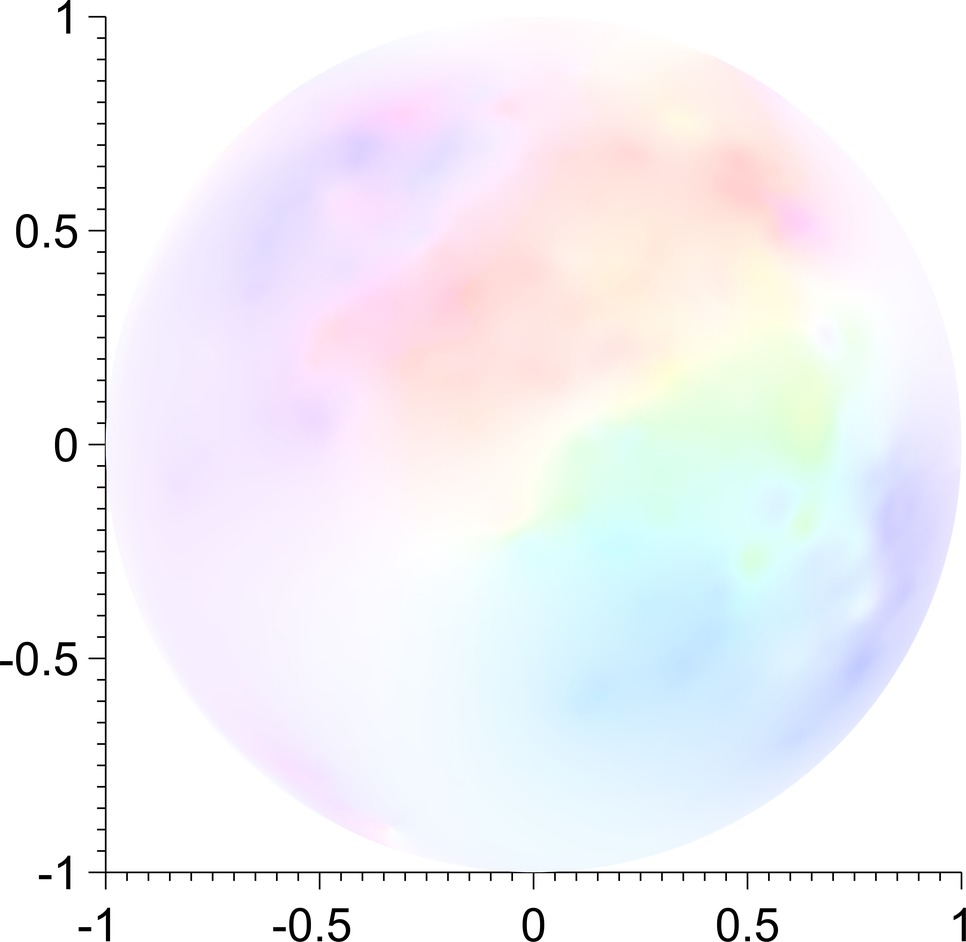} \hfill
	\includegraphics[width=0.32\textwidth]{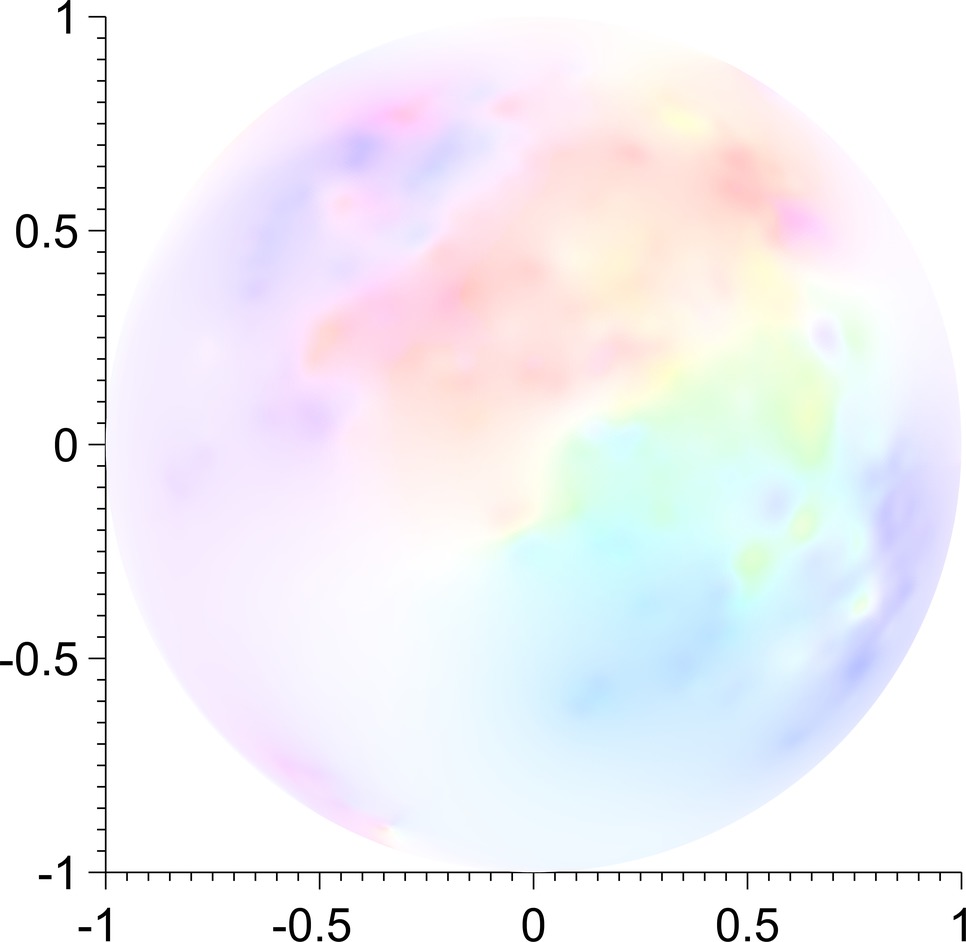} \hfill
	\includegraphics[width=0.32\textwidth]{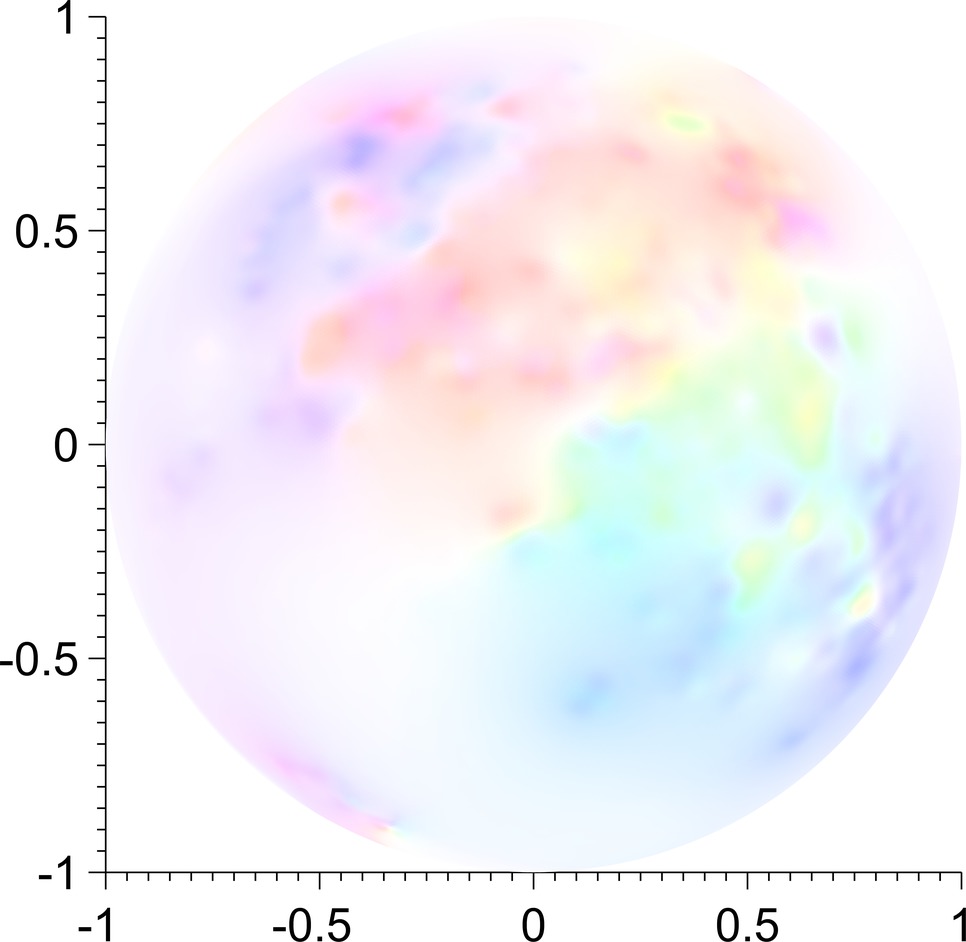} \\
	\smallskip
	\includegraphics[width=0.32\textwidth]{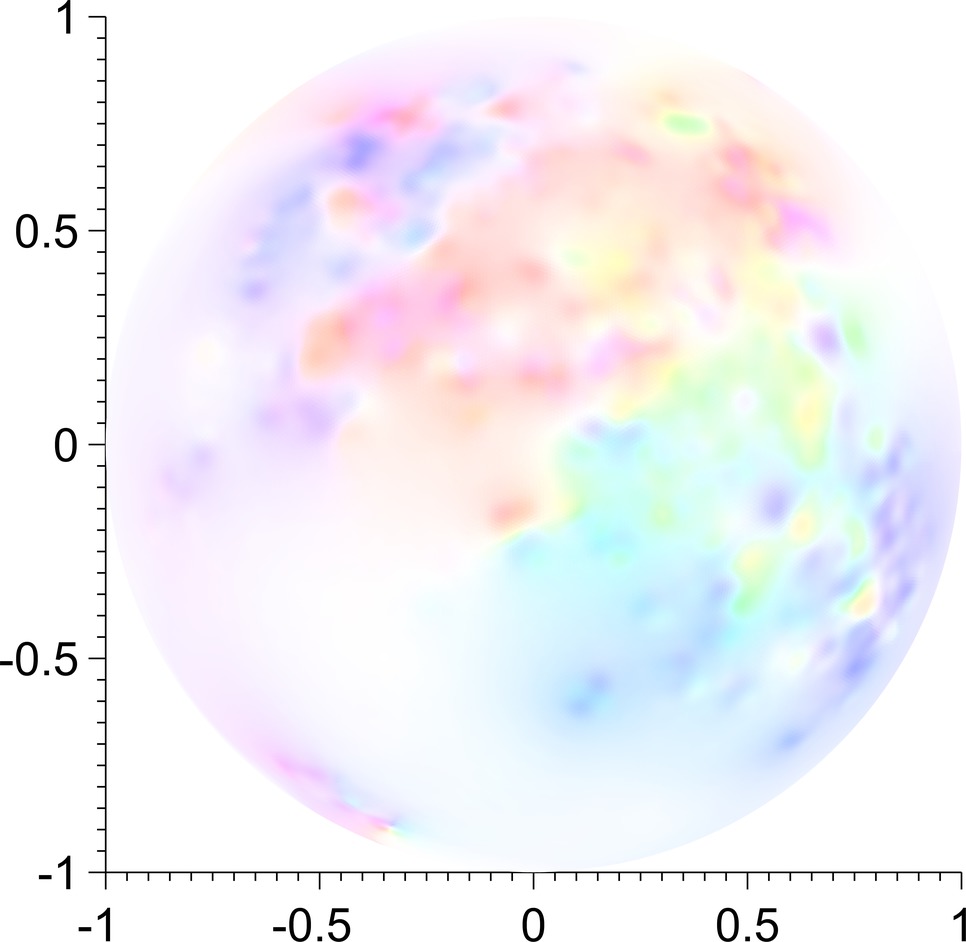} \hfill
	\includegraphics[width=0.32\textwidth]{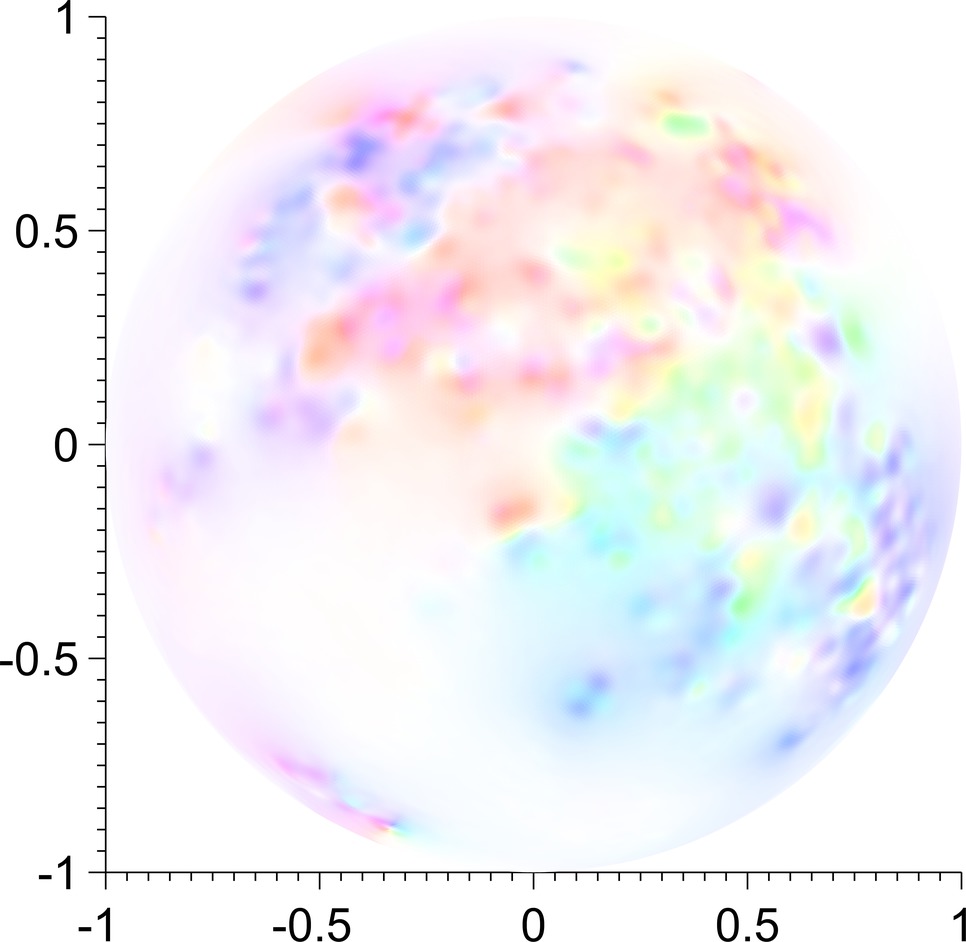} \hfill
	\includegraphics[width=0.32\textwidth]{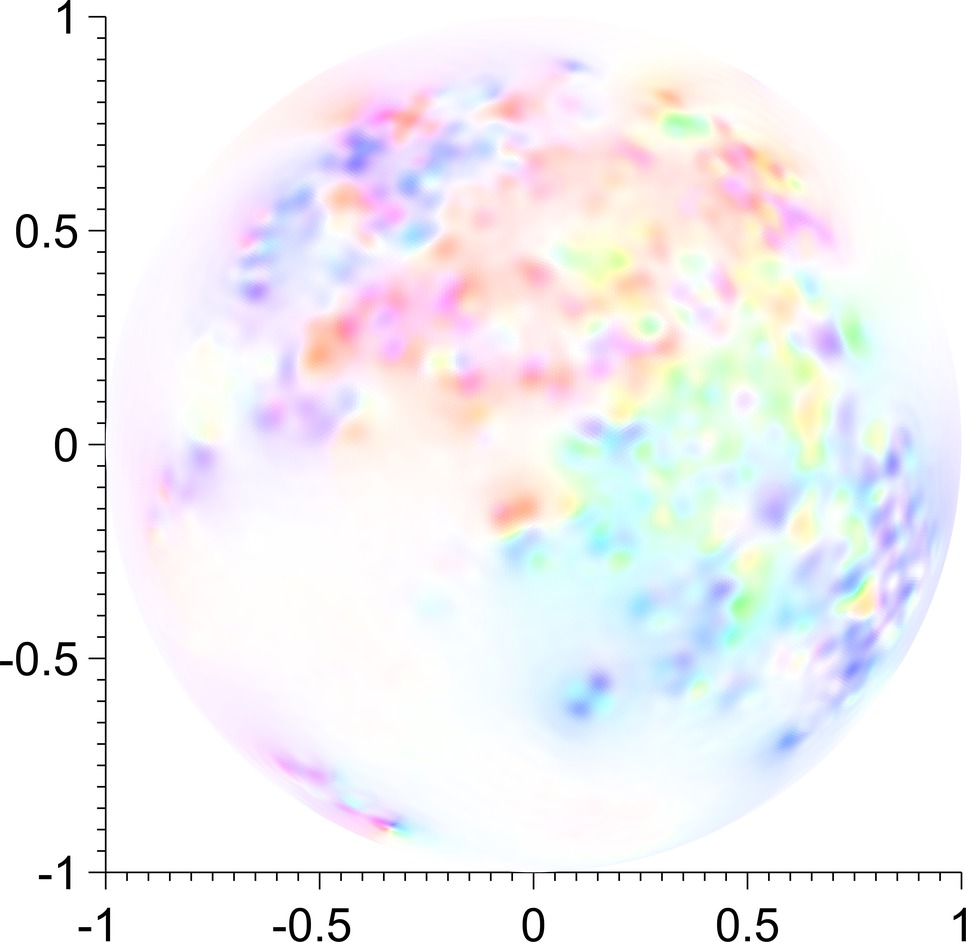} \\
	\smallskip
	\includegraphics[width=0.32\textwidth]{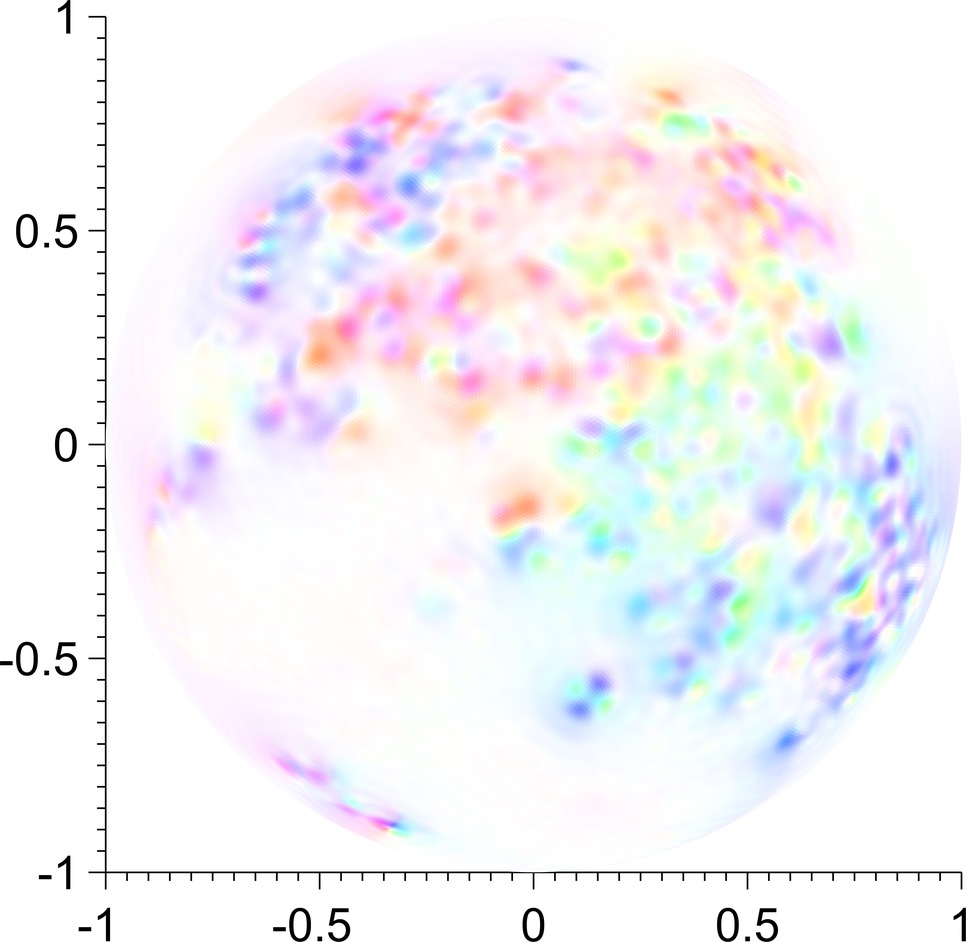} \hfill
	\includegraphics[width=0.32\textwidth]{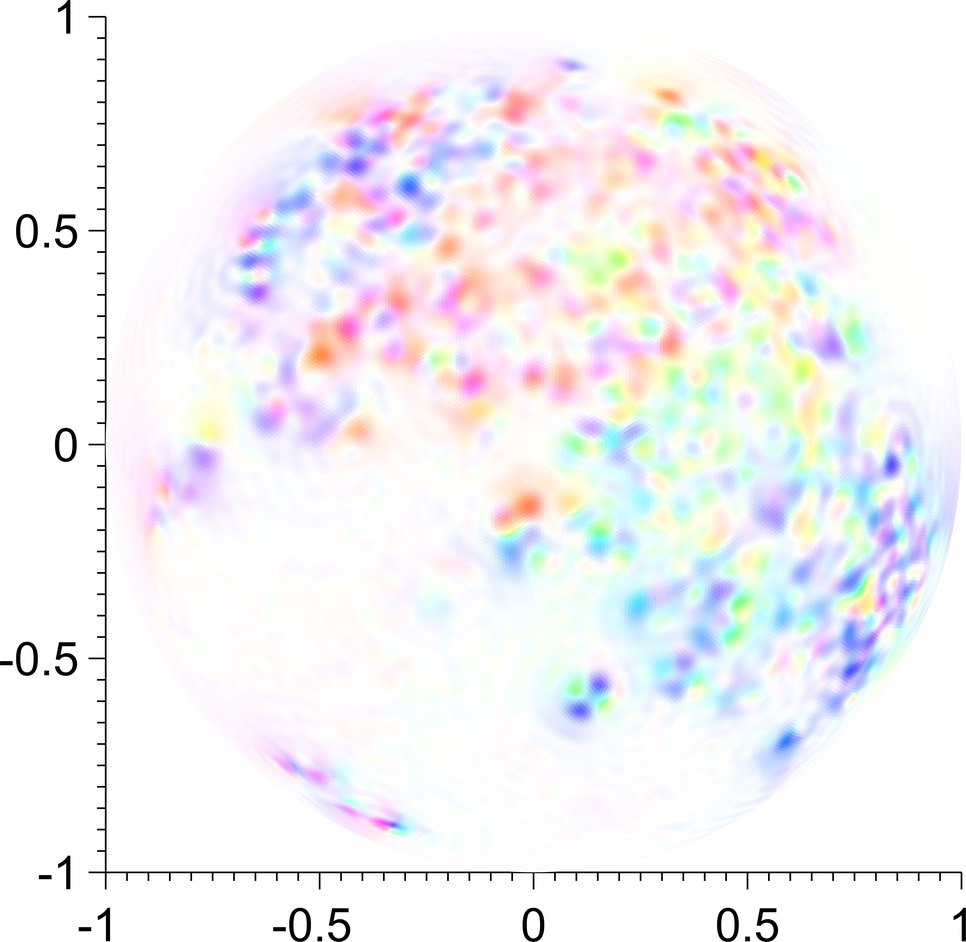} \hfill
	\includegraphics[width=0.32\textwidth]{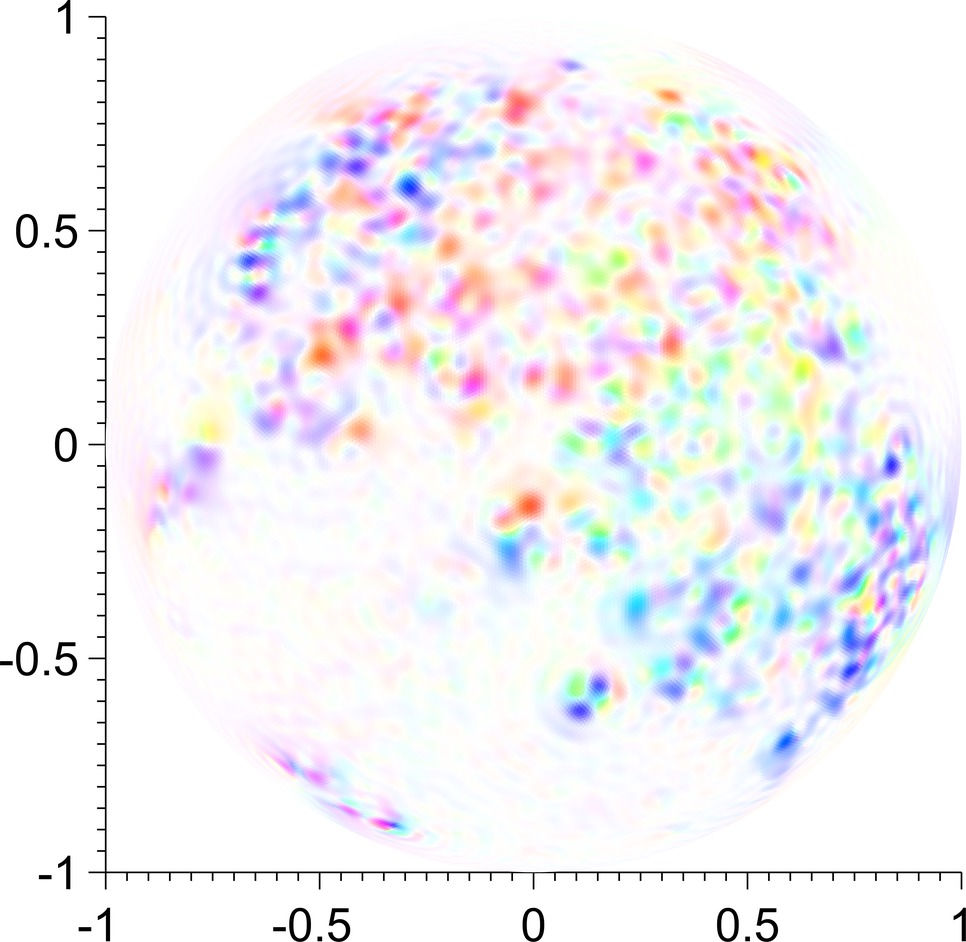}
	\caption{Solutions $(u^{(k)})_{k = 8,\dots,16}$ of velocity fields obtained by the hierarchical decomposition. Images are aligned from left to right and top to bottom. At iteration $k$ the sequence was set to $\mu_{n}^{(k)} = 2^{1-k}\alpha \lambda_{n}^{s}$ with parameters $s = 1$ and $\alpha = 1000$.}
	\label{fig:hierarchical}
\end{figure}

\begin{figure}
	\centering
	\includegraphics[width=0.32\textwidth]{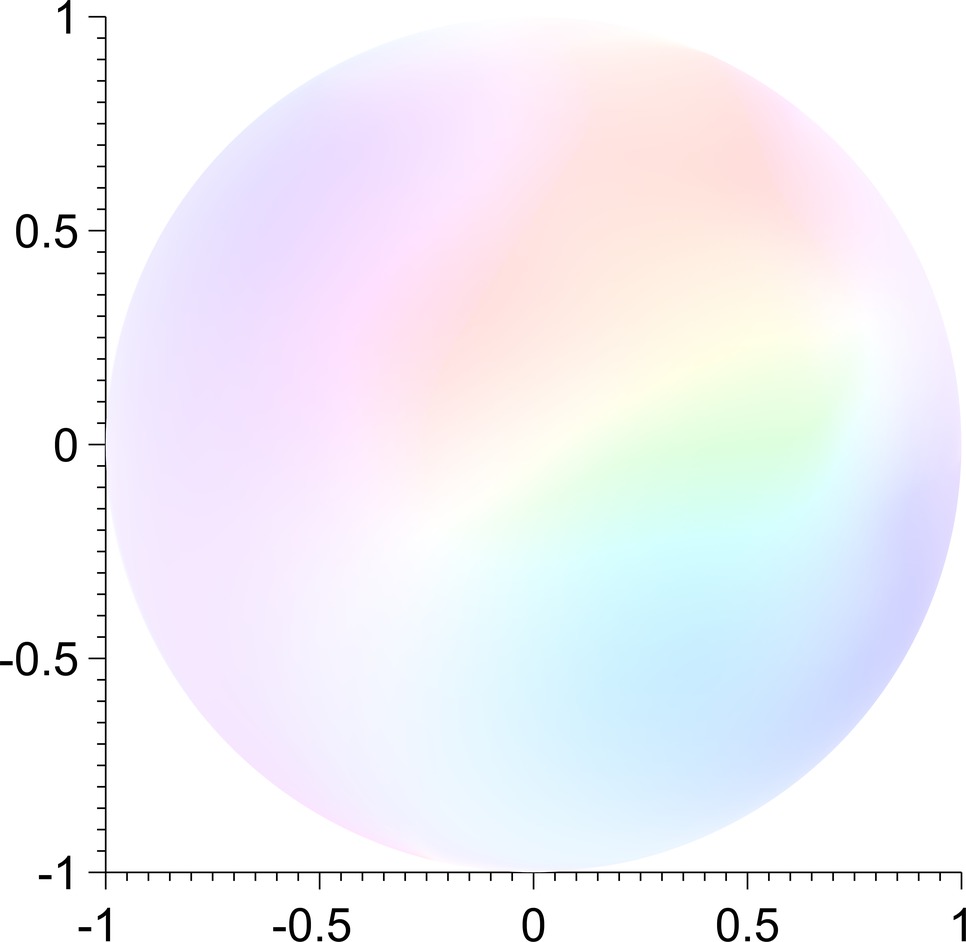} \hfill
	\includegraphics[width=0.32\textwidth]{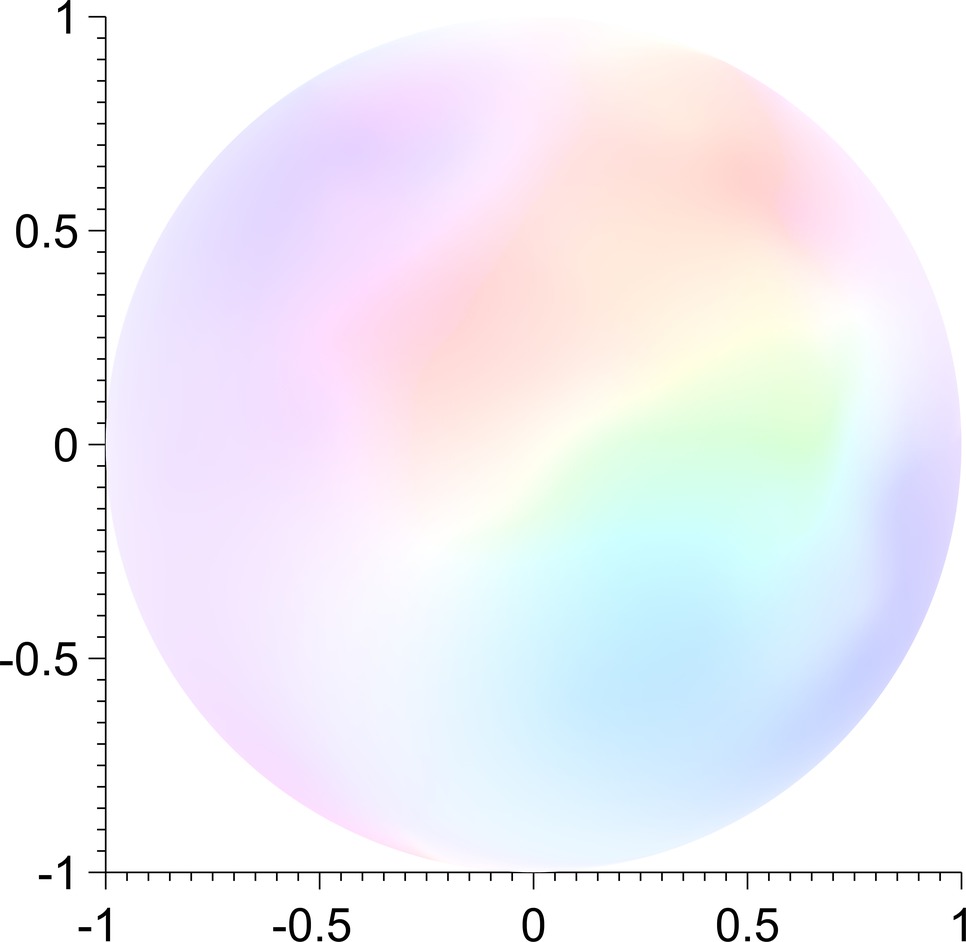} \hfill
	\includegraphics[width=0.32\textwidth]{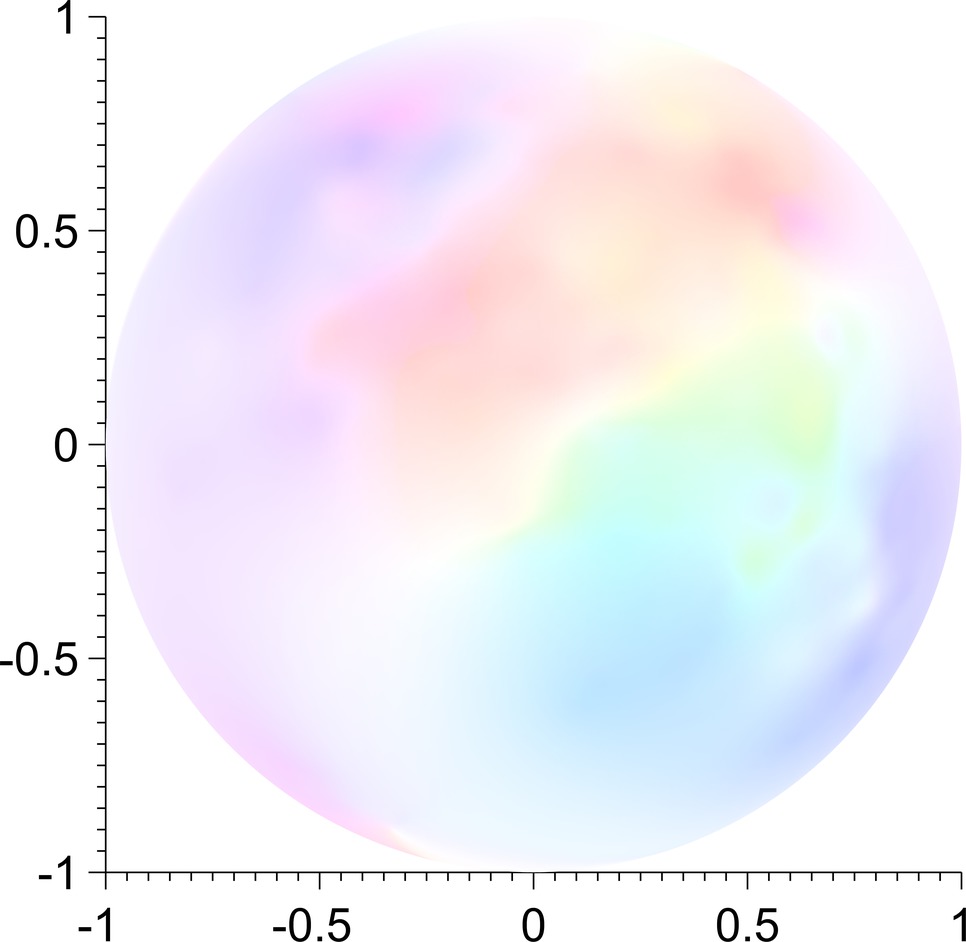} \\
	\smallskip
	\includegraphics[width=0.32\textwidth]{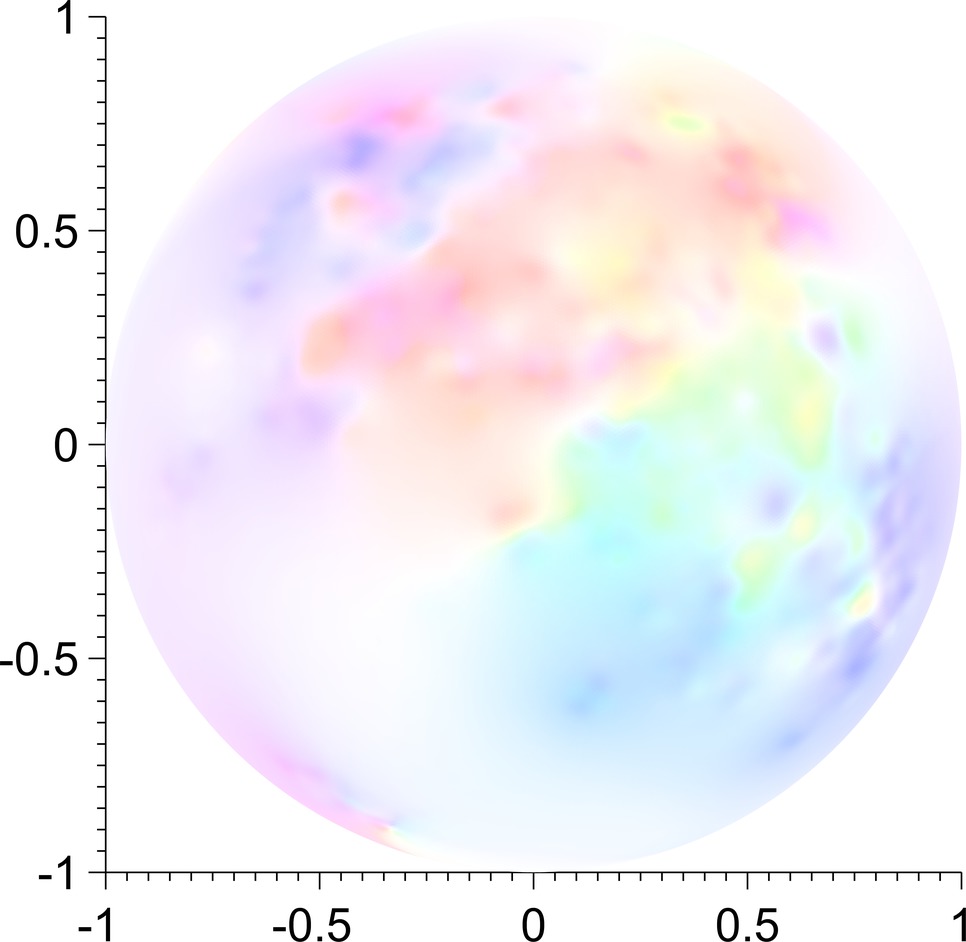} \hfill
	\includegraphics[width=0.32\textwidth]{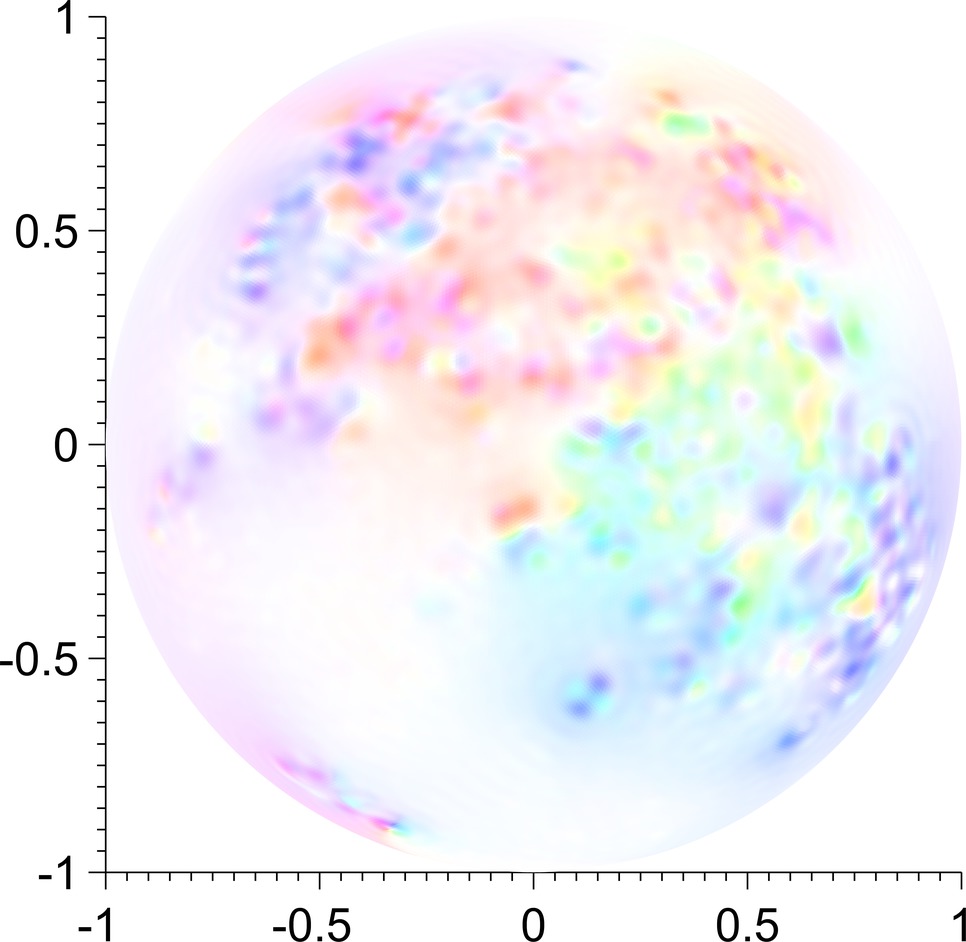} \hfill
	\includegraphics[width=0.32\textwidth]{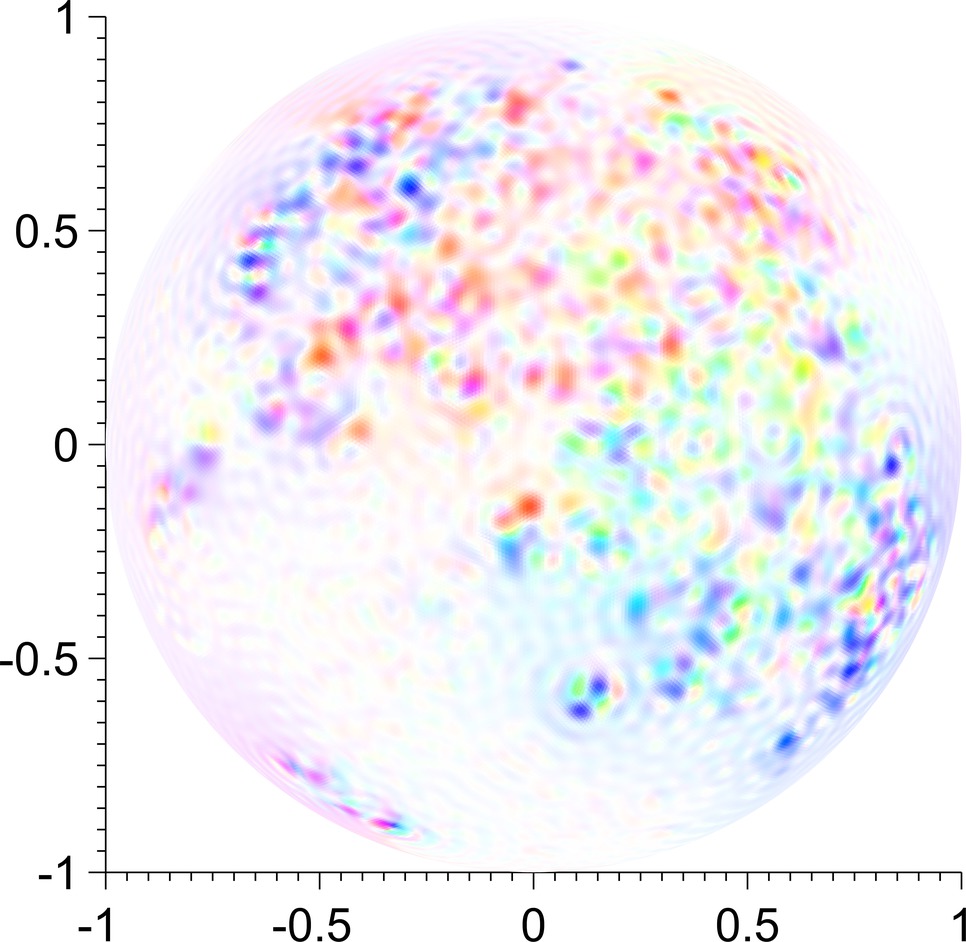}
	\caption{Solutions $(u^{(k)})_{k = 2,\dots,7}$ of velocity fields obtained by the hierarchical decomposition. Images are aligned from left to right and top to bottom. At iteration $k$ the sequence was set to $\mu_{n}^{(k)} = \alpha \lambda_{n}^{s - \frac{k - 1}{4}}$ with parameters $s = 2$ and $\alpha = 1$.}
	\label{fig:hierarchical2}
\end{figure}

As a final experiment, we computed two types of hierarchical decompositions as proposed in Sec.~\ref{sec:decomposition}. First we chose	$\mu_{n}^{(k)} = 2^{1-k}\alpha \lambda_{n}^{s}$ such that $\alpha$ is halved in every iteration. Solutions $(u^{(k)})_{k = 1,\dots,16}$ were obtained using parameters $s = 1$ and $\alpha = 1000$. In Fig.~\ref{fig:hierarchical}, the subsequence $(u^{(k)})_{k = 8,\dots,16}$ is shown. As $k$ increases, the motion field expands on the details. In a second run, $\mu_{n}^{(k)}$ was set to $\alpha \lambda_{n}^{s - \frac{k - 1}{4}}$, decreasing the exponent of $\lambda_{n}$ by $0.25$ in every step. We iteratively computed solutions $(u^{(k)})_{k = 1,\dots,9}$ with parameters $s = 2$ and $\alpha = 1$, which were kept constant this time. Figure~\ref{fig:hierarchical2} depicts the subsequence $(u^{(k)})_{k = 2,\dots,7}$. Note that in Figs.~\ref{fig:hierarchical} and~\ref{fig:hierarchical2} the colour-coded visualisation is relative to the chosen subsequence. As an exception, here we allowed a maximum number of $1000$ iterations instead of $100$ for the linear system solve to ensure a relative residual of 0.025 in the first step of the hierarchical decomposition.

\section{Conclusion}
We provided a set of variational methods for the analysis of motion fields. While their applicability is limited to data given on the sphere, the proposed models have great flexibility in terms of possible regularising functionals. In fact, the chosen numerical method based on tangential vector spherical harmonics allows for a straightforward usage of Sobolev $H^s$ norms for any real $s$. Combined with both $u+v$ and hierarchical decomposition models, which we adapted to the spherical optical flow setting, this flexibility makes it possible to capture different motion characteristics with ease. Feasibility of the proposed models was verified on a microscopy dataset depicting endodermal cells of a zebrafish embryo.

\paragraph{Acknowledgements.}
We thank Pia Aanstad from the University of Innsbruck for sharing her biological insight and for kindly providing the microscopy data. This work has been supported by the Vienna Graduate School in Computational Science (IK I059-N) funded by the University of Vienna. In addition, we acknowledge the support by the Austrian Science Fund (FWF) within the national research networks ``Photoacoustic Imaging in Biology and Medicine" (project S10505-N20, Reconstruction Algorithms for PAI) and ``Geometry + Simulation" (project S11704, Variational Methods for Imaging on Manifolds).

\def\cprime{$'$} \providecommand{\noopsort}[1]{}\def\cprime{$'$}

\end{document}